\documentclass[5p,times]{elsarticle}

\usepackage{subcaption}
\usepackage{amsmath}
\usepackage{hyperref}
\usepackage{xcolor}
\usepackage[format=plain, labelfont={bf}]{caption}
\usepackage{bm}

\begin{document}
\begin{frontmatter}

\title{L1-norm vs.\ L2-norm fitting in optimizing focal multi-channel tES stimulation: linear and semidefinite programming vs.\ weighted least squares}

\author[1]{Fernando Galaz Prieto\corref{cor1}}
\ead{fernando.galazprieto@tuni.fi}
\cortext[cor1]{Corresponding author at: Sähkötalo building, Korkeakoulunkatu 3, Tampere, 33720, FI}

\author[1]{Atena Rezaei}
\author[1]{Maryam Samavaki}
\author[1]{Sampsa Pursiainen}

\address[1]{Computing Sciences Unit, Faculty of Information Technology and Communication Sciences, Tampere University, Tampere, Finland}

\begin{abstract}
\textbf{Background and Objective:} This study focuses on Multi-Channel Transcranial Electrical Stimulation, a non-invasive brain method for stimulating neuronal activity under the influence of low-intensity currents. We introduce mathematical formulation for finding a current pattern which optimizes a L1-norm fit between a given focal target distribution and volume current density inside the brain. L1-norm is well-known to favor well-localized or sparse distributions compared to L2-norm (least-squares) fitted estimates.

\noindent \textbf{Methods:} We present a linear programming approach which performs L1-norm fitting and penalization of the current pattern (L1L1) to control the number of non-zero currents. The optimizer filters a large set of candidate solutions using a two-stage metaheuristic search in from a pre-filtered set of candidates.

\noindent \textbf{Results:} The numerical simulation results, obtained with both a 8- and 20-channel electrode montages, suggest that our hypothesis on the benefits of L1-norm data fitting is valid. As compared to L1-norm regularized L2-norm fitting (L1L2) via semidefinite programming and weighted Tikhonov least-squares method, the L1L1 results were overall preferable with respect to maximizing the focused current density at the target position and the ratio between focused and nuisance current magnitudes.

\noindent \textbf{Conclusions:} We propose the metaheuristic L1L1 optimization approach as a potential technique to obtain a well-localized stimulus with a controllable magnitude at a given target position. L1L1 finds a current pattern with a steep contrast between the anodal and cathodal electrodes meanwhile suppressing the nuisance currents in the brain, hence, providing a potential alternative to modulate the effects of the stimulation, e.g., the sensation experienced by the subject.

\end{abstract}

\end{frontmatter}

\section{Introduction}
\label{sec:intro}
In this numerical simulation study, we consider the task to optimize stimulation currents in the multi-channel version of transcranial electrical stimulation (tES) \cite{herrmann2013transcranial, fernandez2020unification} which is applied non-invasively for stimulating neuronal activity, treating psychiatric disorders and studying neuronal behavior. In tES, a pattern of electric currents are applied through a set of electrodes attached to subject's head. Part of the generated diffusive current field penetrates through the skull into the brain modulating cortical excitability \cite{Nitsche_Paulus2000}. The procedure for adjusting the electrode montage delivering the stimulus varies from one method to another, considering various properties such as the number of active electrodes, physical description (e.g., positioning, shape, permittivity and impedance values), applied stimulus waveform (e.g., amplitude, pulse shape, pulse width, and polarity), the number of stimulation sessions, and time interval in-between \cite{Peterchev2012}.

Transcranial direct current stimulation (tDCS) delivers constant, low-intensity currents injections of 0.5-4.0 mA \cite{Zaghi2010, fregni2005_TDCS, KHADKA202069} over a pair of large saline-soaked 20-35 cm$^2$ electrode patches \cite{MORENODUARTE201435}, with one patch adhered on the scalp, whereas the second patch can be either cephalic or extra-cephalic \cite{Datta2011, neuling2012finite}. The drawback, however, is the limitation of delivering target specific frequencies and the lack of focality. tDCS is a well-known treatment of neuropsychiatric disorders and brain illnesses, for instance, stroke conditions \cite{fregni2005_TDCS, Lindenberg2176}, epilepsy syndromes \cite{fregni2006}, Parkinson's disease \cite{boggio2006, Benninger1105, fregni2005_parkinson, fregni2006_parkinson}, major depression disorder \cite{boggio2007, fregni2006_depre, NITSCHE200914}, tinnitus \cite{fregni20066_tinnitus}, migraine \cite{Antal2003}, and alcoholism \cite{BOGGIO2008}.

Multi-channel tES generally constitutes a task to select a multi-component current pattern to create a sought field in a given location. Selecting such pattern poses an ill-posed inverse problem \cite{kaipio2006statistical,bertero2020introduction}, i.e., it does not have a unique solution and a slight change in the selected current pattern can significantly change the resulting current density in the brain. The problem can be considered as over-determined, i.e., the three-dimensional current field inside the brain is likely to have more degrees of freedom than the current pattern. The total dose of the current pattern must be limited to a given value, which typically reaches up to 2 mA \cite{nitsche2003safety} or 4 mA \cite{KHADKA202069}. Advanced optimization solutions can be obtained via regularized data fitting methods, or projection approaches which target to maximize the current in a given location \cite{dmochowski2011optimized, dmochowski2017optimal, wagner2016optimization, fernandez2020unification, Daubechies2010, khanindividually}. This study aims at finding current patterns that would optimize the L1-norm fit between a given focal (well-localized) target current distribution and, at the same time, minimize the dose given to the subject of stimulation. To achieve this, we introduce an L1-norm fitted and regularized linear programming approach (L1L1) method for finding a focal current distribution since L1-norm based solutions are generally known to be well-localized compared to regularized least squares estimation in inverse modelling \cite{kaipio2006statistical, bertero2020introduction}. While L1-norm has been previously applied to penalize an objective function \cite{wagner2016optimization} and linear programming as a strategy to maximize current density at a given location \cite{dmochowski2011optimized}, our method is one of the first to optimize a global L1-norm fit on this application. We hypothesize that the our method can be advantageous for a configuration where a focal current distribution is sought using a given number of electrodes to deliver the stimulus.

Solutions from our L1L1 algorithm are sensitive to parameter selection, that is, a wide range of regularization parameters and optimization tolerances are required to be covered to allow the search algorithm to effectively obtain a solution. Therefore, we apply a two-stage lattice search algorithm \cite{glover2015metaheuristics} which finds a set of candidate solutions to optimize the current distribution in the head given a focal vector field, that is, the target of stimulus. We consider the resulting parameter optimization problem as a metaheuristic task of computational intelligence \cite{glover2015metaheuristics}, wherein the goal is to find the best fitting solution with respect to one or more metacriteria which, in this document, are referred as Case (A) (\ref{sec:case_a}) for magnitude of the focused current density at the targeted stimulus location, and Case (B) (\ref{sec:case_b}) for ratio between focused and nuisance current intensity.

We compare the performance of the proposed L1-norm regularized L1-norm fitting (L1L1) method with L1-norm regularized L2-norm fitting (L1L2) version, a L2-norm equivalent fitted version obtained via semidefinite programming, and the Tikhonov regularized Least Squares (TLS) method \cite{dmochowski2017optimal}. On each method, the level of regularization and the relative nuisance field weight, which is explicit in L1L1 and L1L2 and embedded in the objective function in TLS, constitute the parameters of the candidate solution set. The comparison was performed by coupling the CVX optimization toolbox \cite{grant2009cvx} with MATLAB-based Zeffiro Interface (ZI) code package \footnote{\url{https://github.com/sampsapursiainen/zeffiro_interface}} \cite{he2019zeffiro, rezaei2021reconstructing, miinalainen2019realistic}, which allows creating a lead field matrix \cite{neuling2012finite} for a multi-compartment volume head model (\ref{app:Lead_Field}) using the Finite Element Method (FEM) together with Complete Electrode Model (CEM) \cite{pursiainen2012, pursiainen2017advanced} boundary conditions (\ref{sec:BConditions}). As a test domain, we used a realistic head model obtained from an openly available MRI data set (\ref{sec:Test_domain}).

The numerical results obtained, over a 8 and 20 active electrode montage, supports our initial hypothesis; compared to L1L2 and TLS, the L1L1 is shown to be advantageous with respect to both metacriteria cases. We propose that the metaheuristic L1L1 optimization approach, presented in this study, provides a potential alternative to determine the stimulation montage \cite{fernandez2020unification}.

This article is organized as follows: methodological details, including the optimization techniques two-stage metaheuristic search, test domain, and target placement are described in Section \ref{sec:Method}. The results can be found in Section \ref{sec:Results} and discussion in Section \ref{sec:Discussion}. The mathematical grounds and principles of tES forward modelling and weighted least-squares are explained in \ref{app:forward_model}.

\section{Methods}
\label{sec:Method}
The inverse problem of tES is to find a current pattern ${\bf y} = (y_1, y_2, \ldots, y_{\ell})$ that can generate a discretized current field ${\bf x} = (x_1, x_2, \ldots, x_{N})$ in the brain. The current pattern is required to meet Kirchhoff's current conservation conditions, that is $\sum_{\ell = 1}^{L} y_\ell = 0$ or  ${\bf 1}^T {\bf y} = 0 $ with ${\bf 1} = (1,1,\ldots,1)^T$, to yield both a total dose $\| {\bf y} \|_1 = \sum_{\ell = 1}^{L} |y_\ell| $ smaller than or equal to a given safety current limit ${\bm \mu}$ \cite{nitsche2003safety, KHADKA202069} and to have an entry-wise upper-limit of less or equal to ${\bm \gamma}$, i.e., $ {\bf y} \preceq \gamma {\bf 1}$. Here, we assume that $\gamma = \mu/2$, that is, the maximum absolute total dose can be achieved by a system with two or more electrodes ---one bearing positive polarity and the other with negative. The fitting between vectors ${\bf y}$ and ${\bf x}$ is enabled by the matrix equation
\begin{equation}
\label{optimization_task}
    {\bf L} {\bf y} = {\bf x},
\end{equation}
where ${\bf L}$ is a linear mapping following from a discretization of Maxwell's equations (\ref{app:forward_model}).

We consider finding optimized current pattern which, when applied into the head model ${\bf \Omega}$ through a given number of active electrodes attached on the scalp, generates a focused volume current distribution matching a synthetic dipolar current at a given orientation and location within the brain, while the nuisance field component remains suppressed. To enable an even comparison between different optimized current patterns, the total dose of each pattern is equaled to $\mu=$ 4 mA ($\gamma =$ 2 mA) \cite{KHADKA202069}.

\subsection{Optimization}
\label{sec:optimization}
To approximately solve (\ref{optimization_task}), we consider a weighted optimization scheme \cite{dmochowski2017optimal}, where the equation (\ref{optimization_task}) is split into two different components as
\begin{equation}
{\bf L} = \begin{pmatrix} {\bf L}_1 \\ {\bf L}_2 \end{pmatrix} \quad \hbox{and} \quad {\bf x} = \begin{pmatrix} {\bf x}_1 \\ {\bf 0} \end{pmatrix}.
\end{equation}
We call the first solution component ${\bf L}_1 {\bf y}$ the \textit{focused field}, i.e., the part that contains the given stimulus target, and the second one ${\bf L}_2 {\bf y}$ the \textit{nuisance field}, i.e., the remaining part of the field which we aim to suppress. To limit the number of non-zero currents in the current pattern ${\bf y}$, the objective function of the optimization task is regularized (penalized) by the norm of the current pattern ${\bf y}$. We control the magnitude difference between ${\bf L}_1 {\bf y}$ and ${\bf L}_2 {\bf y}$ by varying the weight of the nuisance field which is given explicitly for L1L1 and L1L2, and given as a penalty parameter embedded in the objective function for TLS, where explicit constraints are not applicable.

\subsubsection{L1-norm regularized L1-norm fitting}
\label{sec:L1L1}
We propose solving the following L1-norm regularized L1-norm fitting problem (L1L1)
\begin{equation}
    \min_{{\bf y}} \{ \, \| {\bf L}_1 {\bf y} - {\bf x}_1 \|_1 + \| {\bf L}_2 {\bf y}\|_1 + {\alpha} \zeta \| {\bf y} \|_1 \},
\label{objective_function}
\end{equation}
subject to $\| {\bf L}_2 {\bf y}\|_\infty \geq \varepsilon \nu$, ${\bf y}  \preceq \gamma {\bf 1}$, $\| {\bf y}\|_1 \leq \mu$, and $\sum_{\ell = 1}^L y_\ell = 0$. Here ${\bm \alpha}$ is the regularization parameter,  ${\bm \zeta} = \| {\bf L} \|_1$ and ${\bm \nu} = \| {\bf x} \|_\infty$ are scaling factors, and ${\bm \varepsilon}$ is the relative weight (numerical tolerance) of the nuisance field. Problem (\ref{objective_function}) constitutes the following linear programming task \cite[294]{Boyd2004}:
\begin{equation}
    \min_{{\bf y}, {\bf t}^{(1)}, {\bf t}^{(2)},{\bf t}^{(3)}}  \left( \sum_{k = 1}^N t^{(1)}_k + \sum_{m = 1}^M t^{(2)}_m+\alpha \zeta \sum_{\ell = 1}^L t^{(3)}_\ell \right) 
\end{equation}
subject to
{\setlength\arraycolsep{2 pt}
\begin{eqnarray} -\begin{pmatrix} {\bf t}^{(1)} \\ {\bf t}^{(2)} \\ \alpha \zeta {\bf t}^{(3)} \end{pmatrix} & \preceq & \begin{pmatrix} {\bf L}_1 \\ {\bf L}_2 \\  {\bf I} \end{pmatrix}  {\bf y } - \begin{pmatrix} {\bf x }_1 \\ {\bf 0 } \\ {\bf 0} \end{pmatrix} \preceq \begin{pmatrix} {\bf t}^{(1)} \\ {\bf t}^{(2)} \\ \alpha \zeta {\bf t}^{(3)} \end{pmatrix},  \nonumber\\
\begin{pmatrix} {\bf 0} \\ {\varepsilon}  \nu {\bf 1}\end{pmatrix} & \preceq &
\begin{pmatrix} {\bf t}^{(1)} \\ {\bf t}^{(2)}  \end{pmatrix}   \nonumber\\
\begin{pmatrix} {\bf 0} \\ {\bf 0} \end{pmatrix} & \preceq &
\begin{pmatrix} {\bf t}^{(3)} \\ {\bf 1}^T {\bf t}^{(3)} \end{pmatrix}  \preceq  \begin{pmatrix}  \gamma {\bf 1} \\ \mu  \end{pmatrix} \nonumber\\
{\bf 1}^T {\bf y} &  = &  0.
\label{linear_programming_1}
\end{eqnarray}}
Here, ${\bf t}^{(1)}$,  ${\bf t}^{(2)}$ and ${\bf t}^{(3)}$ constitute auxiliary {\textbf{N}}-by-1, {\textbf{M}}-by-1 and {\textbf{L}}-by-1 vectors, respectively. A numerically implementable form of (\ref{linear_programming_1}) with one inequality and equality constraint can be expressed  as follows:
\begin{equation} 
    \min_{{\bf y},  {\bf t}^{(1)}, {\bf t}^{(2)}, {\bf t}^{(3)}} \!\! \begin{pmatrix} {\bf 0}  \\  {\bf 1} \\ {\bf 1} \\  {\bf 1} \end{pmatrix}^T \!\! \begin{pmatrix} {\bf y} \\ {\bf t}^{(1)} \\  {\bf t}^{(2)}  \\ {\bf t}^{(3)}    \end{pmatrix} \end{equation} subject to 
{\setlength\arraycolsep{2pt} \begin{eqnarray}
     \begin{pmatrix} {\bf L}_1  & -{\bf I} & {\bf 0}  & {\bf 0}  \\ 
     {\bf L}_2  & {\bf 0} & -{\bf I} & {\bf 0}  \\ 
     -{\bf I} & {\bf 0} &  {\bf 0} & -{\bf I} \\ 
     -{\bf L}_1  & -{\bf I} &  {\bf 0}& {\bf 0}  \\
     -{\bf L}_2  & {\bf 0} & -{\bf I} & {\bf 0}  \\
     {\bf I} & {\bf 0} & {\bf 0}  & -{\bf I} \\
     {\bf 0} & -{\bf I} & {\bf 0}  & {\bf 0} \\     
     {\bf 0} & {\bf 0} & -{\bf I}  & {\bf 0} \\ 
     {\bf 0} & {\bf 0} & {\bf 0} & -{\bf I} \\
     {\bf 0} & {\bf 0} & {\bf 0} & {\bf I} \\ 
        {\bf 0} & {\bf 0} & {\bf 0} & {\bf 1}^T \\ 
     \end{pmatrix} \begin{pmatrix} {\bf y} \\ {\bf t}^{(1)} \\ {\bf t}^{(2)} \\ {\bf t}^{(3)}  \end{pmatrix}  & \preceq & \begin{pmatrix} {\bf x}_1 \\ {\bf 0} \\ {\bf 0} \\ -{\bf x_1} \\ {\bf 0} \\ {\bf 0} \\ {\bf 0} \\ -\varepsilon \nu {\bf 1} \\ {\bf 0} \\ \gamma  {\bf 1} \\ \mu  \\  \end{pmatrix} \nonumber\\
     {\bf 1}^T {\bf y} & = & 0. 
     \label{eq:LP} 
\end{eqnarray}}
The solution is found via primal-dual interior-point algorithm \cite{Boyd2004,tutuncu2003solving} of the SDPT3 package, accessible via the open CVX toolbox\footnote{\url{http://cvxr.com/cvx/}} \cite{grant2009cvx}.

\subsubsection{L1-norm regularized L2-norm fitting}
\label{sec:L1L2}
For the following L1-norm regularized L2-norm fitting problem (L1L2),
\begin{equation}
    \min_{{\bf y}} \{ \, \| {\bf L}_1 {\bf y} - {\bf x}_1 \|_2 + \| {\bf L}_2 {\bf y}\|_2  +  {\alpha} \zeta \| {\bf y} \|_1  \},
\end{equation}
the L1-norm fitting in (\ref{objective_function}) has been substituted with L2-norm, while the L1-norm regularization and linear constraints (\ref{linear_programming_1}) are the same as in the L1L1 approach. In CVX, the solution is obtained through semidefinite programming incorporating both linear and quadratic constraints of which the latter follow from L2-norm fitting in a straightforward manner \cite{tutuncu2003solving}.

\subsubsection{Tikhonov regularized least-squares}
\label{sec:TLS}
In Tikhonov regularized least-squares (TLS) estimation \cite{dmochowski2011optimized, dmochowski2017optimal}, the optimization problem to be solved is 
\begin{equation}
\min_{{\bf y}}\{ \| {\bf L}_1 {\bf y} - {\bf x}_1 \|^2_2 +  \alpha^2 \delta^2   \| {\bf L}_2 {\bf y} \|^2_2 + \alpha^2 \sigma^2  \| {\bf y}\|^2_2\},
\label{TLS_minimization}
\end{equation} 
where ${\bm \sigma} = \| {\bf L} \|_2$. To enforce focality, targeting the stimulus location and other areas of the brain, the nuisance field weight ${\bm \delta} \geq 0$ is considered as a variable parameter. The solution of (\ref{TLS_minimization}) is given by the linear system
\begin{equation}
\left(  {\bf L}_1^T {\bf L}_1  + \delta^2 \alpha^2 {\bf L}_2^T {\bf L}_2 + \alpha^2 \sigma^2 {\bf I} \right) {\bf y} =  {\bf L}_1^T {\bf x}_1\,.
\label{TLS_formula}
\end{equation}
which one can solve numerically using Matlab's backslash (\texttt{\char`\\}) operator.

It can be shown that the focused current density${\bm \Gamma}$ (unit A/m\textsuperscript{2}),
\begin{equation}
    \label{gamma_def}
    \Gamma = \frac{{\bf x}_1^T {\bf L}_1 {\bf y}}{\| {\bf x}_1 \|_2}
\end{equation} 
in the direction of the targeted brain activity yields its maximum when $\delta = 0$ and that the ratio, ${\bm \Theta}$ (unitless),
\begin{equation}
\label{theta_def}
    \Theta = \frac{ \Gamma}{ \| {\bf L}_2 {\bf y} \|_2/\sqrt{M}}\,. 
\end{equation}
between $\Gamma$ and the average nuisance field magnitude increases along with the value of $\delta$ when $\delta$ is a small positive number (further details can be found in \ref{effect_of_weighting}). Notice that for $\delta >0$, the minimization problem (\ref{TLS_minimization}) can be written in the following alternative form
\begin{equation}
    \min_{\bf y}\{ \kappa \| {\bf L}_1 {\bf y} - {\bf x}_1 \|^2_2 + \| {\bf L}_2 {\bf y} \|^2_2 + (\sigma^2/\delta^2) \| {\bf y}\|^2_2\}
\end{equation}
with ${\bm \kappa} = 1/(\alpha^2 \delta^2)$. Therefore, either focused or nuisance field can be weighted.

\subsection{Candidate solution set}
\label{sec:candidate_solution}
The above formulations of L1L1, L1L2 and TLS optimization problem depends on the regularization parameter and the nuisance field weight, i.e.,  the effect the nuisance field component has on the solution of the optimization problem. In both L1L1 and L1L2 methods, the weight is given explicitly while in TLS the nuisance field component is expressed in a weighted form, leading to a different dependence of the solution on the weight as compared to the previous two methods.

To find the optimal case-wise parameter combination we examine a two-dimensional 36 $\times$ 36  lattice of optimized candidate solutions covering a wide 180 dB dynamical range with 5 dB increments for each parameter value. In L1L1,  $\alpha$ and $\varepsilon$ are varied between -160 and 20 dB, in L1L2 between -140 and 40 dB, and in TLS, the $\alpha$ and $\delta$ varied between -240 and -60 dB and -100 and 80 dB, respectively. Parameter variation is considered necessary to obtain the best possible performance, since the scale of the objective function is affected by both parameters. The lattice resolution was selected so that the total computing time was maintained on an  acceptable level. The maximal uncertainty related to the lattice was estimated after obtaining the candidate solutions.

To compute the candidate solution, we employ a Dell 5820 workstation equipped with Intel Core i9-10900X processor and 256 GB RAM. The total computing time required to evaluate a full lattice of candidate solutions was 7390, 11134, and 138 seconds with L1L1, L1L2, and TLS, respectively. L1L1 and L1L2 utilized a single thread mode while TLS was automatically parallelized by Matlab's interpreter. A relative solver tolerance of 1E-10 was used as a stopping criterion of L1L1 and L1L2.

\subsection{Two-stage metaheuristic search}
\label{sec:two_stage_lattice}
To filter the set of candidate solutions, we perform a two-stage metaheuristic lattice search, where focused current density $\Gamma$ and current ratio $\Theta$ are used as metacriteria.  Of the following two cases, (A) utilizes both criteria, while (B) constitutes a reference for maximizing $\Gamma$.

\subsubsection{Case (A)}
\label{sec:case_a}
The first stage sets a threshold condition $\Gamma \geq 0.11$ A/m\textsuperscript{2}, and the second stage maximizes the thresholded set of candidates with respect to $\Theta$. By using these two criteria, we ensure that the selected candidate will have adequate current magnitude \cite{khanindividually} in the targeted position and appropriately suppressed nuisance field component at the same time.

\subsubsection{Case (B)}
\label{sec:case_b}
For comparison, we consider a simpler scheme in which the focused current density $\Gamma$ alone is maximized over the full candidate set. That is, the search is based on a single criterion and stage.

\subsubsection{Post-optimization with non-fixed vs.\ fixed montage}
\label{sec:post_opt}
Aiming at the best possible optimization outcome, each search run is performed twice: in the first run all the current channels are present in the optimization process while the second one uses a limited montage which is selected based on the first run; the electrodes $\ell$ with the greatest current $|I_{\ell}|$ contribution to the total maximum current value are selected to carry non-zero amplitudes, while the remaining ones are opted out (set to zero) from the second run. In the first run, we apply a cap of 128 electrode positions which are reduced to 8 and 20 active channel montages in the second one. These channel counts are inspired by the commercial state-of-the-art tES systems  \cite{roy2019integration,tost2021choosing}.

\subsection{Synthetic sources and placement}
\label{sec:synthetic_source}
The amplitude of the dipolar target current is related to the corresponding local current density in the brain. As reference, the cortex thickness was set with 4 mm (millimeter) \cite{Fischl11050} and the activity density with 0.77 nAm/mm\textsuperscript{2} (nanoampere per square millimeter) \cite{MURAKAMI201549}. Three 10 nAm dipoles were placed in the following three left hemispheric Brodmann's areas \cite{brodmann2007}: postcentral gyrus (red), superior temporal gyrus (cyan), and occipital lobe (blue) (Fig. \ref{fig:ZI_Parcellations}). Each dipole is oriented normally with respect to the surface of the gray matter tissue to satisfy the normal constraint of the brain activity in the cerebral cortex \cite{creutzfeldt1962influence}. The two-stage metaheuristic search and the optimization routines were conducted using a individual dipole throughout the head model resulting in three independent numerical solutions; in this document, the results are categorized as \textit{Somatosensory} for postcentral gyrus, \textit{Auditory} for superior temporal gyrus, and \textit{Visual} for occipital lobe.

\subsection{Test domain}
\label{sec:Test_domain}
As a test domain of the numerical experiments, we used a multi-compartment volume conductor head model based on openly available anatomical T1-weighted Magnetic Resonance Imaging (MRI) data\footnote{\url{https://brain-development.org/ixi-dataset/}} obtained from a real subject. The data were segmented using FreeSurfer Software Suite{\footnote{\url{https://surfer.nmr.mgh.harvard.edu/}}} which distinguishes different head and brain tissue compartments including scalp, skull, cerebrospinal fluid, gray and white matter as well as subcortical structures such as brain stem, thalamus, amygdala, and ventricles with their own complex geometrical properties \cite{Malmivuo1997}. The volume segmentation was obtained using ZI's Finite Element (FE) mesh generator \cite{he2019zeffiro} which identifies the compartments obtained from the surface segmentation and creates a smoothed and optimized FE mesh composed of these compartments. 

\begin{figure}[ht]
\centering
    \begin{subfigure}[t]{4.2cm}
        \centering
        \includegraphics[width=4.0cm]{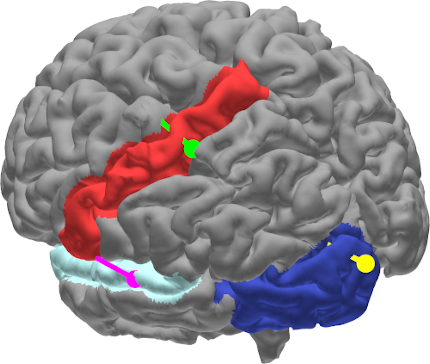}
        \caption{Gray matter and the highlighted regions of interest}.
        \label{fig:ZI_Parcellations}
    \end{subfigure}
    \hskip0.3cm
    \begin{subfigure}[t]{3.8cm}
        \centering
        \includegraphics[width=2.8cm]{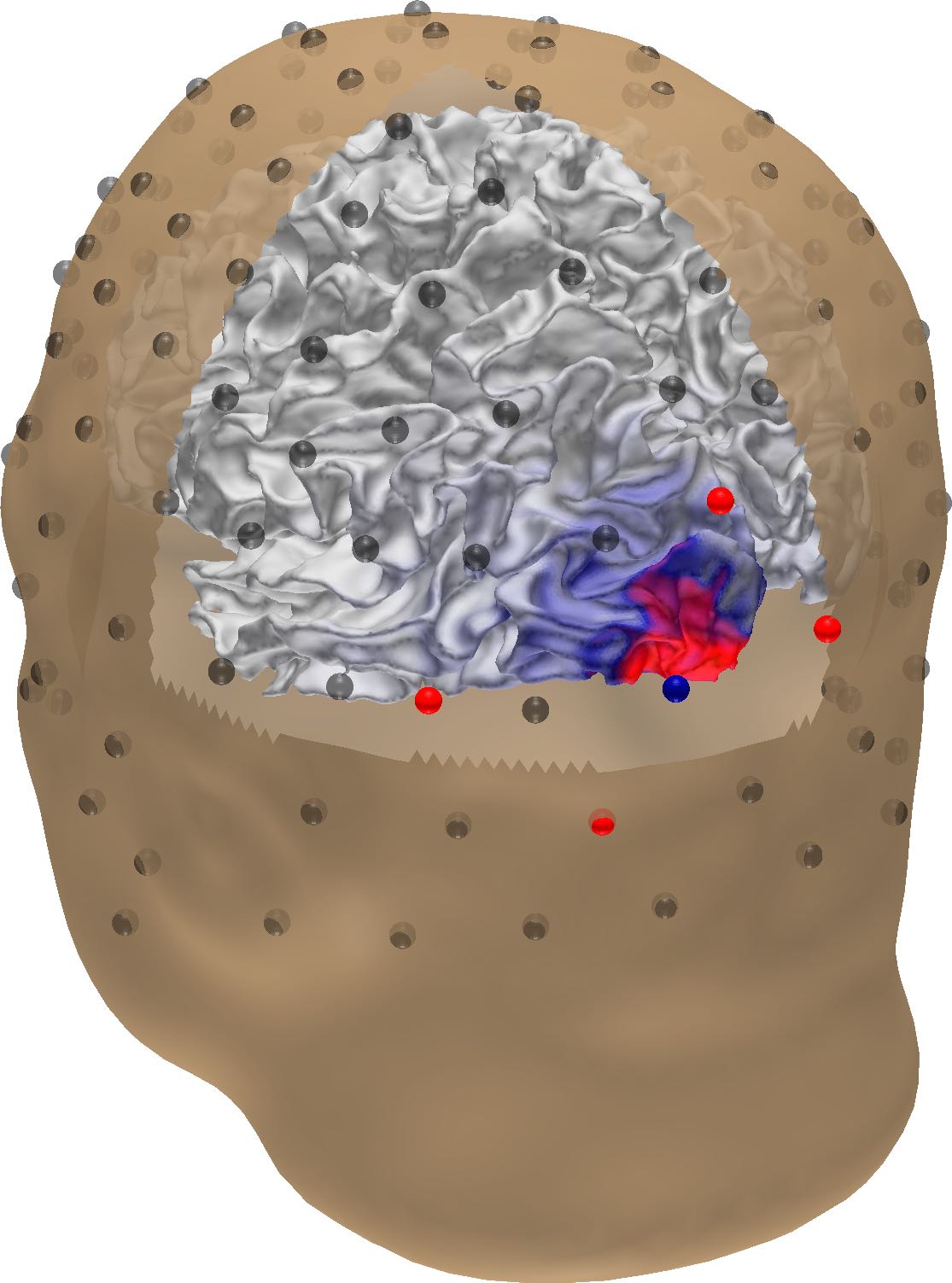}
        \caption{volume current distribution and electrodes.}
        \label{fig:white_matter_example}
    \end{subfigure}
\caption{{\bf (\subref{fig:ZI_Parcellations})}: Left-posterior view of the gray matter of the volume conductor model highlighting the three regions of interest: postcentral gyrus (red), superior temporal gyrus (cyan), and occipital lobe (blue)}, as given by the 36 label Desikan-Killiany atlas. Each region is presented with their own point-like dipolar target current (sphere), and orientation (line). {\bf (\subref{fig:white_matter_example})}: Example of a simple five-channel electrode montage creating a volume current distribution in the left-hemisphere of the occipital lobe. The positive (or anodal) electrode channels injecting the current into the domain are illustrated with red spheres, the negative (or cathodal) with blue, and the disabled/inactive with dark gray.
\label{fig:parcellation_examples}
\end{figure}

To discretize the head mode, we use a finite element mesh resolution of 1 mm to obtain physiologically accurate results \cite{rullmann2009eeg}. The conductivity distribution influences the accuracy of the forward solution \cite{Montes-Restrepo2014}. In our model, the conductivity is constant in each tissue compartment with the values corresponding to the set proposed in \cite{dannhauer2010}. The placement of the 128 EEG/tES electrodes (Fig. \ref{fig:white_matter_example}) follows the International 10-20, 10-10, or 10-5 EEG hardware system \cite{JURCAK20071600}. Electrode impedance was set to be 2 kOhm (kiloohms) uniformly. Impedance modelling was enabled by the incorporation of the complete electrode model into the forward model. The tES lead field matrix ${\bf L}$ was generated as explained in (\ref{app:Lead_Field}) for 1000 uniformly randomly selected spatial set of points contained by the gray matter compartment, including three Cartesian degrees of freedom per point.

\subsubsection{Accuracy analysis}
\label{sec:Accuracy_analysis}
The optimization outcome was examined by evaluating the maximum injected current in the optimized current pattern $\| {\bf y} \|_\infty$, the focused current density $\Gamma$, the current ratio $\Theta$, and the Angle Difference (AD) between the focused and the targeted fields, i.e., 
\begin{equation}
\hbox{AD}(\vec{j}_1,\vec{j}_2) = \arccos\,\left(\frac{\left\langle \vec{j}_1,\vec{j}_2\right\rangle}{\left\|\vec{j}_1\right\|\left\|\vec{j}_2\right\|}\right) \end{equation}
with $\vec{j}_1$ representing the volume current distribution at the target location generated by the injected pattern, and $\vec{j}_2$ the dipolar target current, respectively. The limits for lattice-induced deviation of  $\| {\bf y} \|_\infty$, $\Gamma$, $\Theta$, AD, were estimated by forming a second order Taylor's polynomial approximation in the 3-by-3 lattice centered at the selected candidate solution. These limits were obtained as the maximum deviation of the polynomial within a co-centered 3-by-3 environment of a hypothetical lattice with double the resolution compared to the actual one.

\section{Results}
\label{sec:Results}
\begin{figure*}[h!]
\centering
    \begin{scriptsize}
    \centering
    \begin{minipage}{18.0cm}
        \centering
        \begin{minipage}{0.2cm}
        \centering
           \mbox{}
        \end{minipage}
        \begin{minipage}{16.50cm}
        \centering
            \begin{minipage}{0.25cm}
                \mbox{}
            \end{minipage}
            \begin{minipage}{5.30cm}
            \centering
                \textbf{L1-norm regularized \\ L1-norm fitting (L1L1)}
            \end{minipage}
            \begin{minipage}{5.30cm}
                \centering
                \textbf{L1-norm regularized \\ L2-norm fitting (L1L2)}
            \end{minipage}
            \begin{minipage}{5.30cm}
                \centering
                \textbf{Tikhonov regularized \\ least-squares (TLS)}
            \end{minipage}
        \end{minipage}
    \end{minipage}
    \vskip0.1cm \hrule \vskip0.1cm
    \begin{minipage}{18.0cm}
    \centering
        \begin{minipage}{0.2cm}
        \centering
            \rotatebox{90}{\textbf{Somatosensory}}
        \end{minipage}
        \begin{minipage}{16.50cm}
        \centering
            \begin{minipage}{0.25cm}
            \centering
                \rotatebox{90}{\textbf{Case A}}
            \end{minipage}
            \begin{subfigure}{5.3cm}
            \centering
                \begin{minipage}{2.6cm}
                    \centering
                    \includegraphics[height = 1.75cm]{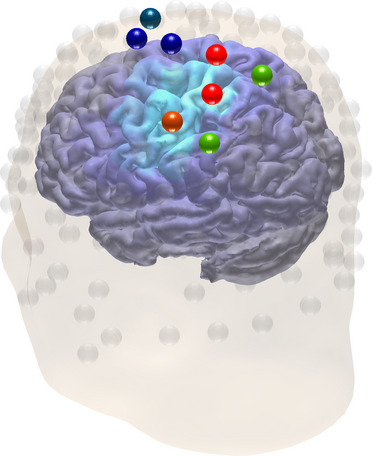} \\
                    \includegraphics[width = 2.30cm,height=1.2cm]{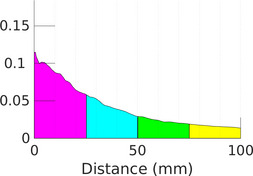} 
                \end{minipage} 
                \begin{minipage}{2.6cm}
                \centering 
                    \includegraphics[height = 2.55cm]{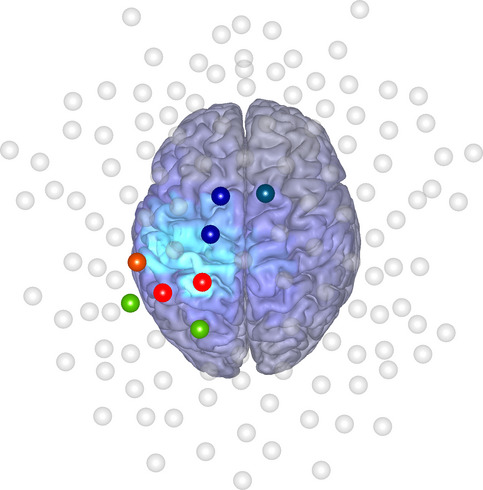}
                \end{minipage} \\
                \label{fig:postcentral_LP_1_8_results}
            \end{subfigure}
            \begin{subfigure}{5.3cm}
            \centering
                \begin{minipage}{2.60cm}
                    \centering
                    \includegraphics[height = 1.75cm]{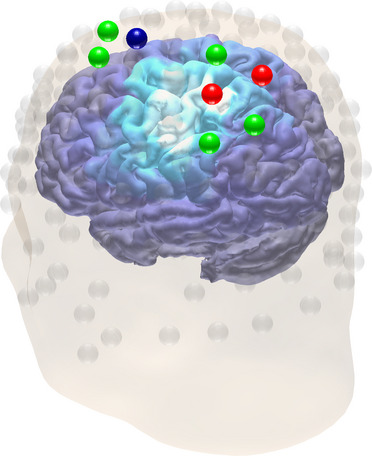} \\
                    \includegraphics[width = 2.30cm,height=1.2cm]{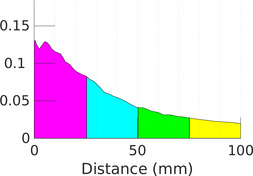} 
                \end{minipage} 
                \begin{minipage}{2.60cm}
                \centering 
                    \includegraphics[height = 2.55cm]{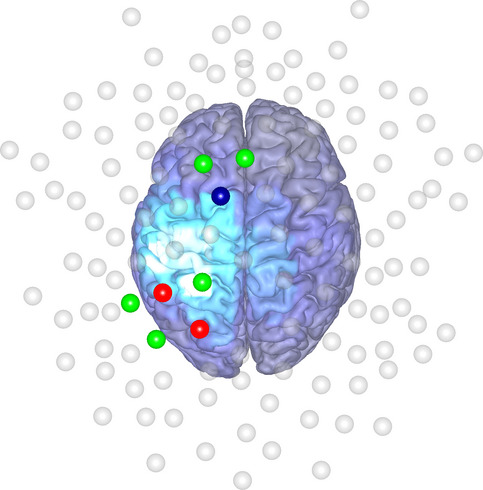}
                \end{minipage}
                \label{fig:postcentral_SDP_1_8_results}
            \end{subfigure}
            \begin{subfigure}{5.3cm}
            \centering
                \begin{minipage}{2.60cm}
                    \centering
                    \includegraphics[height = 1.75cm]{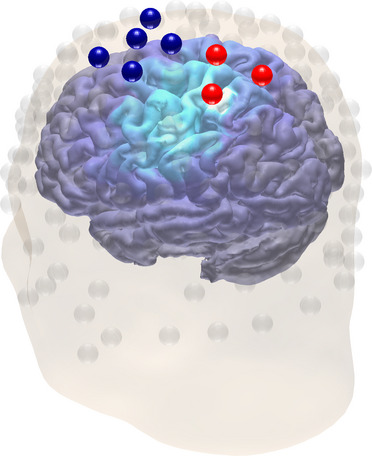} \\
                    \includegraphics[width = 2.30cm,height=1.2cm]{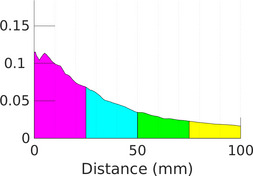}
                \end{minipage} 
                \begin{minipage}{2.60cm}
                \centering 
                    \includegraphics[height = 2.55cm]{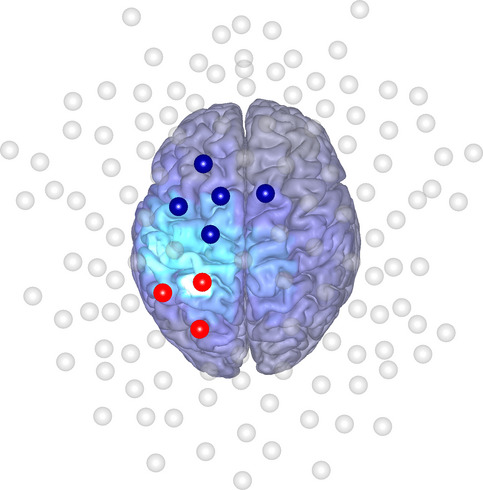}
                \end{minipage}
                \label{fig:postcentral_TLS_1_8_results}
            \end{subfigure}
            \vskip0.15cm
            \begin{minipage}{0.25cm}
            \centering
                \rotatebox{90}{\textbf{Case B}}
            \end{minipage}
            \begin{subfigure}{5.3cm}
            \centering
                \begin{minipage}{2.60cm}
                \centering
                    \includegraphics[height = 1.75cm]{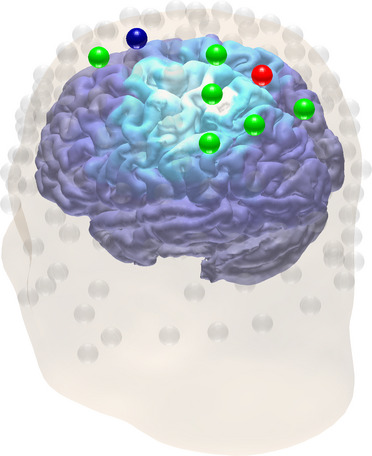} \\
                     \includegraphics[width = 2.30cm,height=1.2cm]{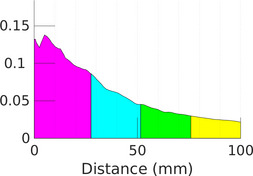}
                \end{minipage}
                \begin{minipage}{2.6cm}
                \centering 
                    \includegraphics[height = 2.55cm]{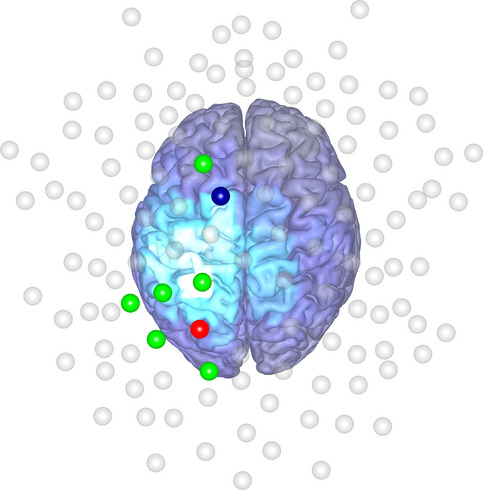}
                \end{minipage}
                \label{fig:postcentral_LP_2_8_results}
            \end{subfigure}
            \begin{subfigure}{5.3cm}
            \centering
                \begin{minipage}{2.60cm}
                    \centering
                    \includegraphics[height = 1.75cm]{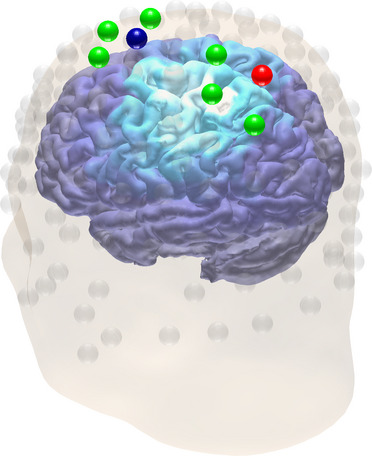} \\
                    \includegraphics[width = 2.30cm,height=1.2cm]{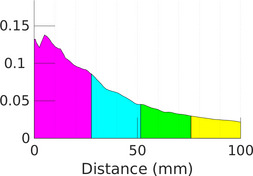}
                \end{minipage} 
                \begin{minipage}{2.60cm}
                \centering 
                    \includegraphics[height = 2.55cm]{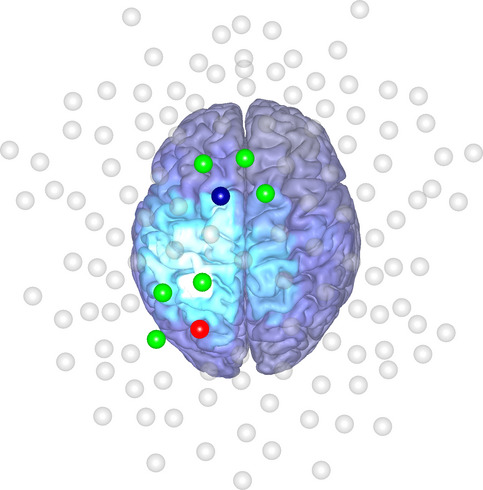}
                \end{minipage}
                \label{fig:postcentral_SDP_2_8_results}
            \end{subfigure}
            \begin{subfigure}{5.3cm}
            \centering
                \begin{minipage}{2.60cm}
                    \centering
                    \includegraphics[height = 1.75cm]{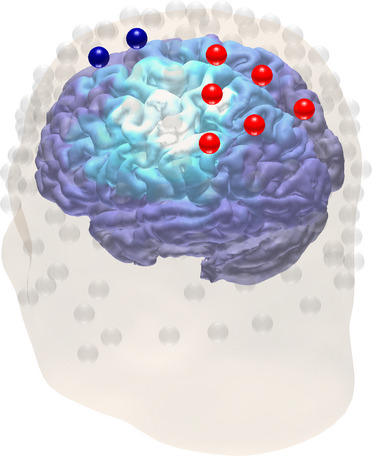} \\
                     \includegraphics[width = 2.30cm,height=1.2cm]{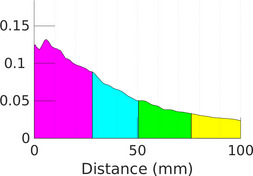} 
                \end{minipage} 
                \begin{minipage}{2.60cm}
                \centering 
                    \includegraphics[height = 2.55cm]{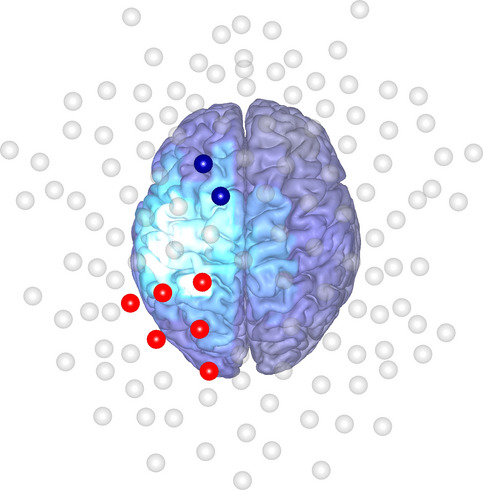}
                \end{minipage}
                \label{fig:postcentral_TLS_2_8_results}
            \end{subfigure}
        \end{minipage}
    \end{minipage}
    \vskip0.1cm \hrule \vskip0.1cm
    \begin{minipage}{18.0cm}
    \centering
        \begin{minipage}{0.2cm}
        \centering
            \rotatebox{90}{\textbf{Auditory}}
        \end{minipage}
        \begin{minipage}{16.50cm}
        \centering
            \begin{minipage}{0.25cm}
            \centering
                \rotatebox{90}{\textbf{Case A}}
            \end{minipage}
            \begin{subfigure}{5.3cm}
            \centering
                \begin{minipage}{2.6cm}
                    \centering
                    \includegraphics[height = 1.75cm]{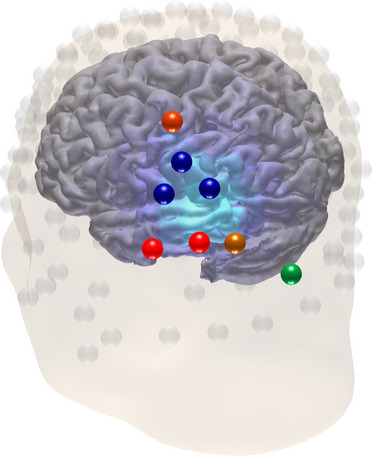} \\
                    \includegraphics[width = 2.30cm,height=1.2cm]{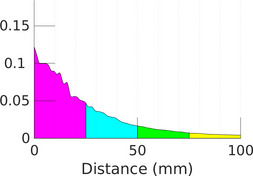} 
                \end{minipage} 
                \begin{minipage}{2.6cm}
                \centering 
                    \includegraphics[height = 2.55cm]{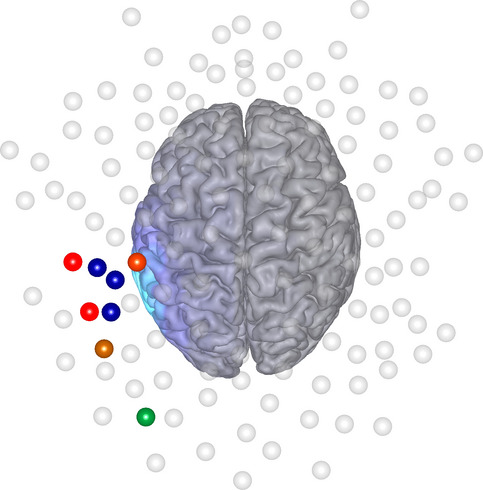}
                \end{minipage}
                \label{fig:auditory_LP_1_8_results}
            \end{subfigure}
            \begin{subfigure}{5.3cm}
            \centering
                \begin{minipage}{2.60cm}
                    \centering
                    \includegraphics[height = 1.75cm]{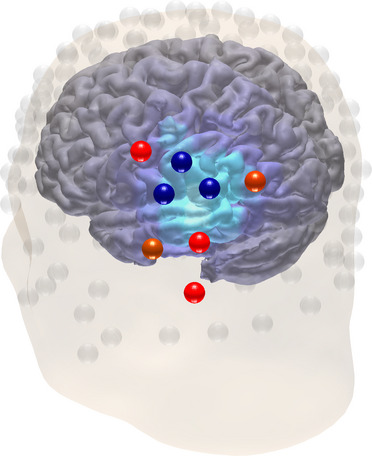} \\
                    \includegraphics[width = 2.30cm,height=1.2cm]{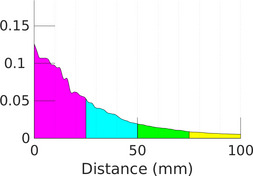} 
                \end{minipage}
                \begin{minipage}{2.60cm}
                \centering 
                    \includegraphics[height = 2.55cm]{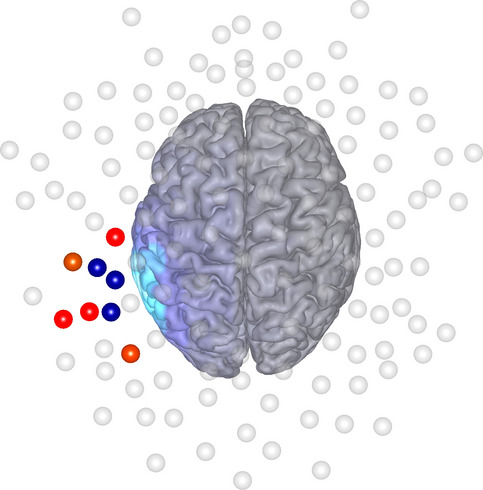}
                \end{minipage}
                \label{fig:auditory_SDP_1_8_results}
            \end{subfigure}
            \begin{subfigure}{5.3cm}
            \centering
                \begin{minipage}{2.60cm}
                    \centering
                    \includegraphics[height = 1.75cm]{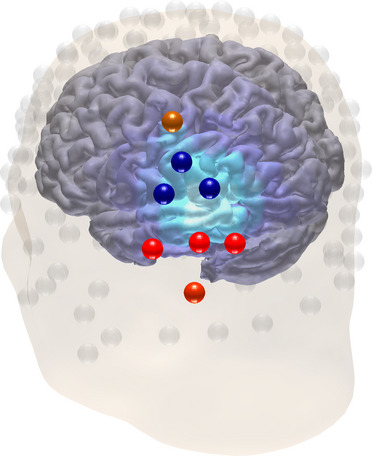} \\
                    \includegraphics[width = 2.30cm,height=1.2cm]{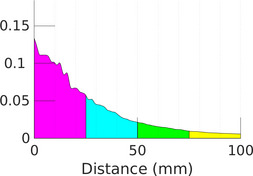}
                \end{minipage} 
                \begin{minipage}{2.60cm}
                \centering 
                    \includegraphics[height = 2.55cm]{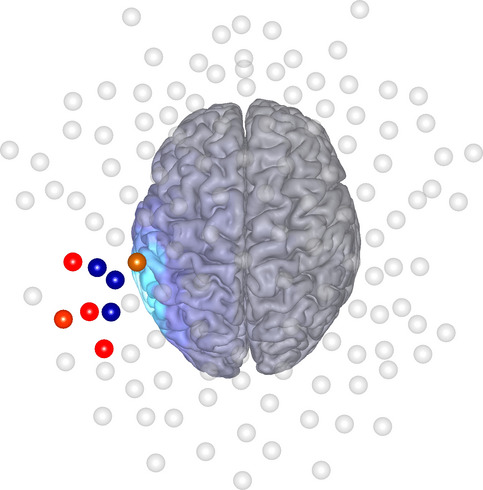}
                \end{minipage}
                \label{fig:auditory_TLS_1_8_results}
            \end{subfigure}
            \vskip0.15cm
            \begin{minipage}{0.25cm}
            \centering
                \rotatebox{90}{\textbf{Case B}}
            \end{minipage}
            \begin{subfigure}{5.3cm}
            \centering
                \begin{minipage}{2.60cm}
                    \centering
                    \includegraphics[height = 1.75cm]{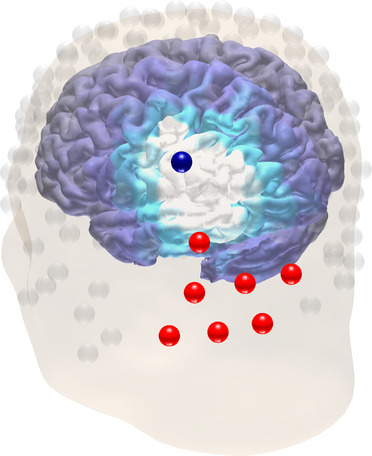} \\
                    \includegraphics[width = 2.30cm,height=1.2cm]{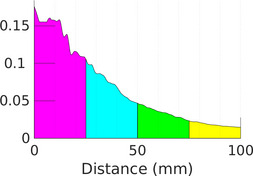}
                \end{minipage}
                \begin{minipage}{2.60cm}
                \centering 
                    \includegraphics[height = 2.55cm]{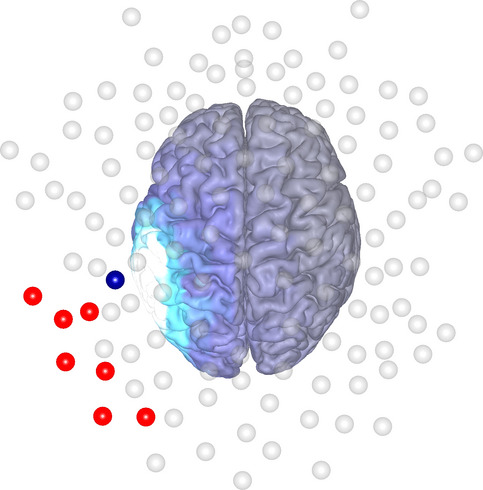}
                \end{minipage}
                \label{fig:auditory_LP_2_8_results}
            \end{subfigure}
            \begin{subfigure}{5.3cm}
            \centering
                \begin{minipage}{2.60cm}
                    \centering
                    \includegraphics[height = 1.75cm]{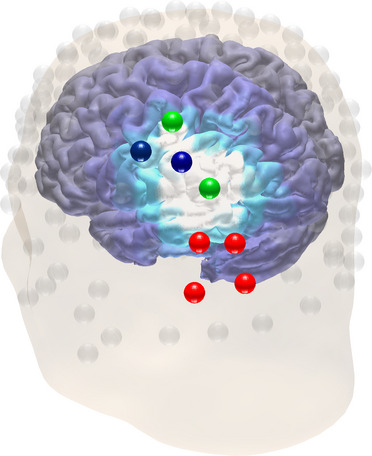} \\
                    \includegraphics[width = 2.30cm,height=1.2cm]{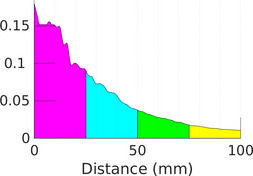}
                \end{minipage} 
                \begin{minipage}{2.60cm}
                \centering 
                    \includegraphics[height = 2.55cm]{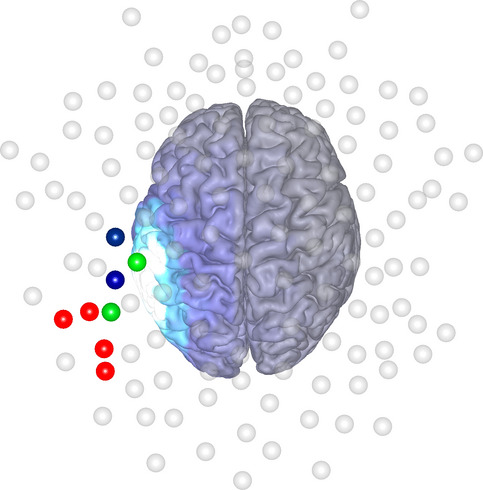}
                \end{minipage}
                \label{fig:auditory_SDP_2_8_results}
            \end{subfigure}
            \begin{subfigure}{5.3cm}
            \centering
                \begin{minipage}{2.60cm}
                    \centering
                    \includegraphics[height = 1.75cm]{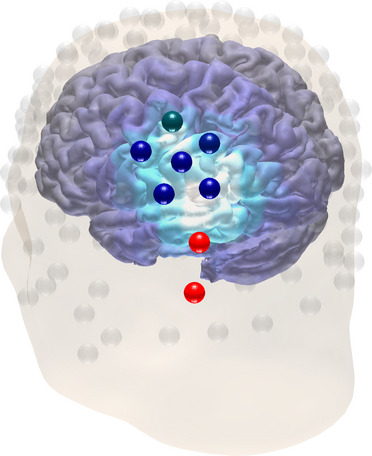} \\
                    \includegraphics[width = 2.30cm,height=1.2cm]{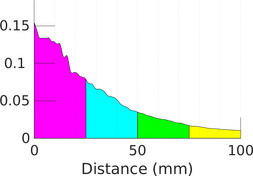} 
                \end{minipage} 
                \begin{minipage}{2.60cm}
                \centering 
                    \includegraphics[height = 2.55cm]{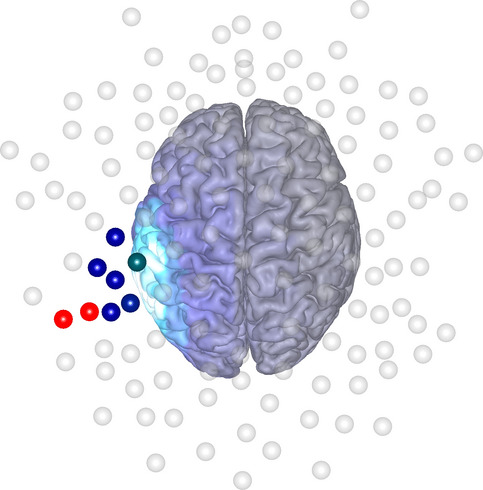}
                \end{minipage}
                \label{fig:auditory_TLS_2_8_results}
            \end{subfigure}
        \end{minipage}
    \end{minipage}
    \vskip0.1cm \hrule \vskip0.1cm
    \begin{minipage}{18.0cm}
    \centering
        \begin{minipage}{0.2cm}
        \centering
            \rotatebox{90}{\textbf{Visual}}
        \end{minipage}
        \begin{minipage}{16.50cm}
        \centering
            \begin{minipage}{0.25cm}
            \centering
                \rotatebox{90}{\textbf{Case A}}
            \end{minipage}
            \begin{subfigure}{5.3cm}
            \centering
                \begin{minipage}{2.60cm}
                \centering
                    \includegraphics[height = 1.75cm]{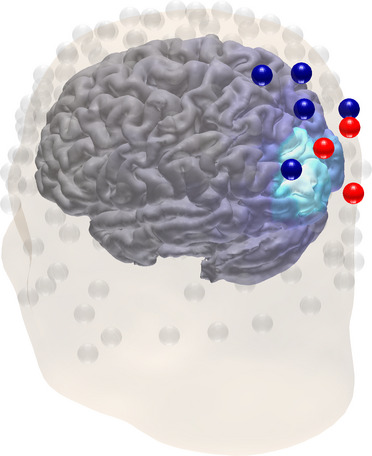} \\
                    \includegraphics[width = 2.30cm,height=1.2cm]{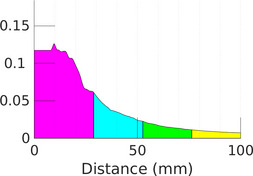} 
                \end{minipage} 
                \begin{minipage}{2.60cm}
                \centering 
                    \includegraphics[height = 2.55cm]{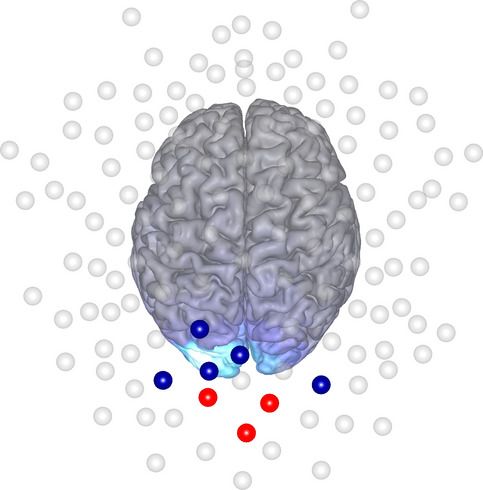}
                \end{minipage}
                \label{fig:occipital_LP_1_8_results}
            \end{subfigure}
            \begin{subfigure}{5.3cm}
            \centering
                \begin{minipage}{2.60cm}
                    \centering
                    \includegraphics[height = 1.75cm]{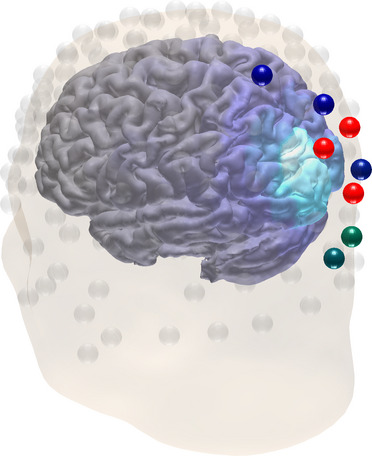} \\
                    \includegraphics[width = 2.30cm,height=1.2cm]{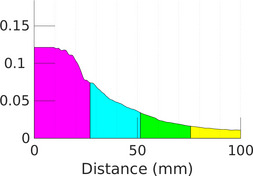} 
                \end{minipage} 
                \begin{minipage}{2.60cm}
                \centering 
                    \includegraphics[height = 2.55cm]{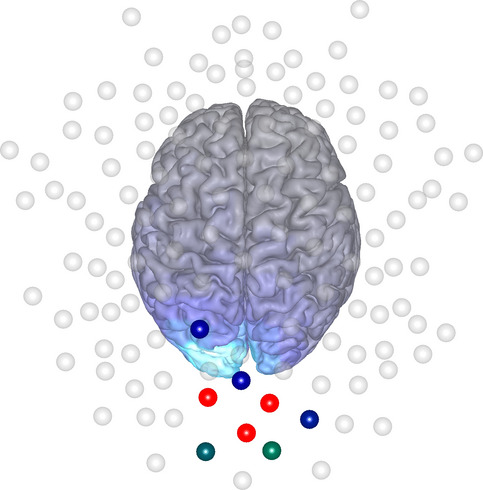}
                \end{minipage}
                \label{fig:occipital_SDP_1_8_results}
            \end{subfigure}
            \begin{subfigure}{5.3cm}
            \centering
                \begin{minipage}{2.60cm}
                    \centering
                    \includegraphics[height = 1.75cm]{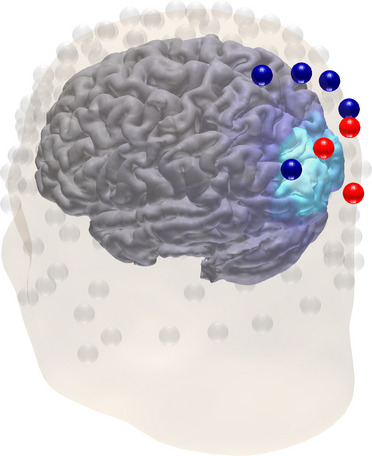} \\
                    \includegraphics[width = 2.30cm,height=1.2cm]{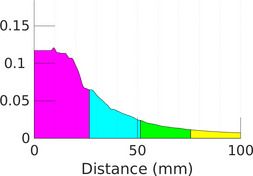}
                \end{minipage} 
                \begin{minipage}{2.60cm}
                \centering 
                    \includegraphics[height = 2.55cm]{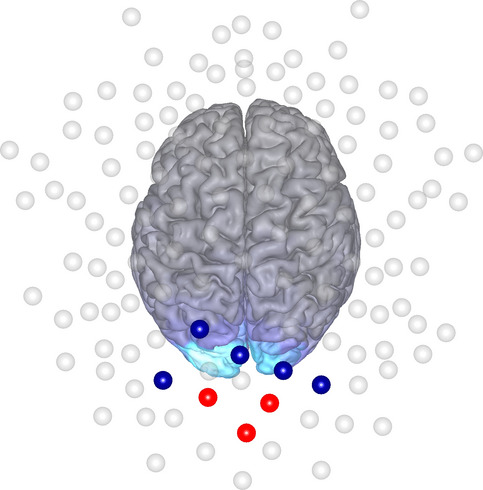}
                \end{minipage}
                \label{fig:occipital_TLS_1_8_results}
            \end{subfigure}
            \vskip0.15cm
            \begin{minipage}{0.25cm}
            \centering
                \rotatebox{90}{\textbf{Case B}}
            \end{minipage}
            \begin{subfigure}{5.3cm}
            \centering
                \begin{minipage}{2.60cm}
                    \centering
                    \includegraphics[height = 1.75cm]{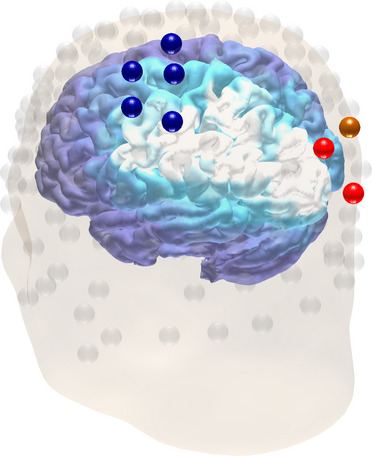} \\
                    \includegraphics[width = 2.30cm,height=1.2cm]{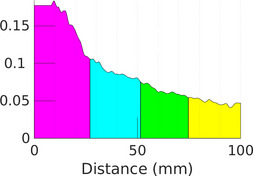}
                \end{minipage}
                \begin{minipage}{2.60cm}
                \centering 
                    \includegraphics[height = 2.55cm]{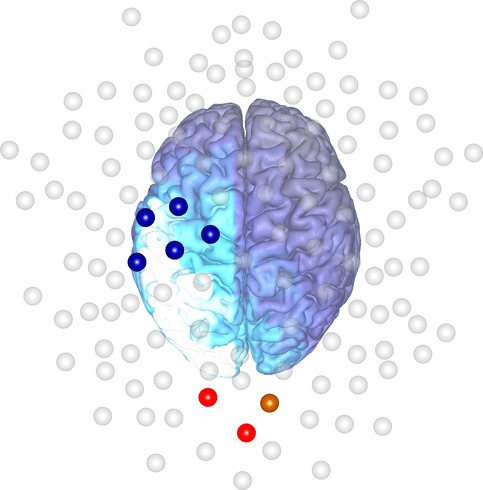}
                \end{minipage}
                \label{fig:occipital_LP_2_8_results}
            \end{subfigure}
            \begin{subfigure}{5.3cm}
            \centering
                \begin{minipage}{2.60cm}
                    \centering
                    \includegraphics[height = 1.75cm]{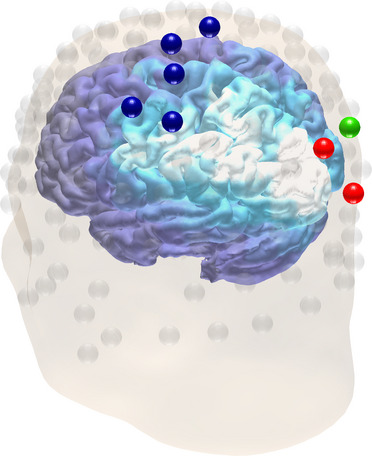} \\
                    \includegraphics[width = 2.30cm,height=1.2cm]{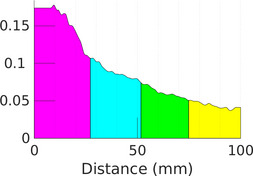}
                \end{minipage} 
                \begin{minipage}{2.60cm}
                \centering 
                    \includegraphics[height = 2.55cm]{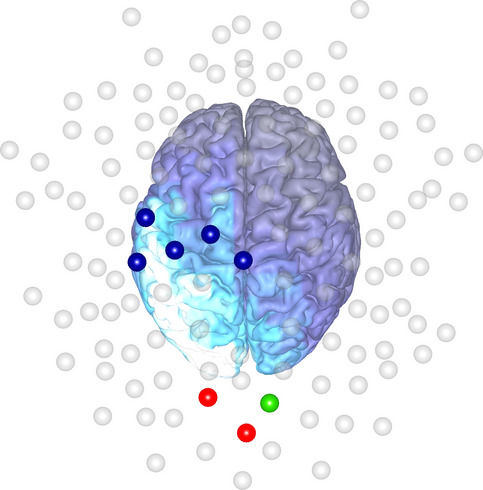}
                \end{minipage}
                \label{fig:occipital_SDP_2_8_results}
            \end{subfigure}
            \begin{subfigure}{5.3cm}
            \centering
                \begin{minipage}{2.60cm}
                \centering
                    \includegraphics[height = 1.75cm]{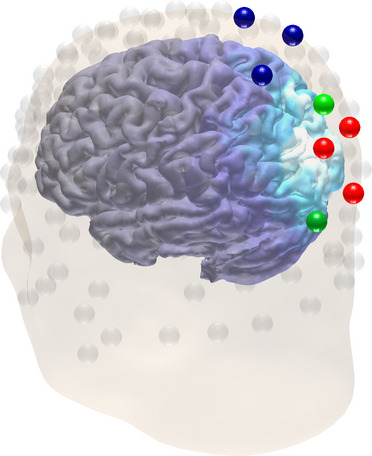} \\
                    \includegraphics[width = 2.30cm,height=1.2cm]{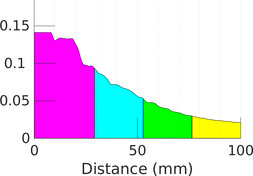} 
                \end{minipage} 
                \begin{minipage}{2.60cm}
                \centering 
                    \includegraphics[height = 2.55cm]{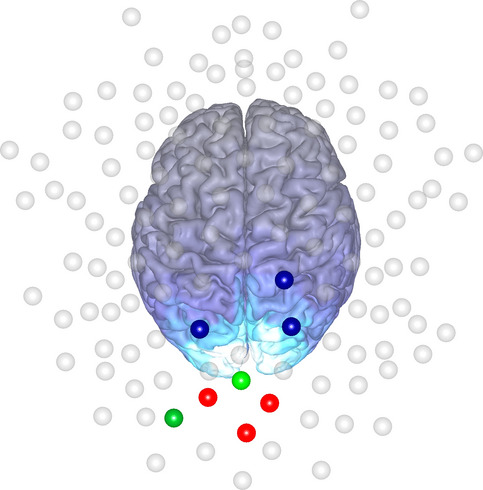}
                \end{minipage}
                \label{fig:occipital_TLS_2_8_results}
            \end{subfigure}
        \end{minipage}
    \end{minipage}
    \vskip0.1cm \hrule \vskip0.1cm
    \begin{minipage}{18.0cm}
    \centering
        \hskip0.75cm
        \includegraphics[height = 0.65cm, width = 3.5cm]{images_jpg/revision/thumbnail_volume_current_colorbar_horizontal.png}
        \hskip4.25cm
        \includegraphics[height = 0.65cm, width = 3.5cm]{images_jpg/revision/electrode_colorbar_2mA_hor.png}
    \end{minipage}
\end{scriptsize}
\caption{The 8-electrode montage and current pattern (mA) together with corresponding volume current density A/m\textsuperscript{2} obtained through two consecutive runs of the two-stage metaheuristic optimization process with respect to non-fixed and fixed montage of 8 active electrodes, respectively. Average magnitude in the direction of the target dipole is shown as a function of distance (mm) from the dipole position. The colorbar of the current pattern shows a color gradient for the interval from -0.25 to 0.25 mA to enhance the visibility of small variations in the pattern.}
\label{fig:results_8_channel}
\end{figure*}

Numerical comparison between L1L1, L1L2, and TLS have been included in Table \ref{table:method_comparison}. The outcome of their metaheuristic search process are shown in Figure \ref{fig:imagesc_L1L1}, \ref{fig:imagesc_L1L2} and \ref{fig:imagesc_TLS}, respectively. Stimulation accuracy in Figure \ref{fig:stem_comparison}, optimized current patterns (as well as their projections in the brain) for the 20 and 8 channel montage in Figure \ref{fig:results_20_channel} and \ref{fig:results_8_channel}, respectively.

\begin{table*}[ht]
\centering
    \caption{Optimization results for a 20- and 8-channel electrode montage with the target current placed at the somatosensory, auditory and visual regions of interest. The maximum injected current  $\| {\bf y} \|_\infty$ (mA), focused current $\Gamma$ (A/m\textsuperscript{2}) at the target location, angle difference AD (relative), current ratio $\Theta$ (relative) between focused and nuisance current. L1L1, L1L2 and TLS optimization method have been applied in finding the candidate solution for the metaheuristic search process. In case (A), $\Theta$ has been  maximized over those candidate solutions for which $\Gamma \geq 0.11$ A/m\textsuperscript{2}. In (B), the global maximizer of $\Gamma$ has been found. Each estimate has been associated with an estimate of maximal deviation due to lattice resolution as described in Section \ref{sec:Accuracy_analysis}.}
    \begin{scriptsize}
        \begin{tabular}{lllcccccccc} 
        \hline
        &  & \multicolumn{2}{c}{\textbf{Maximum}}   & \multicolumn{2}{c}{\textbf{Current}}   &  \multicolumn{2}{c}{\textbf{Angle}}   &  \multicolumn{2}{c}{\textbf{Current}}  \\
        &  & \multicolumn{2}{c}{\textbf{current $\| {\bf y} \|_\infty$ (mA)}}       & \multicolumn{2}{c}{\textbf{density $\Gamma$  (A/m\textsuperscript{2}) }}      & \multicolumn{2}{c}{\textbf{difference (deg)}}   & \multicolumn{2}{c}{\textbf{ratio $\Theta$ (rel.)} }     \\
        \textbf{Method}   &   \textbf{Chn}             & \textbf{Value}   &  \textbf{Deviation}      &  \textbf{Value} & \textbf{Deviation}  & \textbf{Value} & \textbf{Deviation}  & \textbf{Value} &  \textbf{Deviation} \\
        \hline   
        \multicolumn{10}{c}{\bf Somatosensory}  \\
        \hline
        L1L1 (A)  & 20 &  1.13  &  9.9E-02  &  0.110   &  7.3E-03  & 11.7   &  1.7E+00  &  4.1 &  1.4E-01 \\ 
        L1L1 (B)  & 20 &  2.00  &  1.9E-02  &  0.131   &  7.3E-03  &  7.2   &  2.9E+00  &  3.5 &  4.1E-01 \\    
        L1L2 (A) & 20 &  2.00      &  1.4E-06        &    0.129      &  1.8E-08  &   7.3     &  8.6E-06  &  3.8    & 9.7E-07  \\ 
        L1L2 (B) & 20 &   2.00     &  1.2E-06        &   0.131       &  3.4E-08  &  7.2      &  1.1E-05  &    3.5  &  1.6E-06\\ 
        TLS (A) & 20 &  0.28  &  5.4E-04  &  0.110   &  1.1E-05  &  4.6   &  3.4E-03  &  2.3 &  9.4E-04 \\
        TLS (B) & 20 &  0.27  &  2.6E-04  &  0.110   &   1.1E-05 &  7.8   &  3.4E-03  &  2.2 &  9.5E-04 \\  \hline
        L1L1 (A)  & 8  &  1.28  &  1.9E-01  &  0.114   &  1.4E-02  &  7.8   &  1.6E+01  &  4.7 &  2.2E-01\\
        L1L1 (B)  & 8  &  2.00  &  1.9E-02  &  0.131   &  7.2E-03  &  7.2   &  2.9E+00  &  3.5 &  3.6E-01 \\
        L1L2 (A) & 8  &   2.00     &  7.3E-07       &    0.129      &  2.8E-08  &  7.3      &  8.5E-06      &   3.8   & 5.6E-07   \\ 
        L1L2 (B) &  8 &   2.00     &  8.1E-08       & 0.131          &  1.5E-08  &  7.2      &  2.0E-06      &    3.5  & 2.8E-07 \\
        TLS (A) & 8  &  0.79  &  5.3E-02  &  0.115   &  4.2E-03  &  3.4   &  1.3E+00  &  4.0 &  3.1E-01 \\
        TLS (B) & 8  &  1.04  &  1.4E-03  &  0.124   &  6.1E-06  &  5.0   &  7.3E-03  &  3.0 &  3.6E-04 \\
        \hline 
        \multicolumn{10}{c}{\bf Auditory}  \\
        \hline
        L1L1 (A)  & 20 &  1.66  &  7.5E-02  &  0.113   &  1.2E-03  &  26.6  &  1.0E+00  &  8.2 &  8.2E-02 \\
        L1L1 (B)  & 20 &  2.00  &  4.5E-02  &  0.165   &  5.1E-03  &  13.8  &  1.1E+00  &  6.2 &  1.5E-01 \\
        L1L2 (A) & 20 &   1.85     &  1.7E-01       &    0.123      &  9.5E-04  &  30.1    &  7.8E-01     &   7.5   & 2.3E-01  \\ 
        L1L2 (B) & 20 &    1.43    &  2.7E-02       &    0.151      &  2.0E-03 &   95 &   1.2E+01    & 5.1     &  3.1E-01\\
        TLS (A) & 20 &  0.56  &  5.6E-02  &  0.116   &  8.4E-03  &   5.6  &  1.8E+00  &  6.8 &  7.0E-01 \\
        TLS (B) & 20 &  0.54  &  3.7E-02  &  0.133   &  2.1E-03  &   5.8  &  7.2E-01  &  5.2 &  3.3E-01 \\ \hline
        L1L1 (A)  & 8  &  1.32  &  1.2E-01  &  0.114   &  2.5E-02  &   6.7  &  1.5E+01  &  9.1 &  9.0E-01 \\
        L1L1 (B)  & 8  &  2.00  &  3.4E-02  &  0.167   &  2.9E-03  &   9.8  &  1.1E+01 &  4.8 &  5.7E-01 \\
        L1L2 (A) & 8  &   1.19     &  3.8E-02       &   0.113       &  3.1E-03  &  19.1     &  1.3E+01      &  7.6   &  9.4E-01 \\ 
        L1L2 (B) & 8  &  1.73      &  4.9E-02      &  0.163         &  1.5E-03  & 17.8        &  7.6E+00      &    5.9  &   7.6E-01 \\ 
        TLS (A) & 8  &  0.84  &  6.4E-02  &  0.126   &  1.9E-02  &   6.6  &  1.0E+01  &  7.8 &  6.5E-01 \\
        TLS (B) & 8  &  1.00  &  5.2E-02  &  0.145   &  9.3E-04 &  10.4  &  1.2E+00  &  5.6 &  1.2E-01 \\ 
        \hline
        \multicolumn{10}{c}{\bf Visual}  \\
        \hline
        L1L1 (A)  & 20 &  0.98  &  1.2E-02  &  0.113   &  3.2E-03  &  16.0  &  1.0E+00  &  7.3 &  2.8E-01 \\
        L1L1 (B)  & 20 &  1.00  &  2.8E-02  &  0.149   &  8.9E-03 &  26.4  &  3.0E+00  &  2.7 &  1.2E-01 \\
        L1L2 (A) & 20 &  0.83      &  1.1E-07       &   0.131       &   2.7E-06 &  10.5     & 9.0E-02     &  4.5  &  1.2E-02 \\ 
        L1L2 (B) & 20 &  1.06      &  4.4E-03    &    0.149      &   1.1E-03 & 27.4      & 1.3E+01       &  2.8    &   1.1E-01 \\ 
        TLS (A) & 20 &  0.87  &  6.0E-02  &  0.113   &  8.7E-03  &  13.8  &  4.2E+00  &  6.5 &  3.3E-01 \\
        TLS (B) & 20 &  0.60  &  3.1E-02  &  0.127   &  1.6E-03  &  12.4  &  3.7E-01  &  3.0 &  3.1E-01 \\ \hline
        L1L1 (A)  & 8  &  1.19  &  3.5E-02  &  0.113   &  1.1E-02  &  15.8  &  1.2E+01  &  7.5 &  7.0E-01 \\
        L1L1 (B)  & 8  &  1.18  &  4.8E-02  &  0.151   &  1.7E-03 &  31.3  &  6.6E+00  &  2.8 &  6.9E-02 \\
        L1L2 (A) & 8  &   1.19     &  3.8E-02        &    0.112      &  6.3E-03 &  23.1     &    1.8E+00   &  5.4    &  3.2E-01 \\ 
        L1L2 (B) & 8  &   1.20     &  3.8E-02        &   0.151       &  1.1E-03 &    29.6    &    2.7E+00    &  3.0    &  1.0E-08 \\ 
        TLS (A) & 8  &  0.97  &  8.0E-02  &  0.111   &  4.0E-03  &  18.7  &  2.4E+00  &  7.2 &  5.1E-01 \\
        TLS (B) & 8  &  0.81  &  2.7E-02  &  0.134   &  2.4E-03  &  17.9  &  8.4E-01  &  4.0 &  3.1E-01     \\ 
        \hline  
        \end{tabular}
   \end{scriptsize}
\label{table:method_comparison}
\end{table*}

\subsection{Results of the metaheuristic search}
\label{sec:Lattice_Search_Results}

\begin{figure}[h!]
    \centering
    \begin{scriptsize}
        \begin{minipage}{8cm}
            \centering
            \hrule
            \vskip0.1cm
            \begin{minipage}{0.2cm}
                \mbox{}
            \end{minipage}
            \begin{minipage}{2.4cm}
                \centering 
                \textbf{Somatosensory}
            \end{minipage}
            \begin{minipage}{2.4cm}
                \centering
                \textbf{Auditory}
            \end{minipage}
            \begin{minipage}{2.4cm}
                \centering
                \textbf{Visual}
            \end{minipage}
            \vskip0.1cm
            \hrule 
            \vskip0.1cm
            \begin{minipage}{0.2cm}
                \rotatebox{90}{$\Gamma$ (A/m\textsuperscript{2})}
            \end{minipage}
            \begin{minipage}{2.4cm}
                \centering 
                \includegraphics[height=2.0cm]{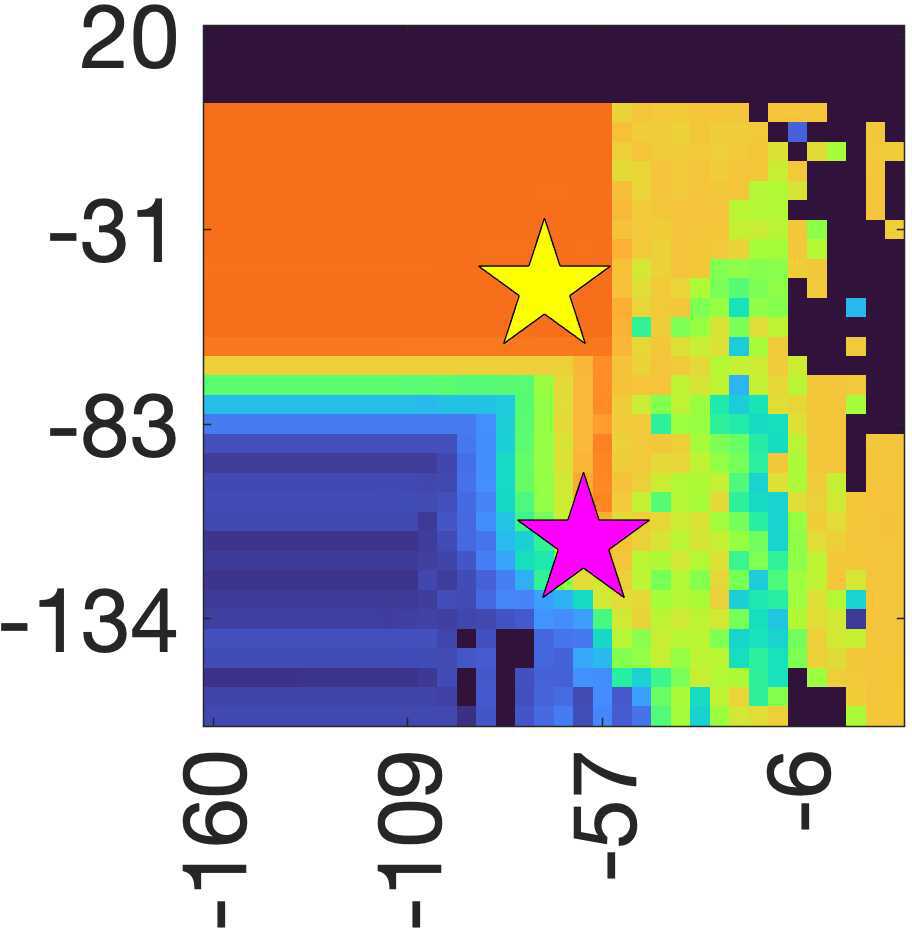}
            \end{minipage}
            \begin{minipage}{2.4cm}
                \centering
               \includegraphics[height=2.0cm]{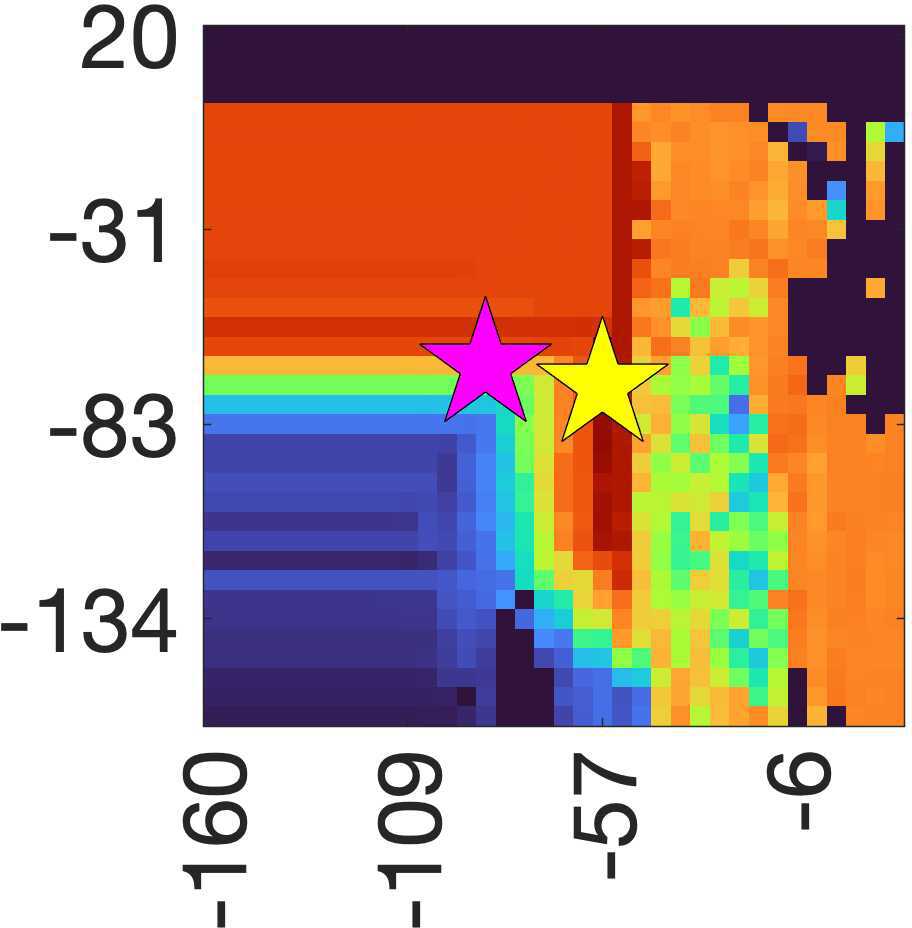}
            \end{minipage}
            \begin{minipage}{2.4cm}
                \centering
                \includegraphics[height=2.0cm]{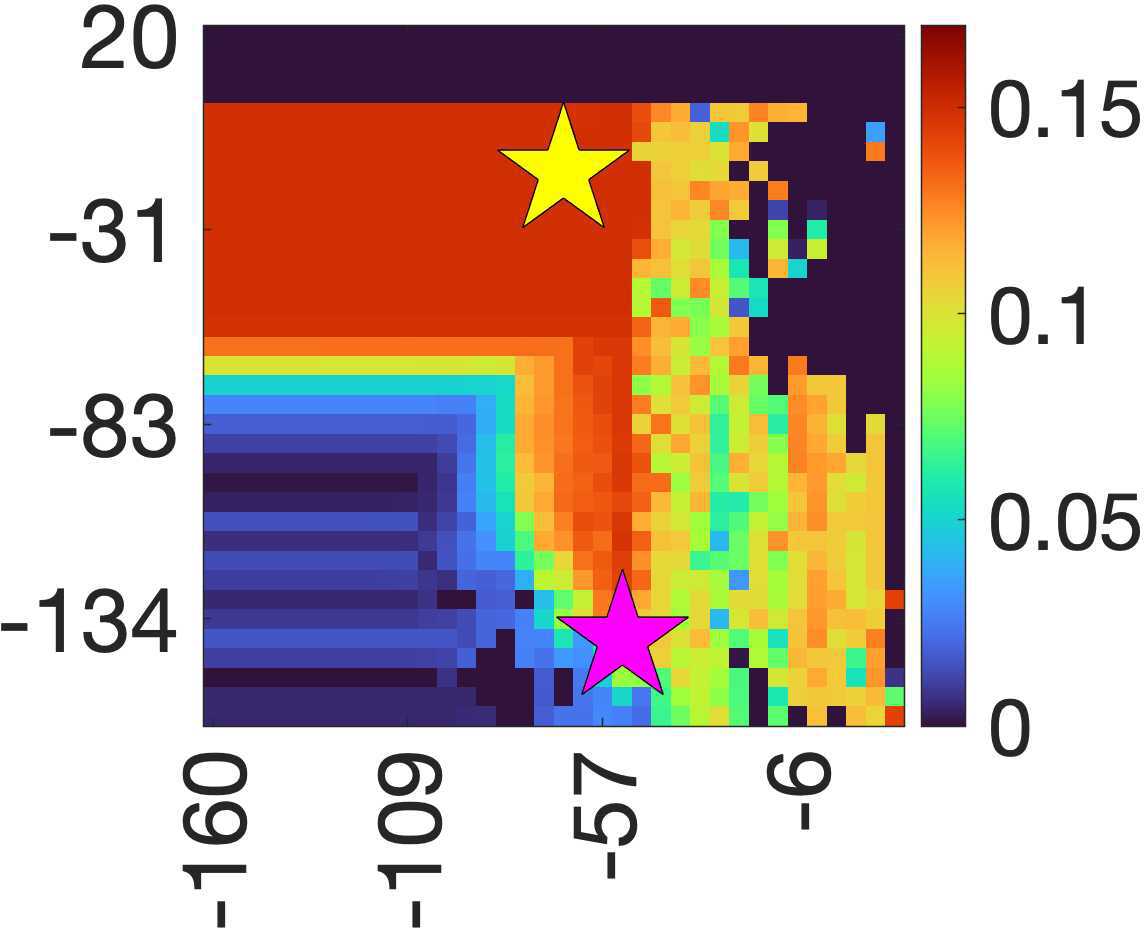}
            \end{minipage}
            \vskip0.1cm
            \begin{minipage}{0.2cm}
                \rotatebox{90}{$\Theta$}
            \end{minipage}
            \begin{minipage}{2.4cm}
                \centering 
               \includegraphics[height=2.0cm]{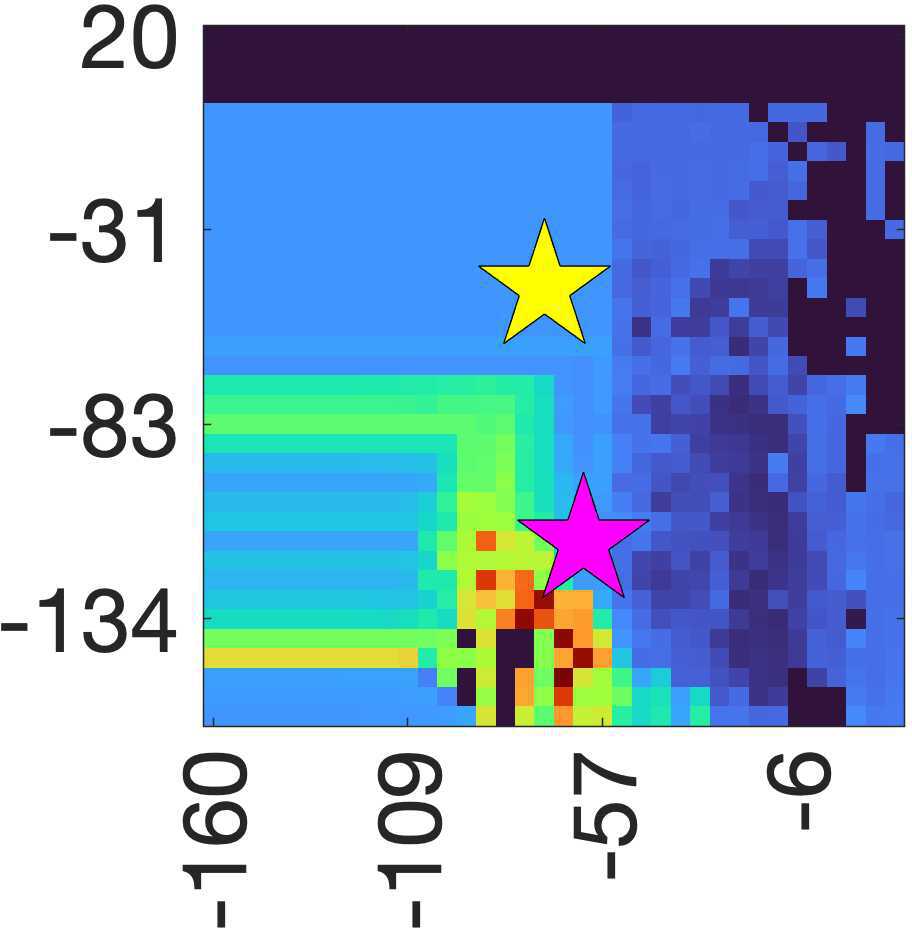}
            \end{minipage}
            \begin{minipage}{2.4cm}
                \centering
               \includegraphics[height=2.0cm]{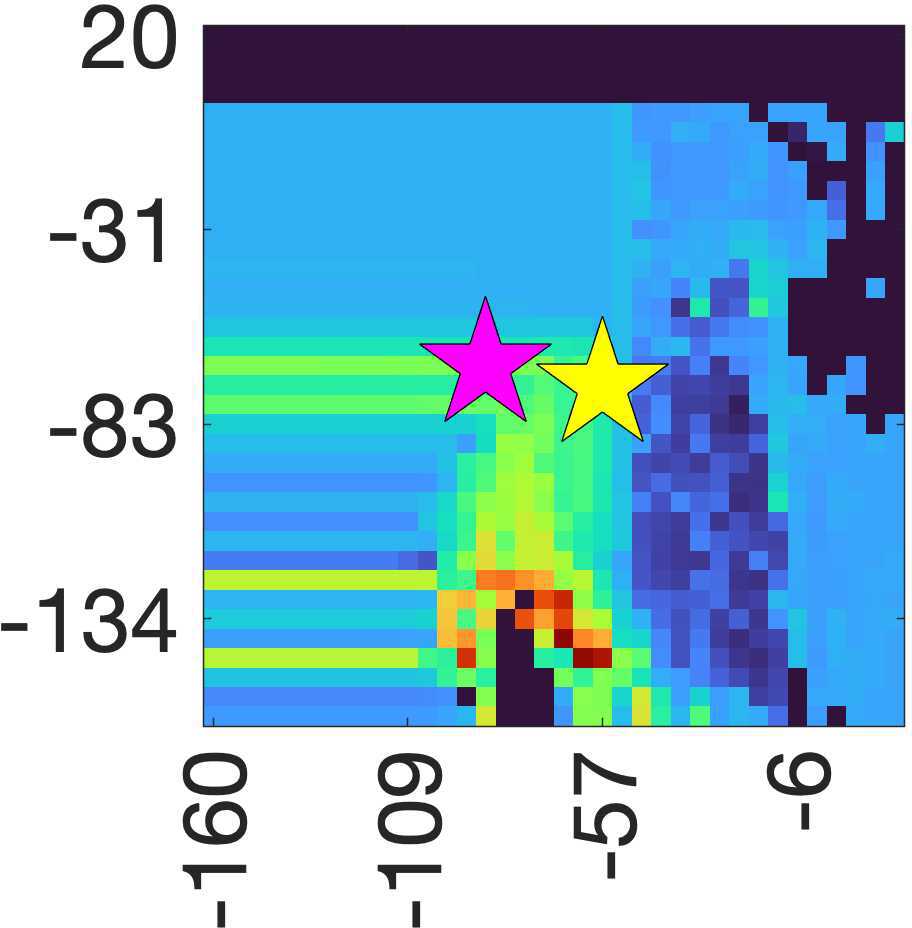}
            \end{minipage}
            \begin{minipage}{2.4cm}
                \centering
                \includegraphics[height=2.0cm]{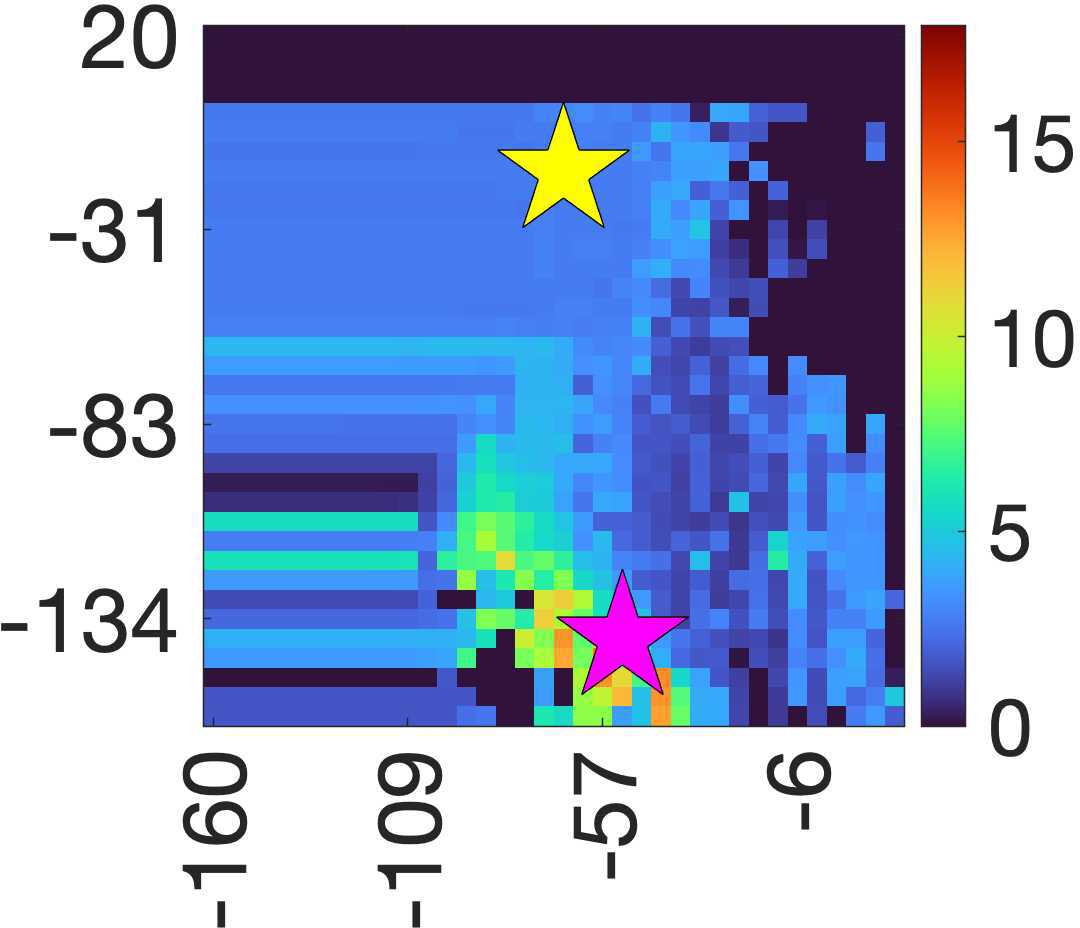}
            \end{minipage}
            \vskip0.1cm
            \begin{minipage}{0.2cm}
                \rotatebox{90}{AD (deg)}
            \end{minipage}
            \begin{minipage}{2.4cm}
                \centering 
                \includegraphics[height=2.0cm]{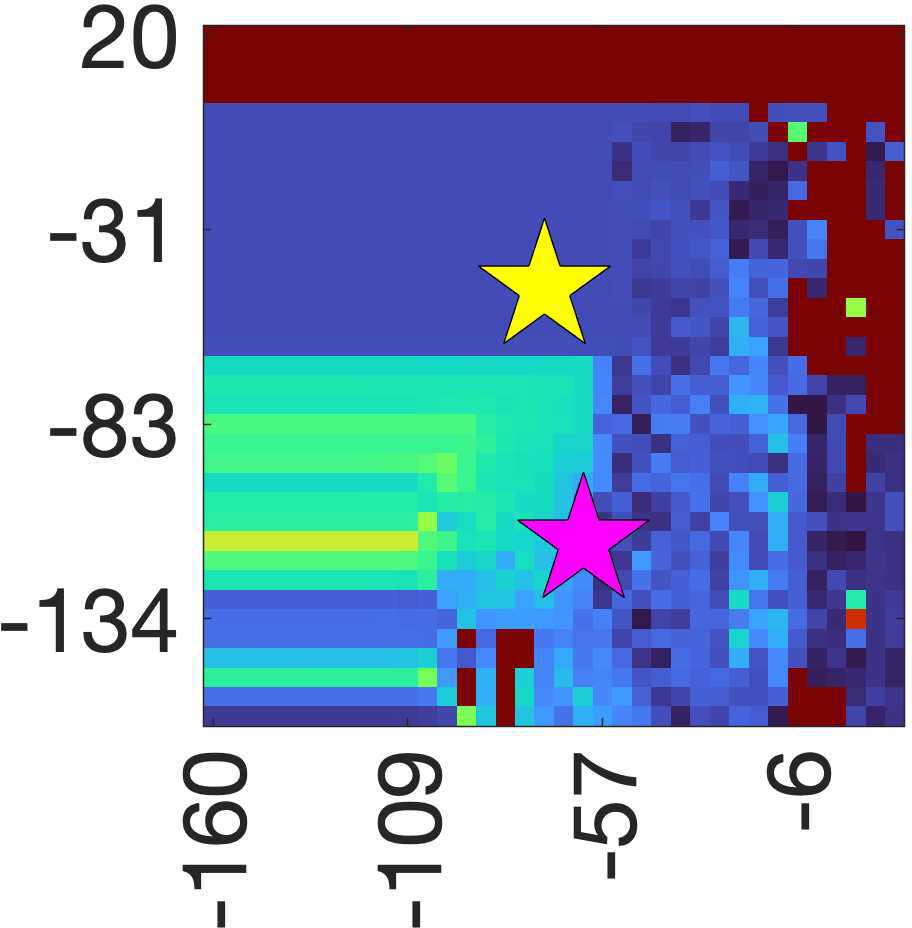}
            \end{minipage}
            \begin{minipage}{2.4cm}
                \centering
                \includegraphics[height=2.0cm]{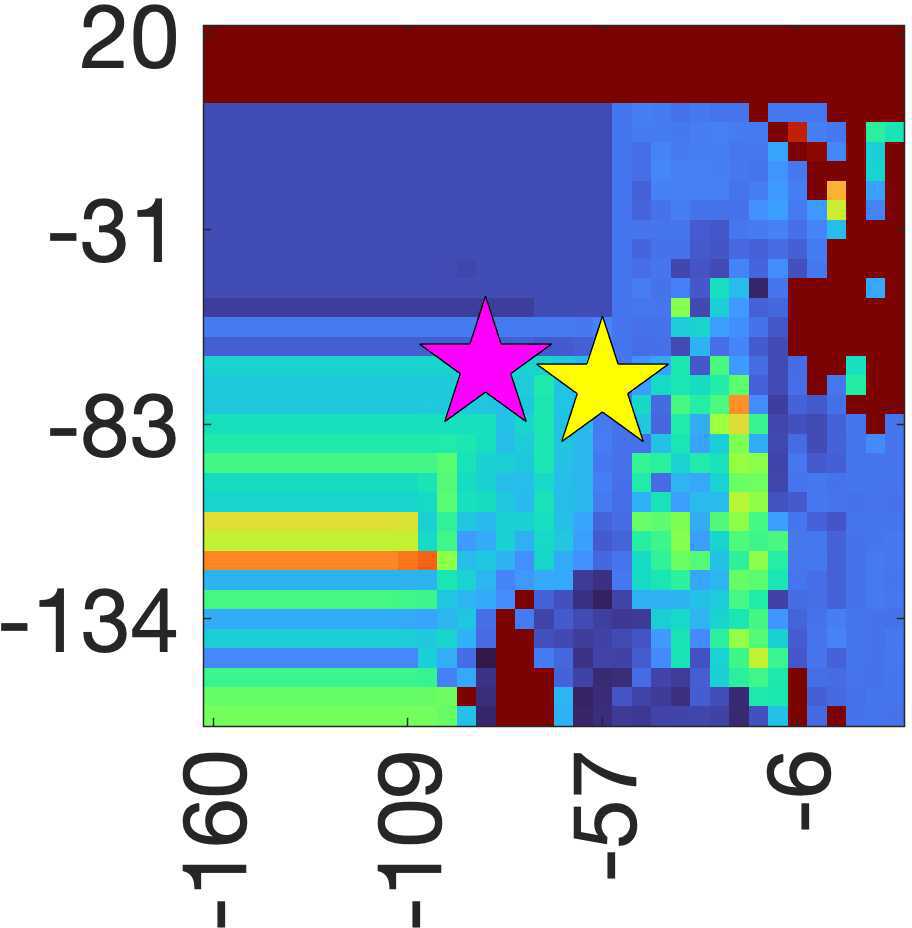}
            \end{minipage}
            \begin{minipage}{2.4cm}
                \centering
                \includegraphics[height=2.0cm]{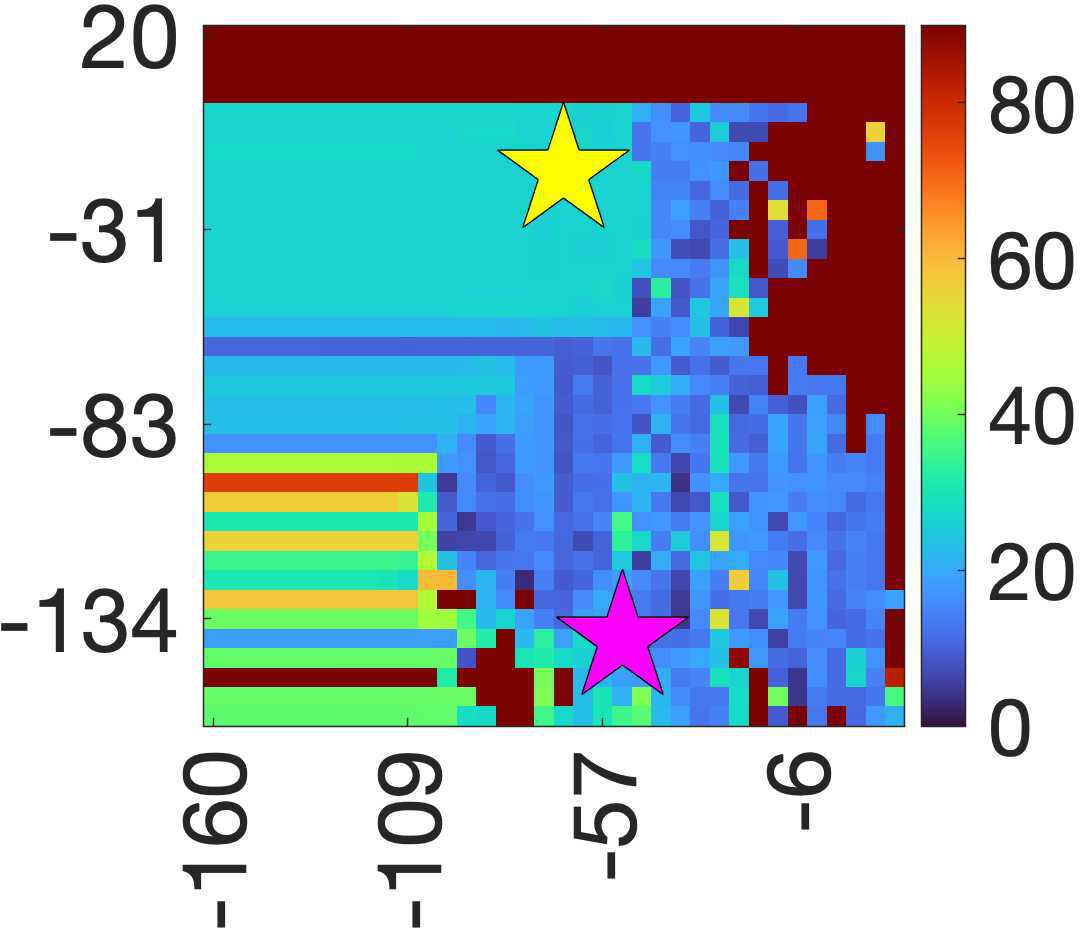}
            \end{minipage}
            \vskip0.1cm
            \begin{minipage}{0.2cm}
                \rotatebox{90}{$\|{\bf y}\|_\infty$ (mA)}
            \end{minipage} 
            \begin{minipage}{2.4cm}
                \centering 
                \includegraphics[height=2.0cm]{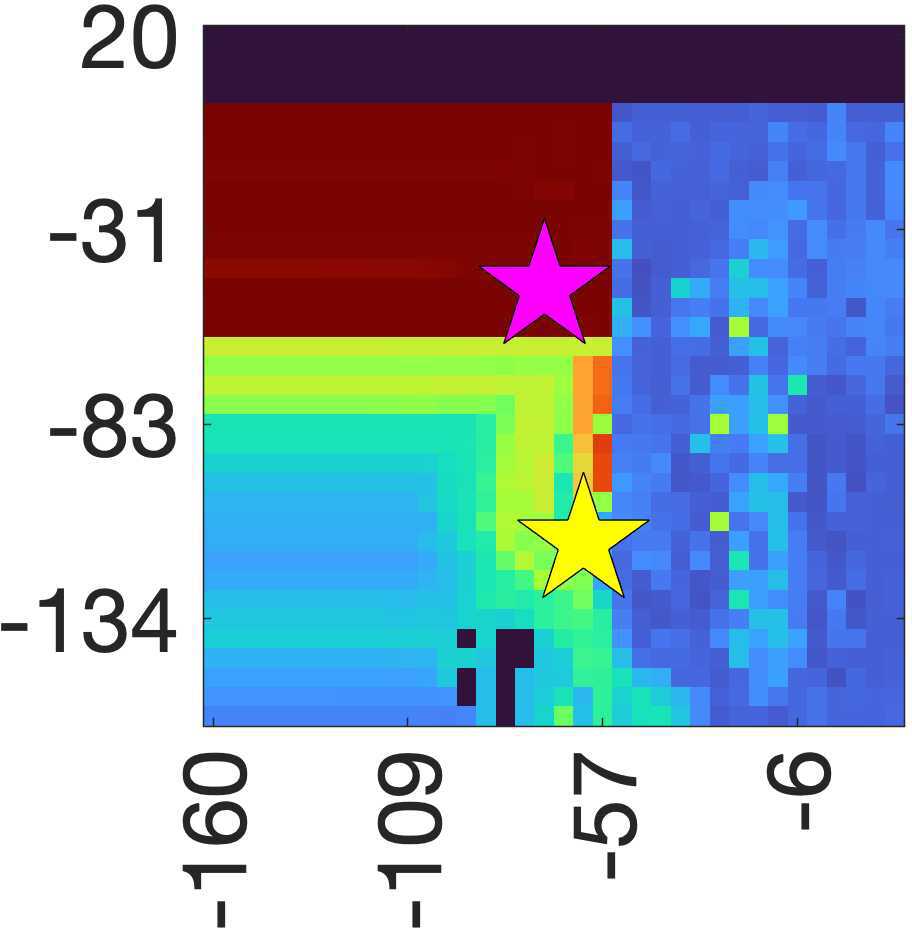}
            \end{minipage}
            \begin{minipage}{2.4cm}
                \centering
                \includegraphics[height=2.0cm]{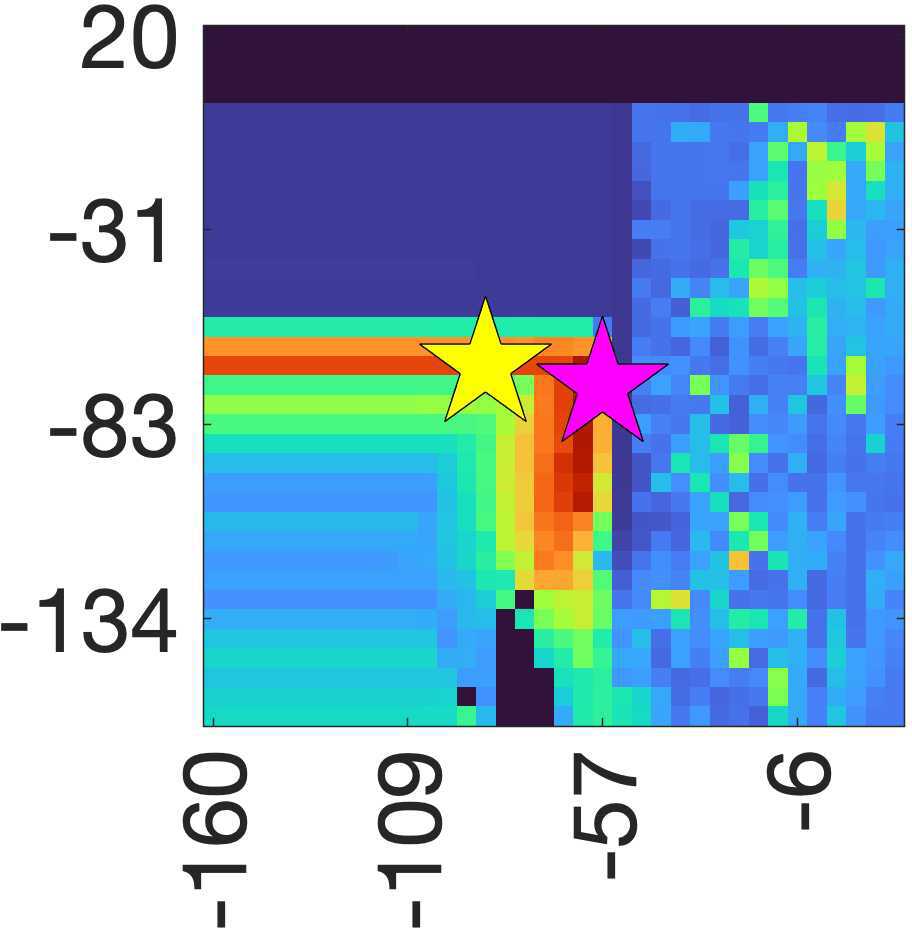}
            \end{minipage}
            \begin{minipage}{2.4cm}
                \centering
                \includegraphics[height=2.0cm]{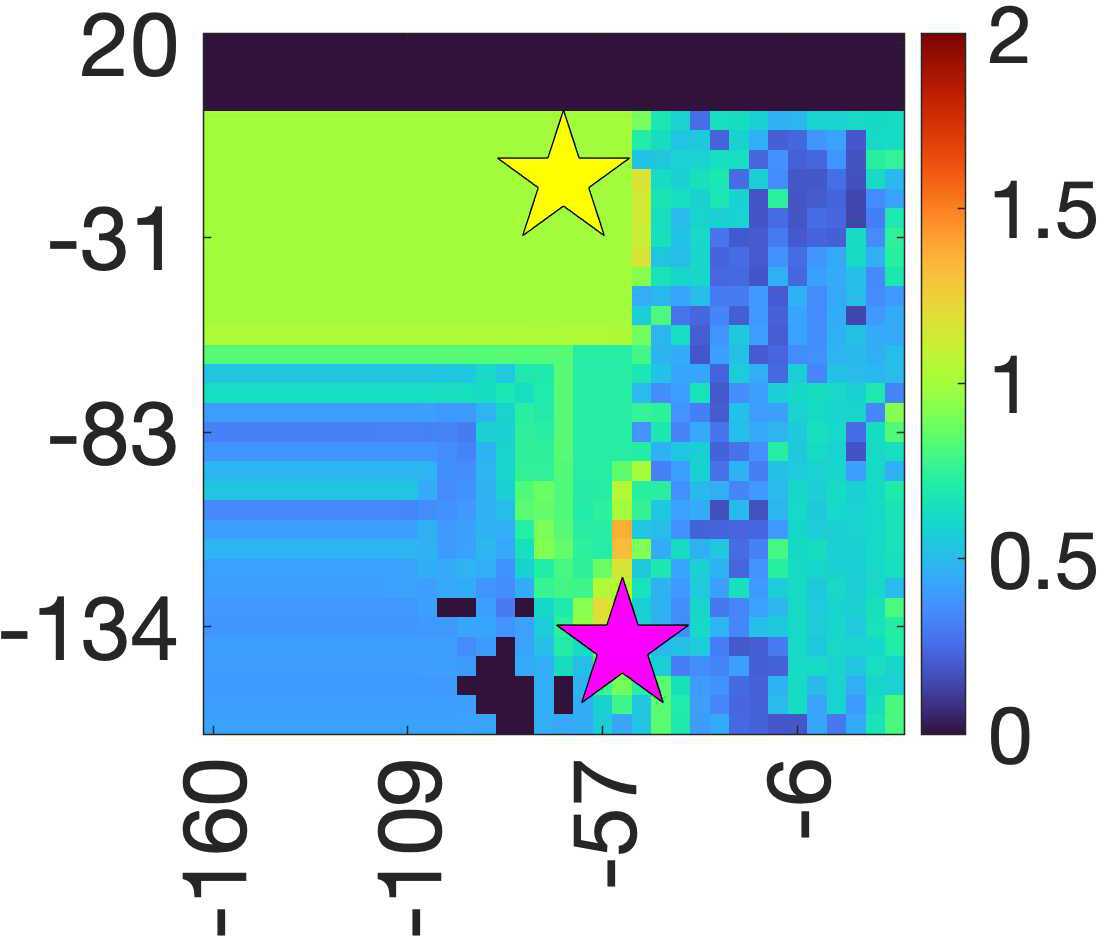}
            \end{minipage}
            \vskip0.1cm
            \hrule
        \end{minipage}
    \end{scriptsize}
\caption{L1-norm regularized L1-norm fitting (L1L1) charts showing $\Gamma$ (A/m\textsuperscript{2}), $\Theta$ (relative), AD (deg), and $\| {\bf y} \|_\infty$ (mA) for each point of the lattice search in the first run of the two-stage metaheuristic optimization process. Horizontal axis corresponds to the regularization parameter $\alpha$ and the vertical to the nuisance field weight $\varepsilon$. The optimal solution for the case (A) is represented by a purple star, while for the case (B) as a yellow star. Axis are decibel (dB) scaled.
\label{fig:imagesc_L1L1}}
\end{figure}

\begin{figure}[h!]
    \centering
    \begin{scriptsize}
        \begin{minipage}{8cm}
            \centering
            \hrule
            \vskip0.1cm
            \begin{minipage}{0.2cm}
                \mbox{}
            \end{minipage}
            \begin{minipage}{2.4cm}
                \centering 
                \textbf{Somatosensory}
            \end{minipage}
            \begin{minipage}{2.4cm}
                \centering
                \textbf{Auditory}
            \end{minipage}
            \begin{minipage}{2.4cm}
                \centering
                \textbf{Visual}
            \end{minipage}
            \vskip0.1cm
            \hrule
            \vskip0.1cm
            \begin{minipage}{0.2cm}
                \rotatebox{90}{$\Gamma$ (A/m\textsuperscript{2})}
            \end{minipage}
            \begin{minipage}{2.4cm}
                \centering 
                \includegraphics[height=2.0cm]{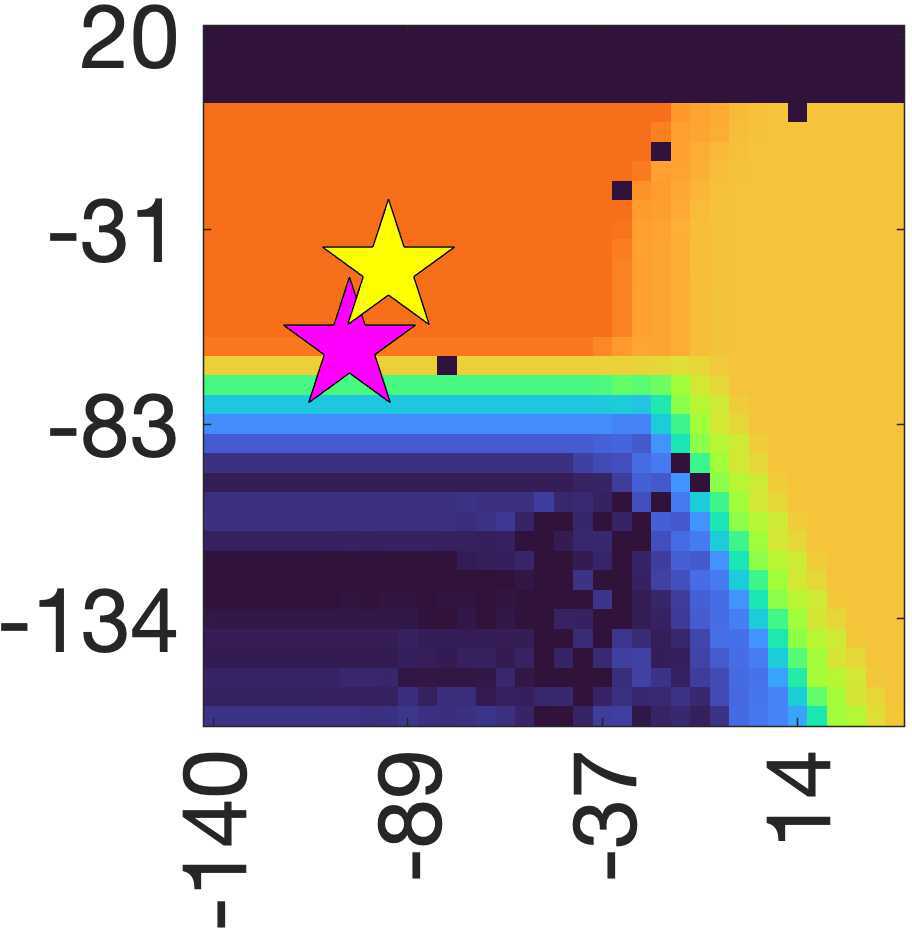}
            \end{minipage}
            \begin{minipage}{2.4cm}
                \centering
                \includegraphics[height=2.0cm]{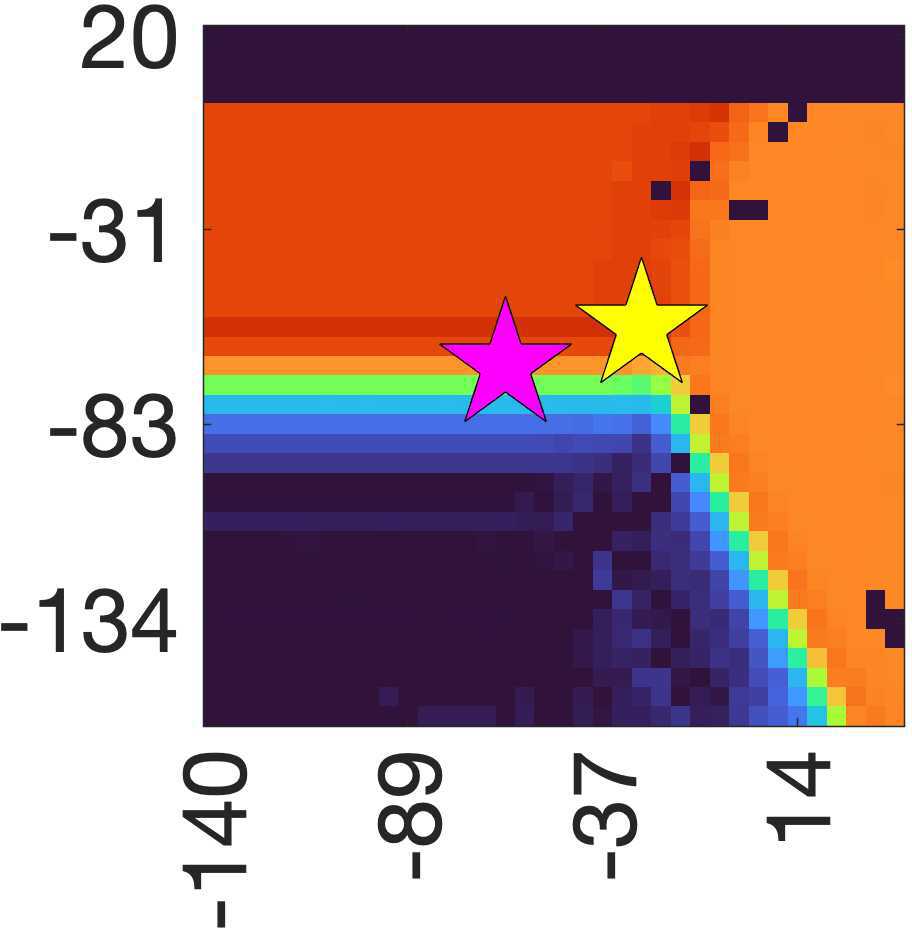}
            \end{minipage}
            \begin{minipage}{2.4cm}
                \centering
                \includegraphics[height=2.0cm]{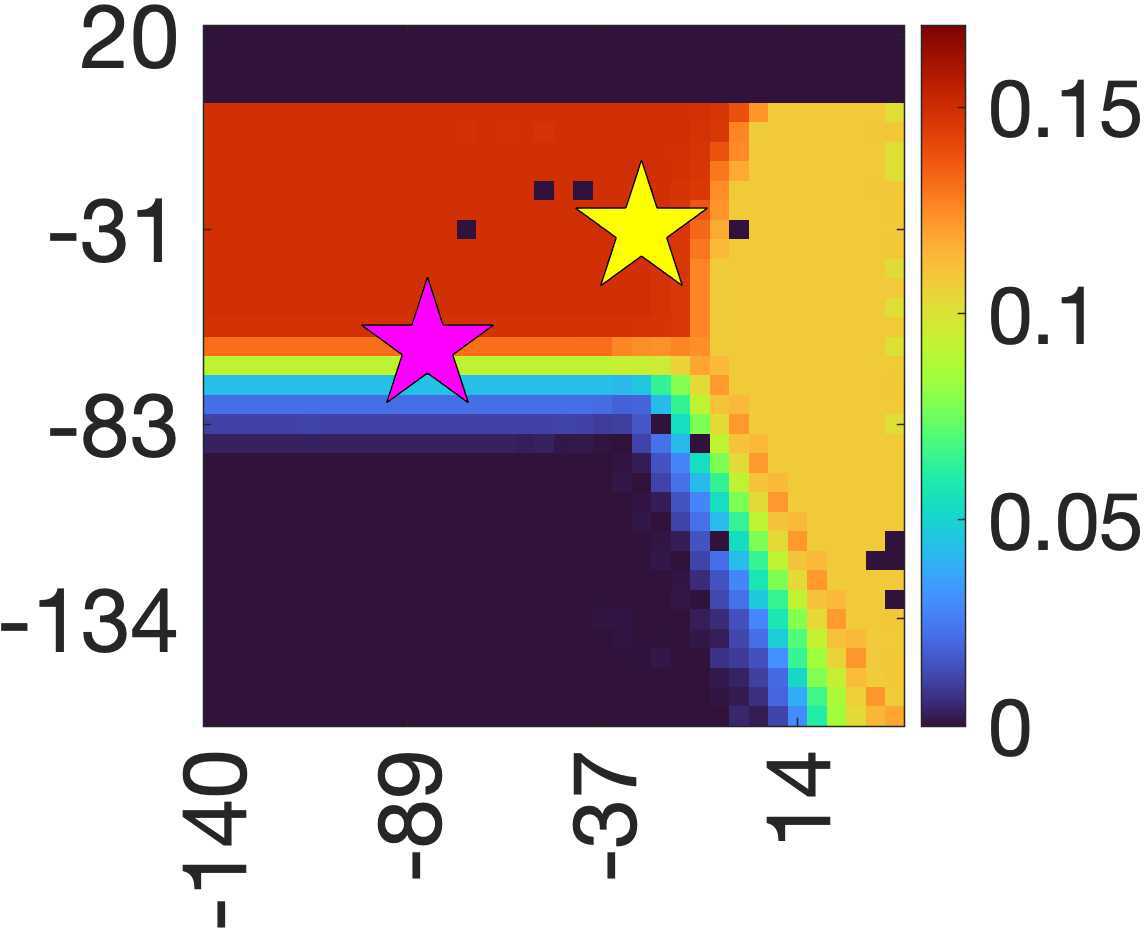}
            \end{minipage}
            \vskip0.1cm
            \begin{minipage}{0.2cm}
                \rotatebox{90}{$\Theta$}
            \end{minipage}
            \begin{minipage}{2.4cm}
                \centering 
                \includegraphics[height=2.0cm]{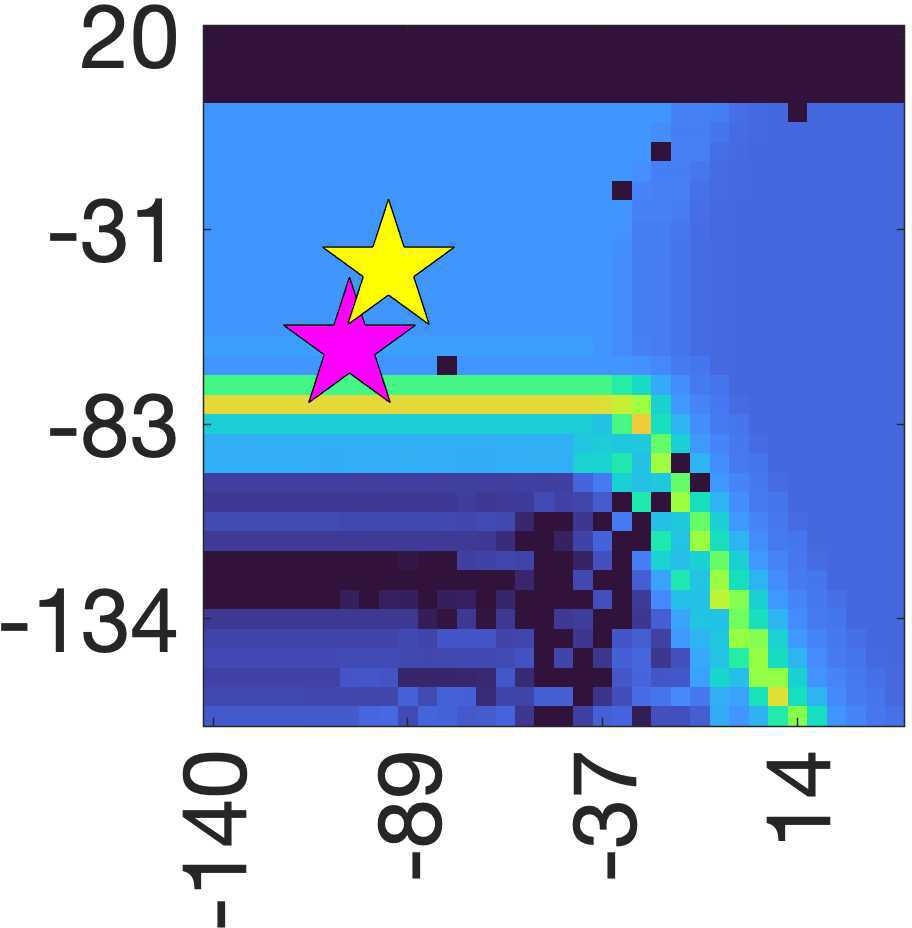}
            \end{minipage}
            \begin{minipage}{2.4cm}
                \centering
                 \includegraphics[height=2.0cm]{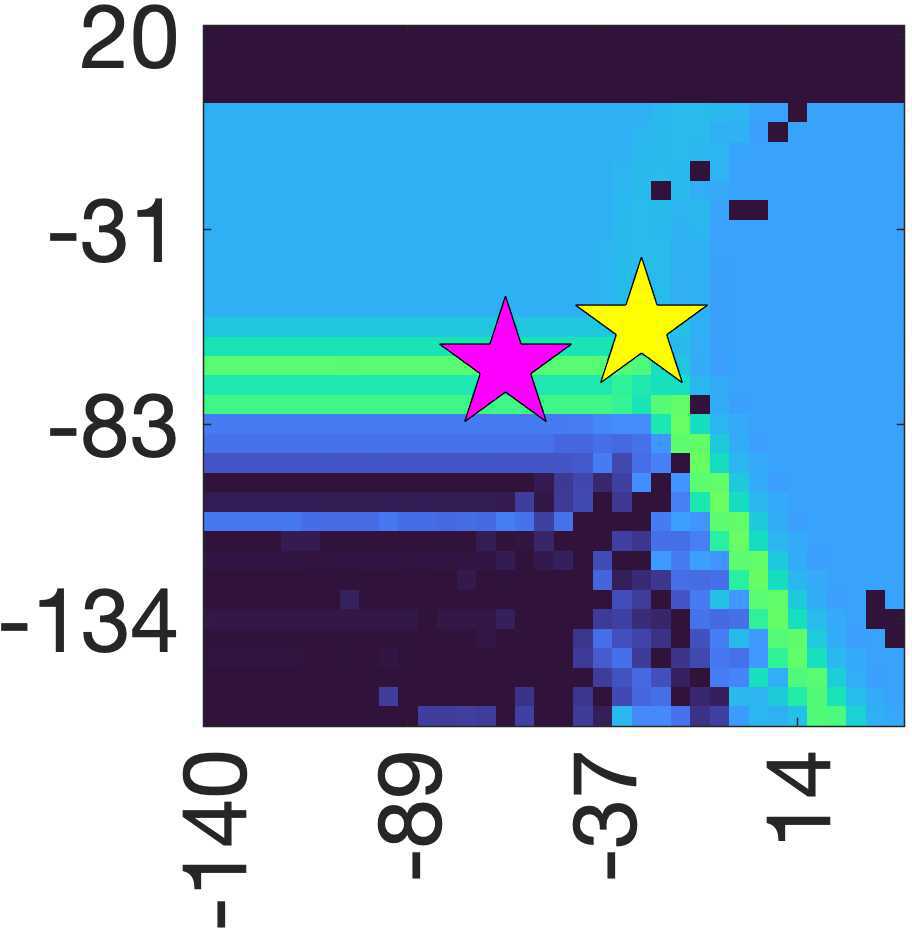}
            \end{minipage}
            \begin{minipage}{2.4cm}
                \centering
                \includegraphics[height=2.0cm]{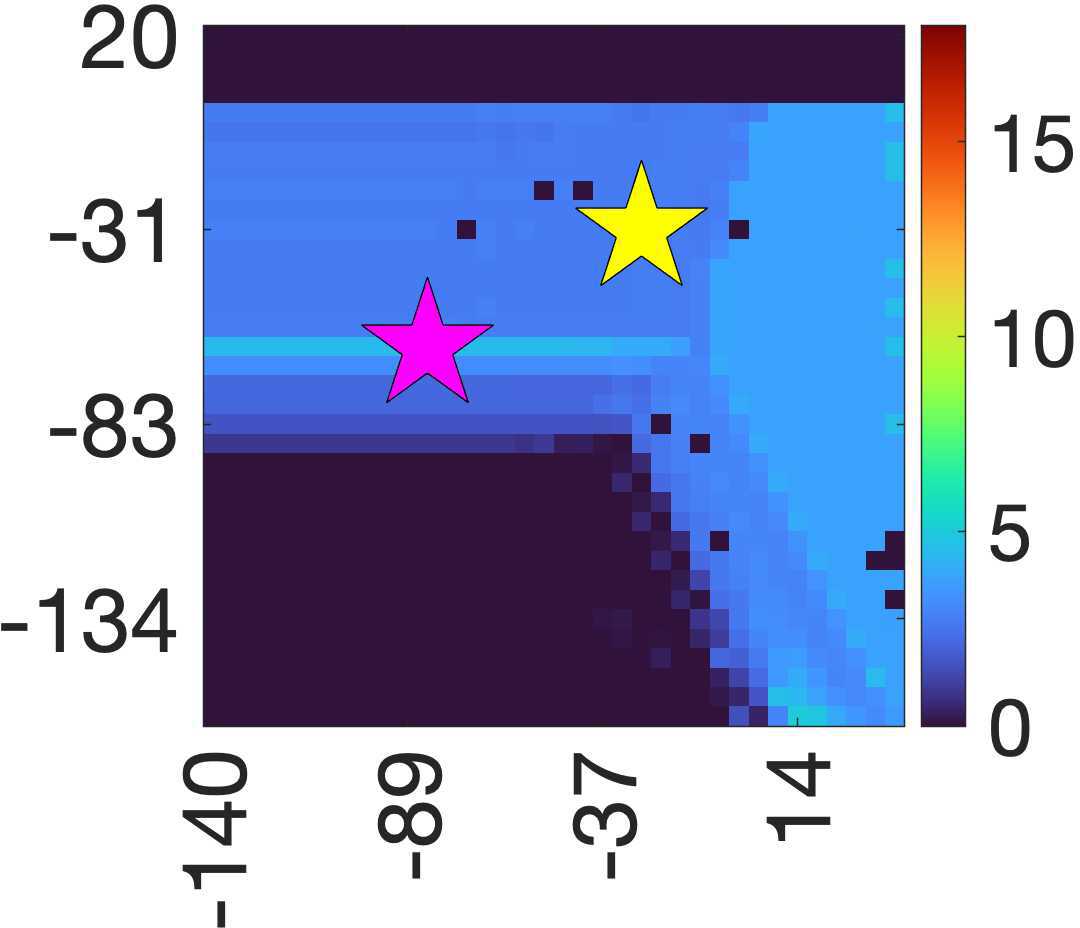}
            \end{minipage}
            \vskip0.1cm
            \begin{minipage}{0.2cm}
                \rotatebox{90}{AD (deg)}
            \end{minipage}
            \begin{minipage}{2.4cm}
                \centering 
                \includegraphics[height=2.0cm]{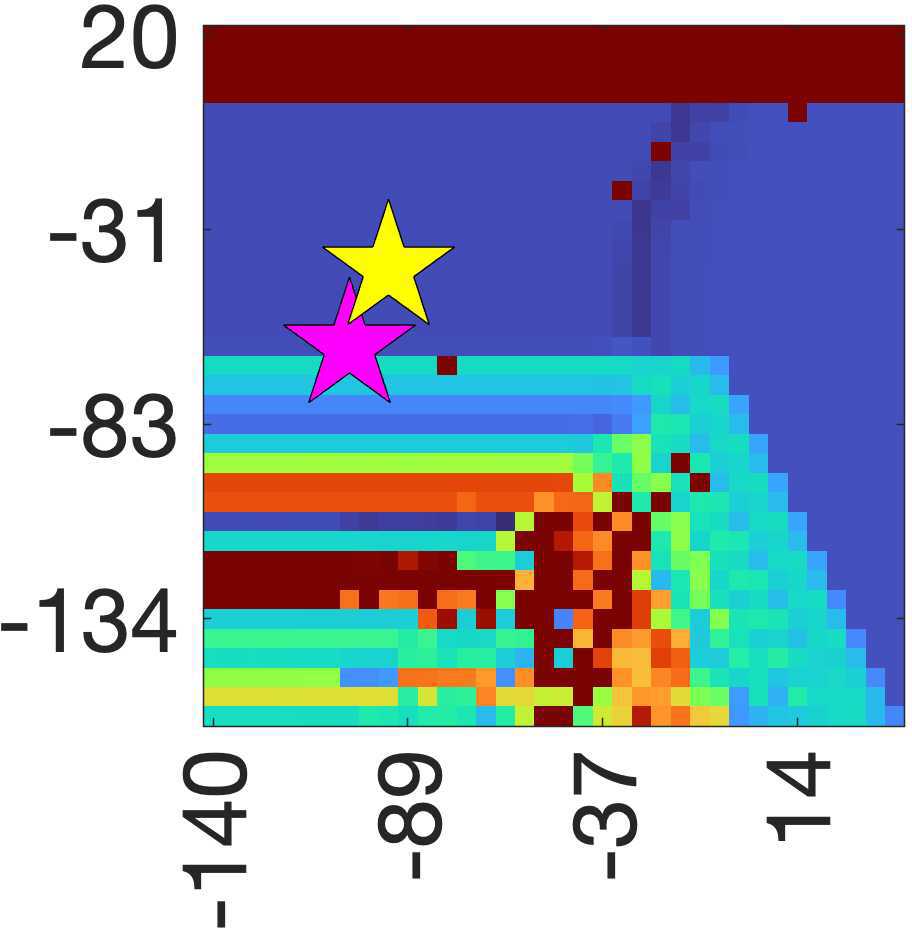}
            \end{minipage}
            \begin{minipage}{2.4cm}
                \centering
                \includegraphics[height=2.0cm]{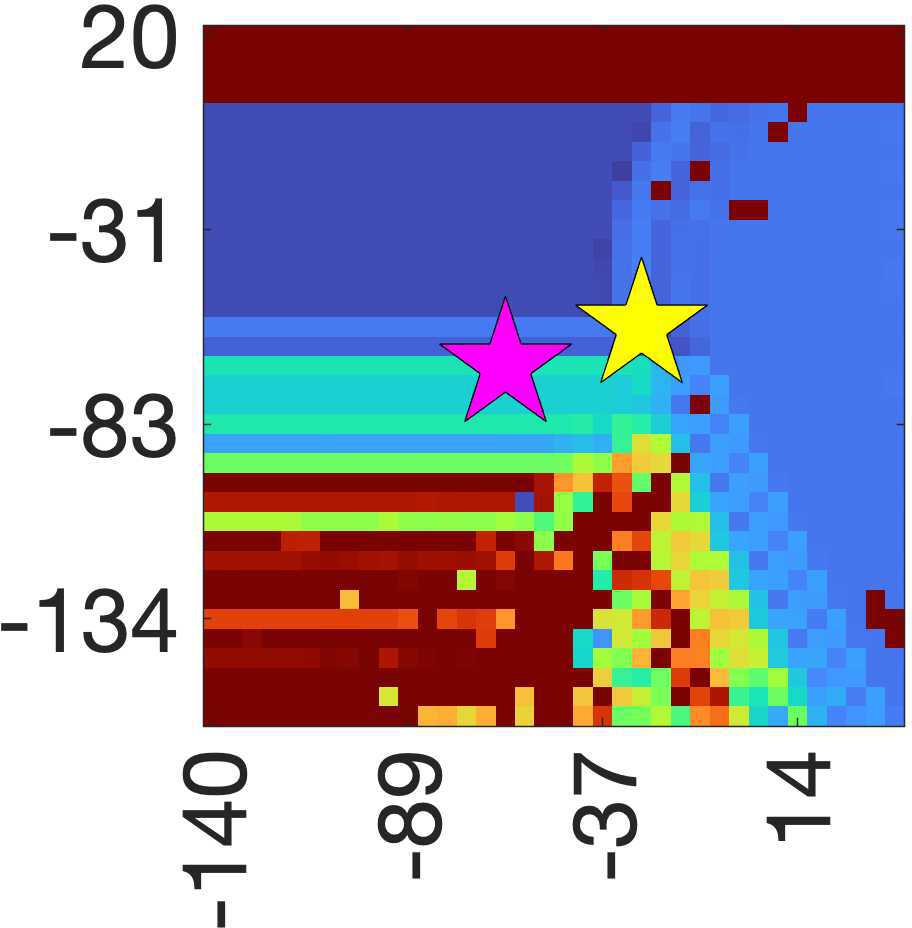}
            \end{minipage}
            \begin{minipage}{2.4cm}
                \centering
               \includegraphics[height=2.0cm]{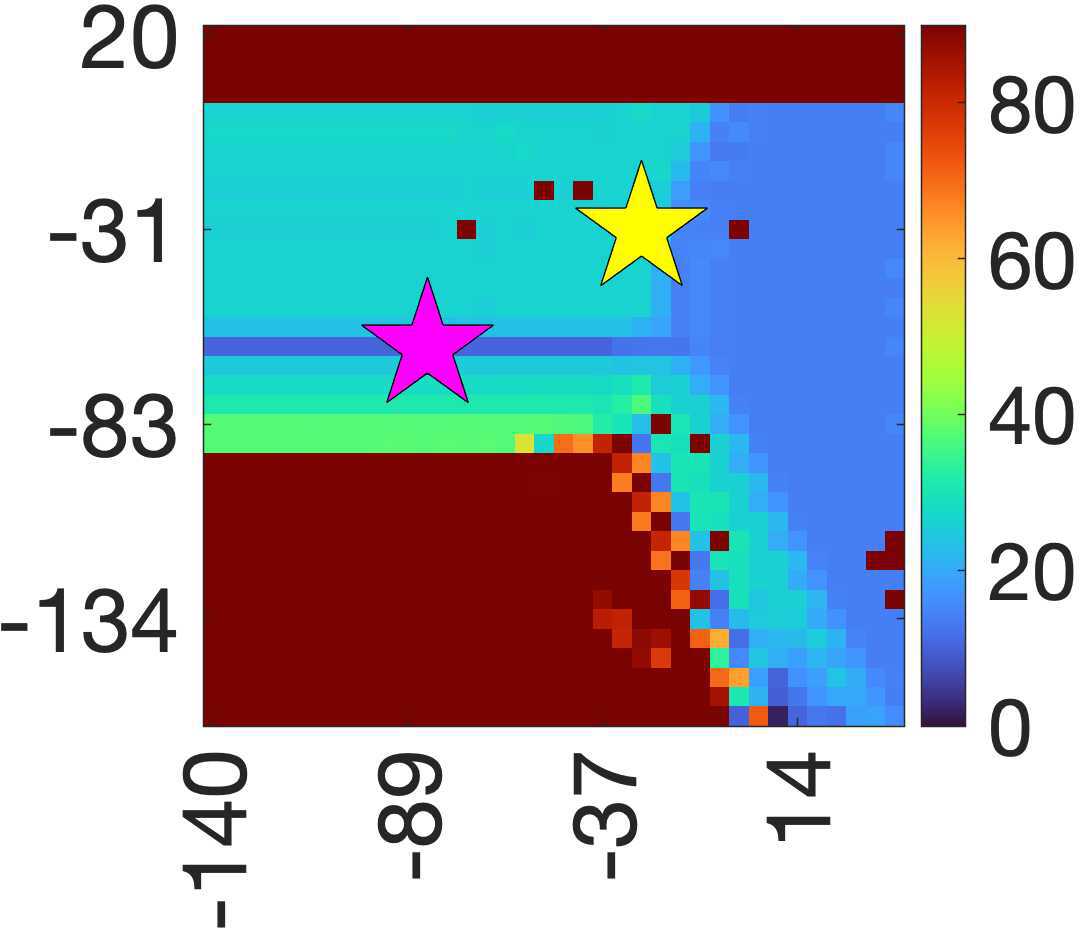}
            \end{minipage}
            \vskip0.1cm
            \begin{minipage}{0.2cm}
                \rotatebox{90}{$\| {\bf y}  \|_\infty$ (mA)}
            \end{minipage}
            \begin{minipage}{2.4cm}
                \centering 
                \includegraphics[height=2.0cm]{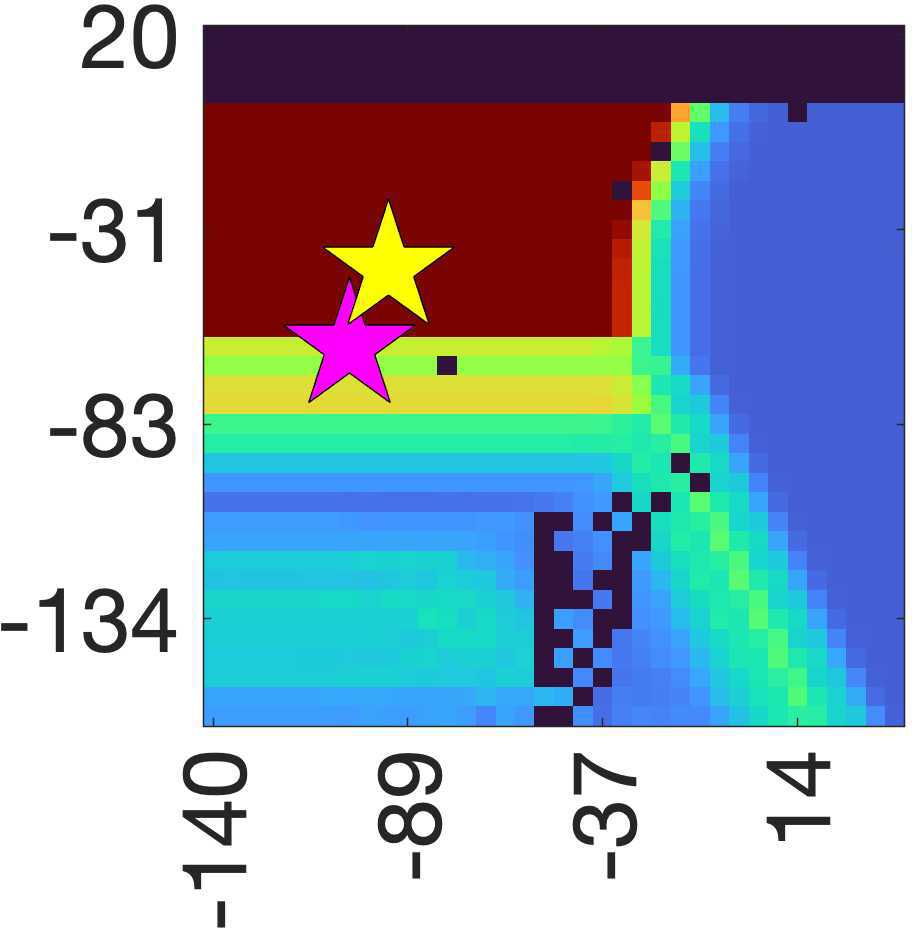}
            \end{minipage}
            \begin{minipage}{2.4cm}
                \centering
                 \includegraphics[height=2.0cm]{images_jpg/revision/auditory_SDP_angle.jpg}
            \end{minipage}
            \begin{minipage}{2.4cm}
                \centering
                \includegraphics[height=2.0cm]{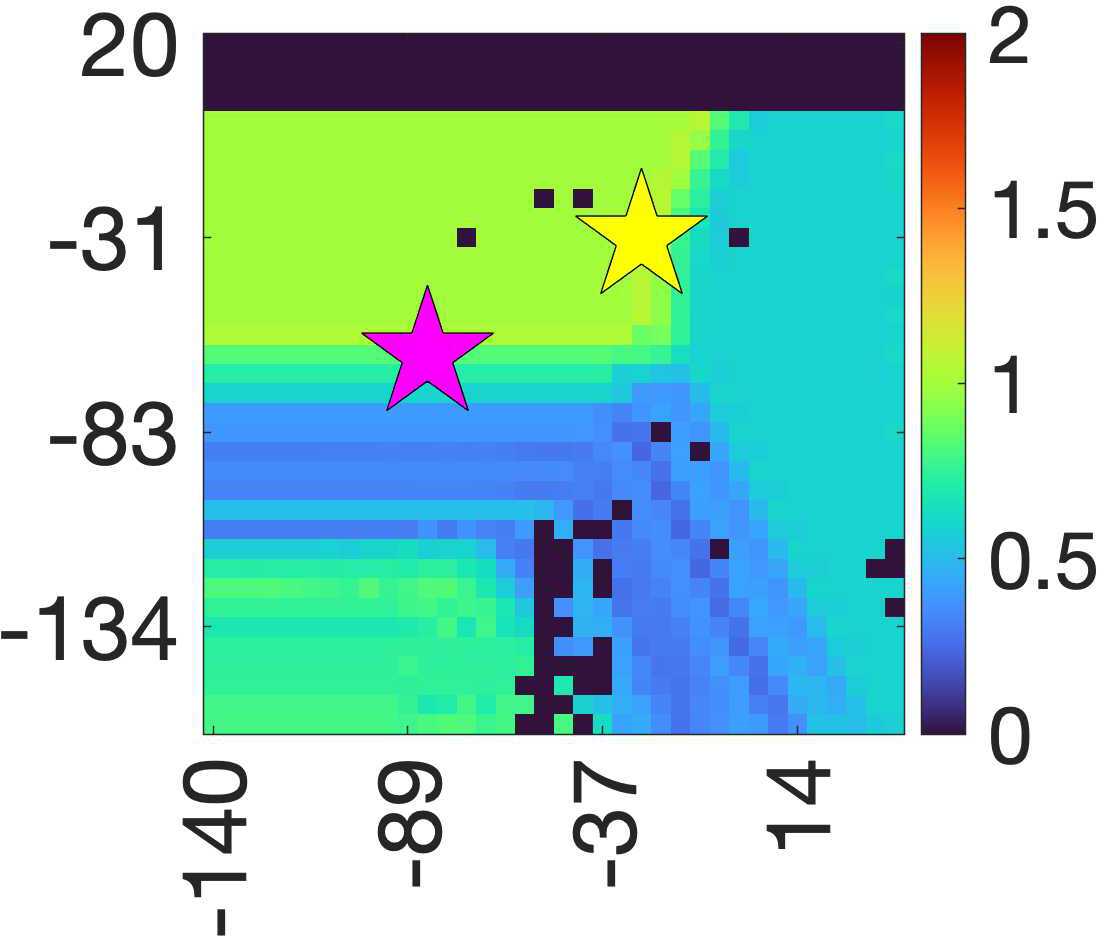}
            \end{minipage}
            \vskip0.1cm
            \hrule
        \end{minipage}
    \end{scriptsize}
\caption{L1-norm regularized L2-norm fitting (L1L2) charts showing $\Gamma$ (A/m\textsuperscript{2}), $\Theta$ (relative), AD (deg), and $\| {\bf y} \|_\infty$ (mA) for each point of the lattice search in the first run of the two-stage metaheuristic optimization process. Horizontal axis corresponds to the regularization parameter $\alpha$ and the vertical to the nuisance field weight $\varepsilon$. The optimal solution for the case (A) is represented by a purple star, while for the case (B) as a yellow star. Axis are decibel (dB) scaled.
\label{fig:imagesc_L1L2}}
\end{figure}

\begin{figure}[h!]
    \centering
    \begin{scriptsize}
        \begin{minipage}{8cm}
            \centering
            \hrule
            \vskip0.1cm
            \begin{minipage}{0.2cm}
                \mbox{}
            \end{minipage}
            \begin{minipage}{2.4cm}
                \centering 
                \textbf{Somatosensory}
            \end{minipage}
            \begin{minipage}{2.4cm}
                \centering
                \textbf{Auditory}
            \end{minipage}
            \begin{minipage}{2.4cm}
                \centering
                \textbf{Visual}
            \end{minipage}
            \vskip0.1cm
            \hrule
            \vskip0.1cm
            \begin{minipage}{0.2cm}
                \rotatebox{90}{$\Gamma$ (A/m\textsuperscript{2})}
            \end{minipage}
            \begin{minipage}{2.4cm}
                \centering 
                \includegraphics[height=2.0cm]{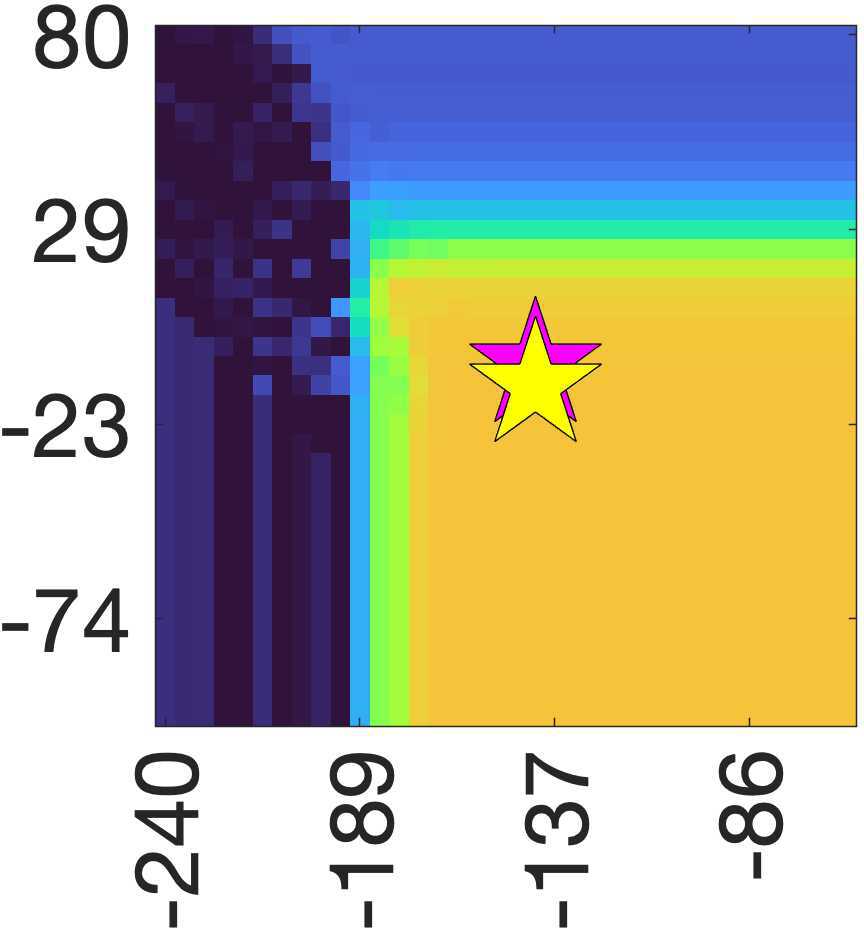}
            \end{minipage}
            \begin{minipage}{2.4cm}
                \centering
                \includegraphics[height=2.0cm]{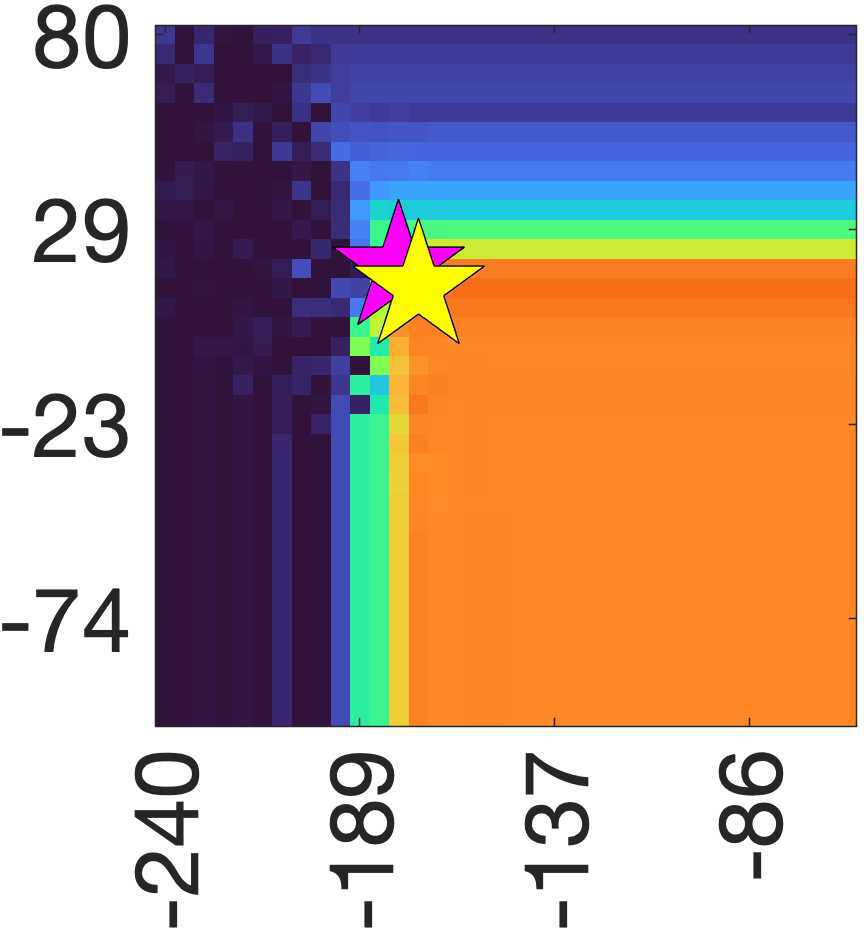}
            \end{minipage}
            \begin{minipage}{2.4cm}
                \centering
                \includegraphics[height=2.0cm]{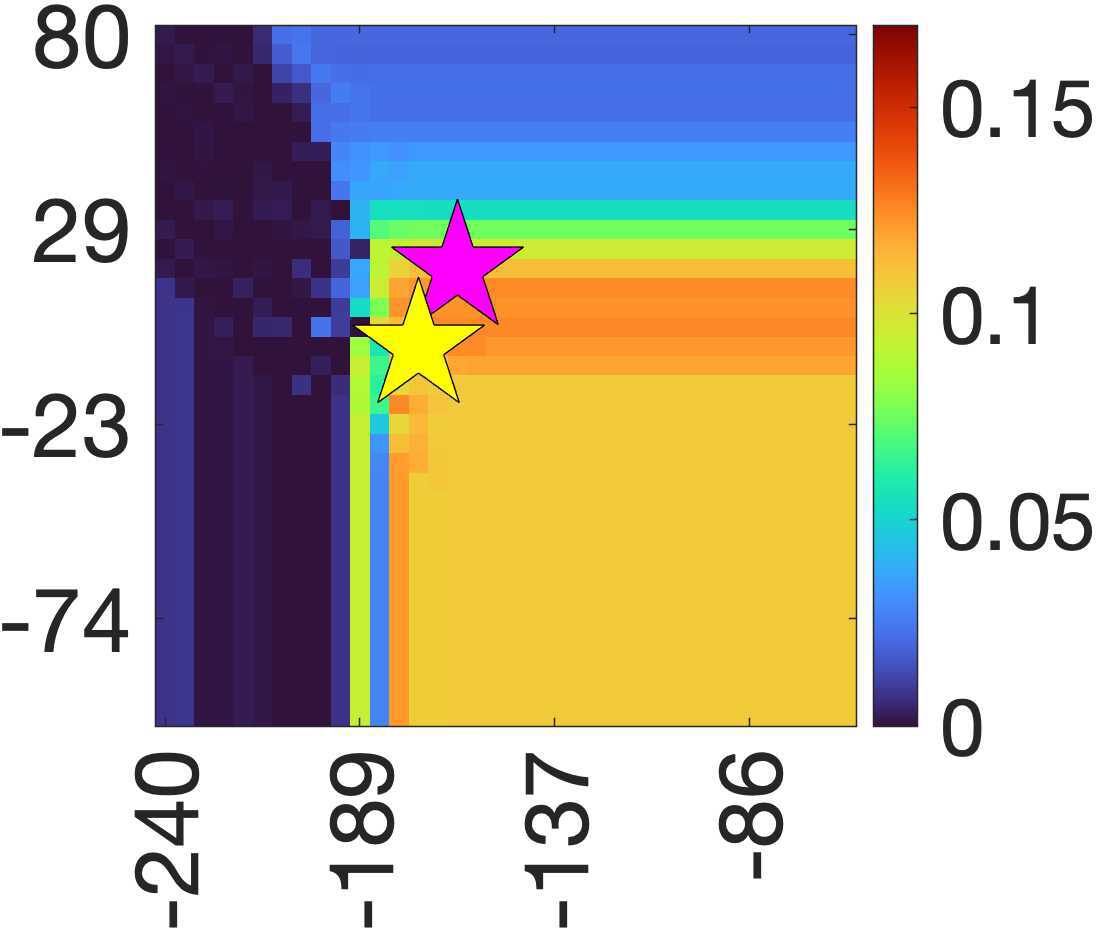}
            \end{minipage}
            \label{fig:imagesc_TLS_CD}
            \vskip0.1cm
            \begin{minipage}{0.2cm}
                \rotatebox{90}{$\Theta$}
            \end{minipage}
            \begin{minipage}{2.4cm}
                \centering 
                \includegraphics[height=2.0cm]{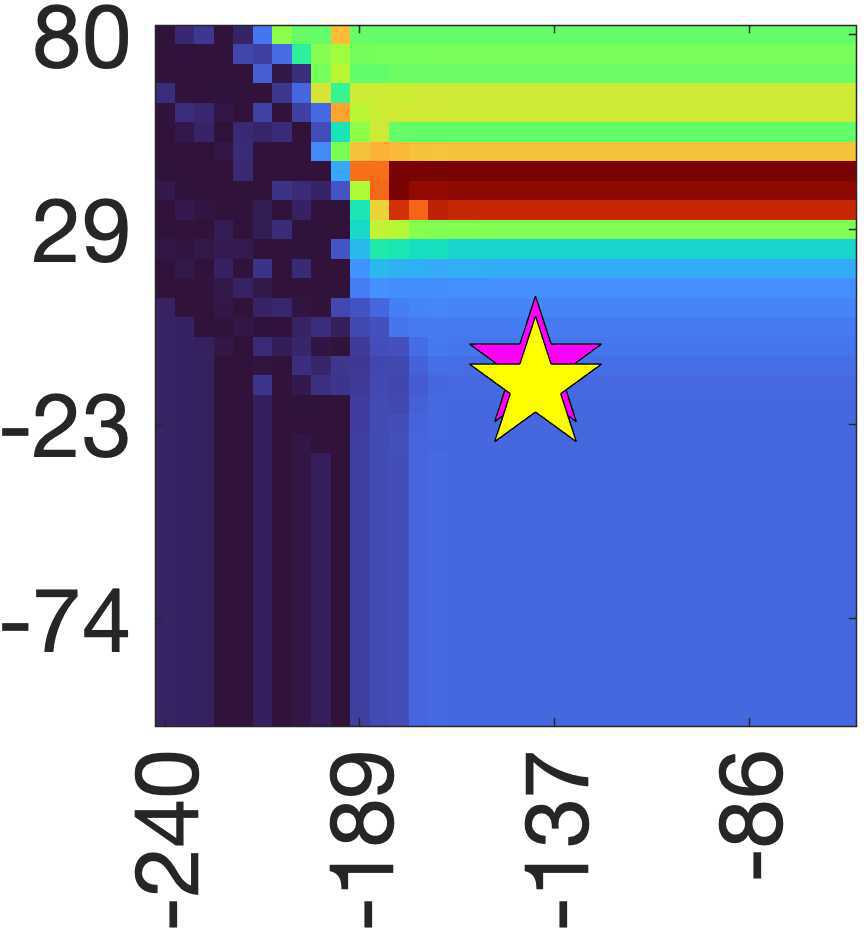}
            \end{minipage}
            \begin{minipage}{2.4cm}
                \centering
                 \includegraphics[height=2.0cm]{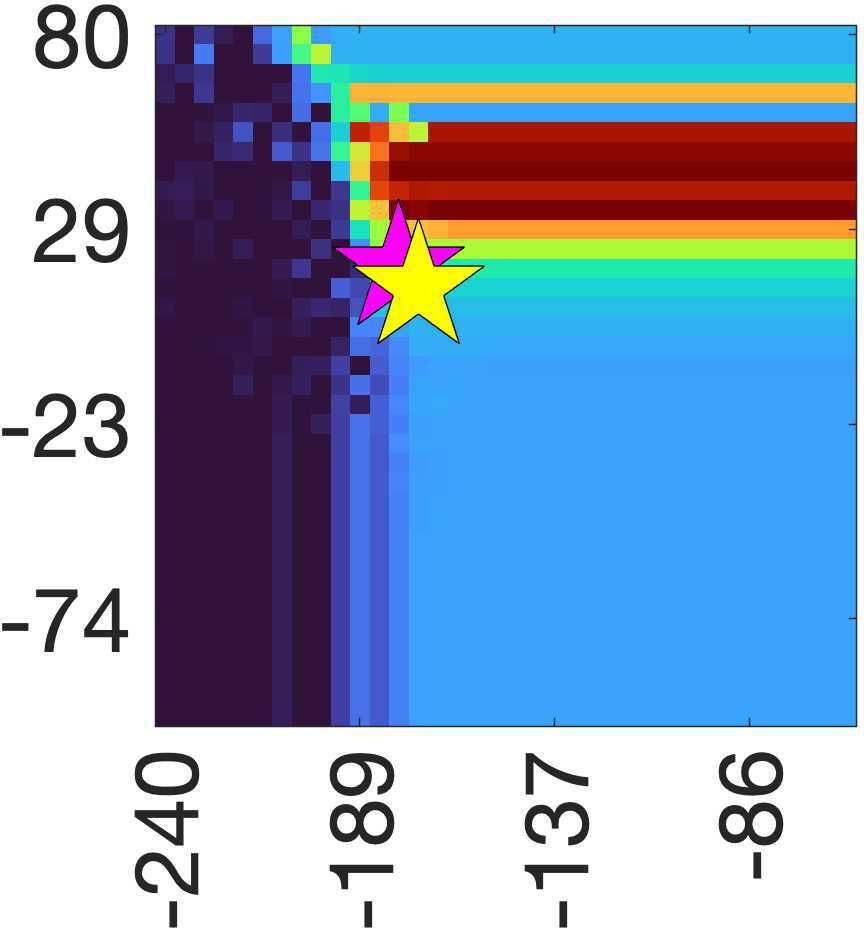}
            \end{minipage}
            \begin{minipage}{2.4cm}
                \centering
                \includegraphics[height=2.0cm]{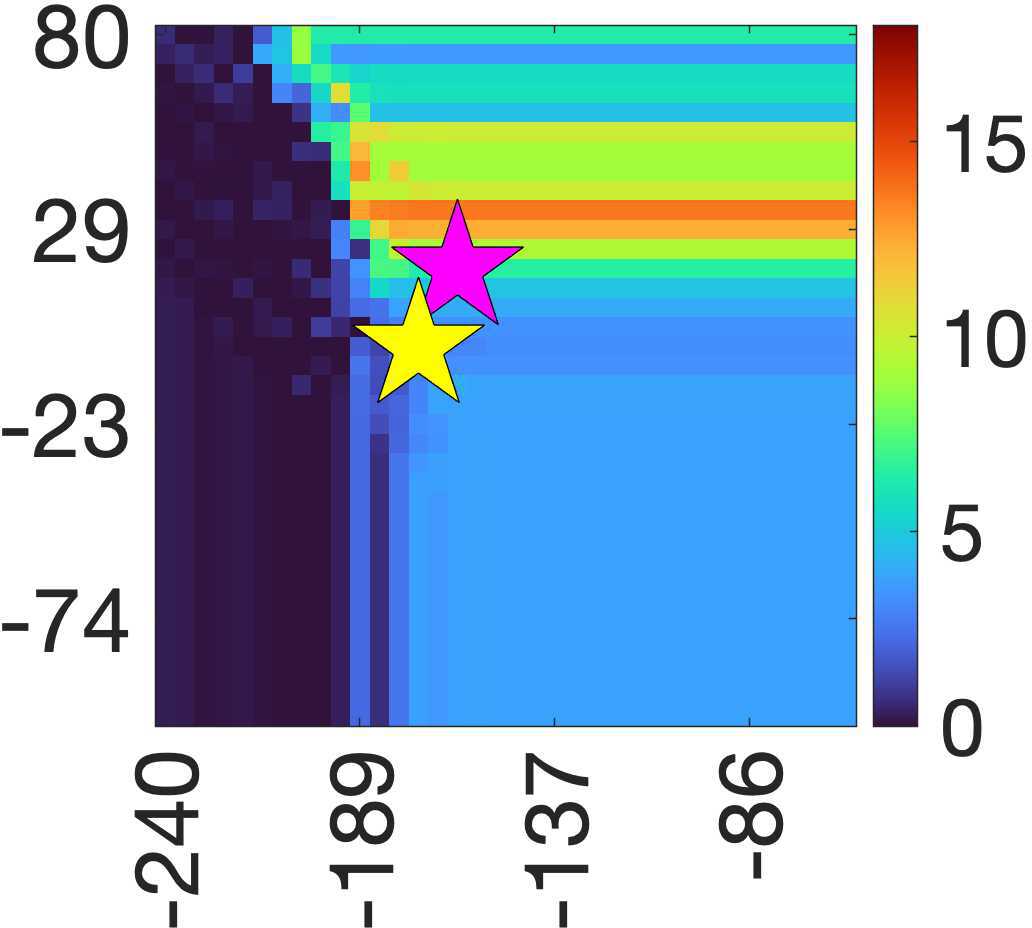}
            \end{minipage}
            \vskip0.1cm
            \begin{minipage}{0.2cm}
                \rotatebox{90}{AD (deg)}
            \end{minipage}
            \begin{minipage}{2.4cm}
                \centering 
                \includegraphics[height=2.0cm]{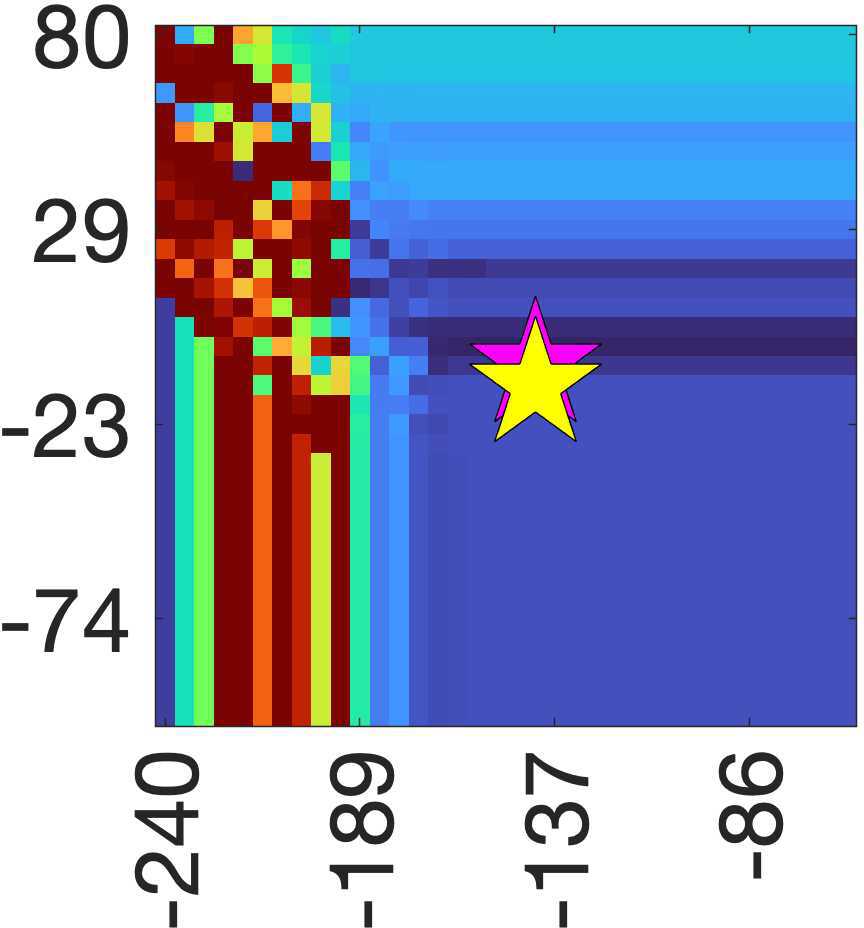}
            \end{minipage}
            \begin{minipage}{2.4cm}
                \centering
                \includegraphics[height=2.0cm]{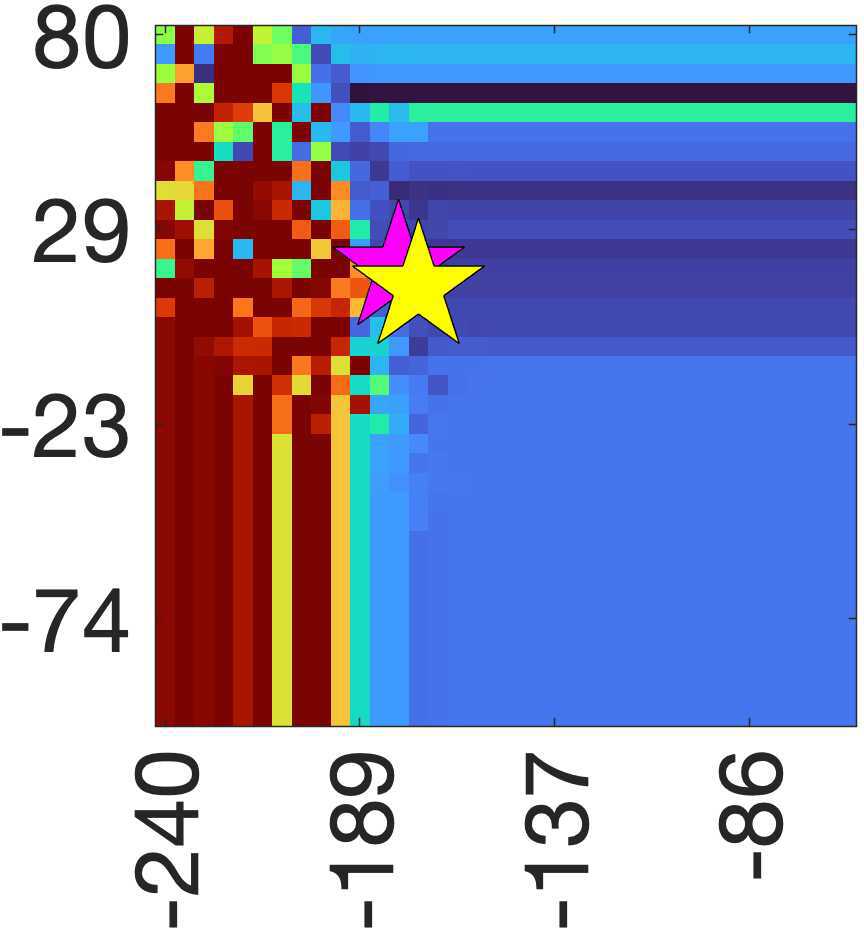}
            \end{minipage}
            \begin{minipage}{2.4cm}
                \centering
               \includegraphics[height=2.0cm]{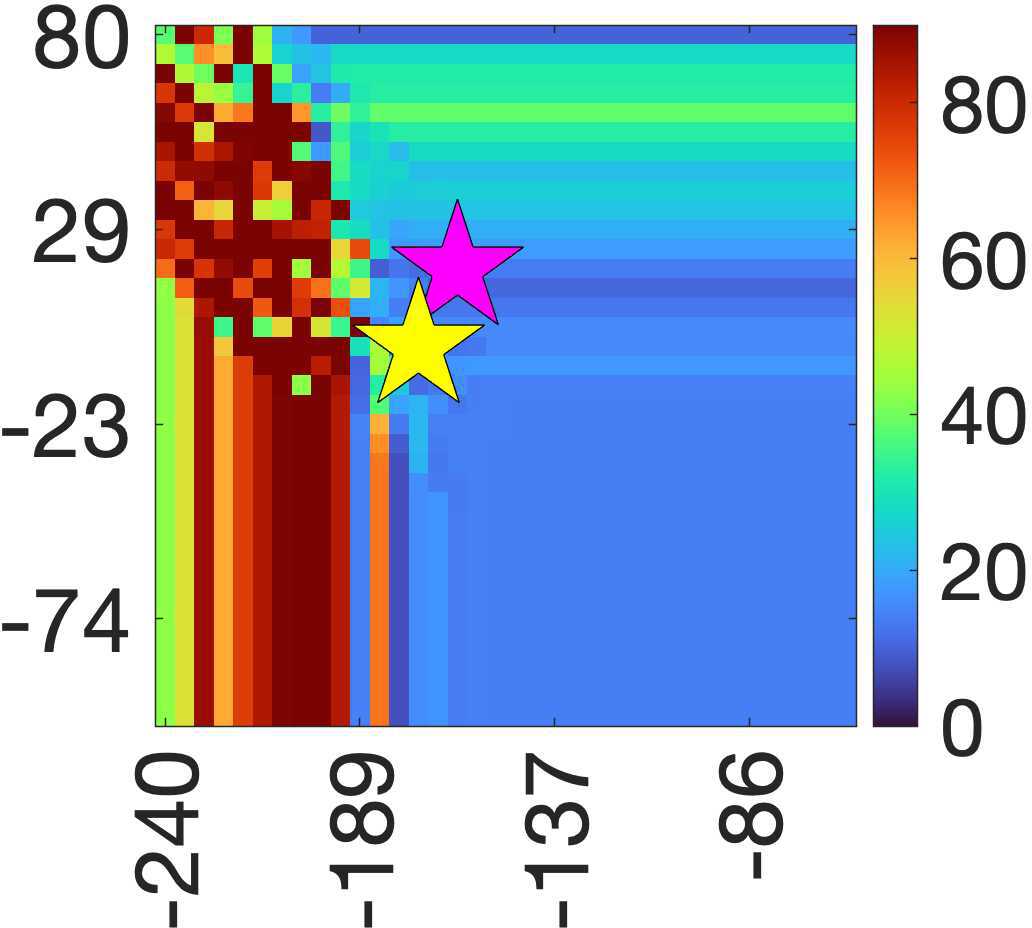}
            \end{minipage}
            \vskip0.1cm
            \begin{minipage}{0.2cm}
                \rotatebox{90}{$\| {\bf y}  \|_\infty$ (mA)}
            \end{minipage}
            \begin{minipage}{2.4cm}
                \centering 
                \includegraphics[height=2.0cm]{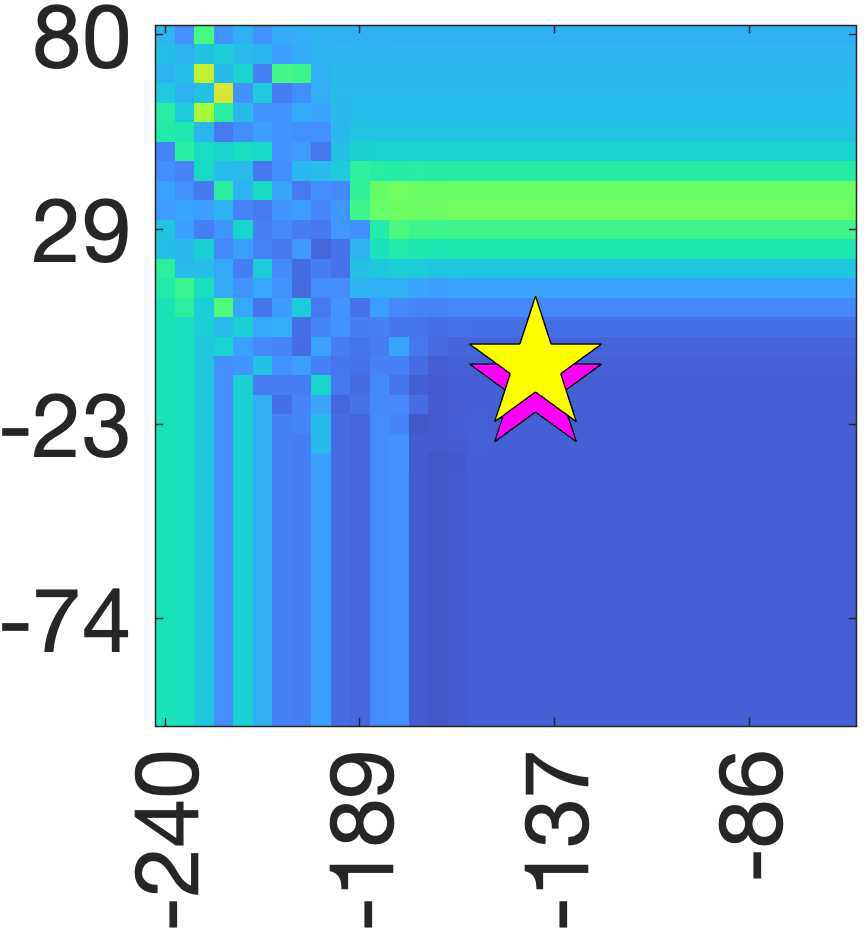}
            \end{minipage}
            \begin{minipage}{2.4cm}
                \centering
                 \includegraphics[height=2.0cm]{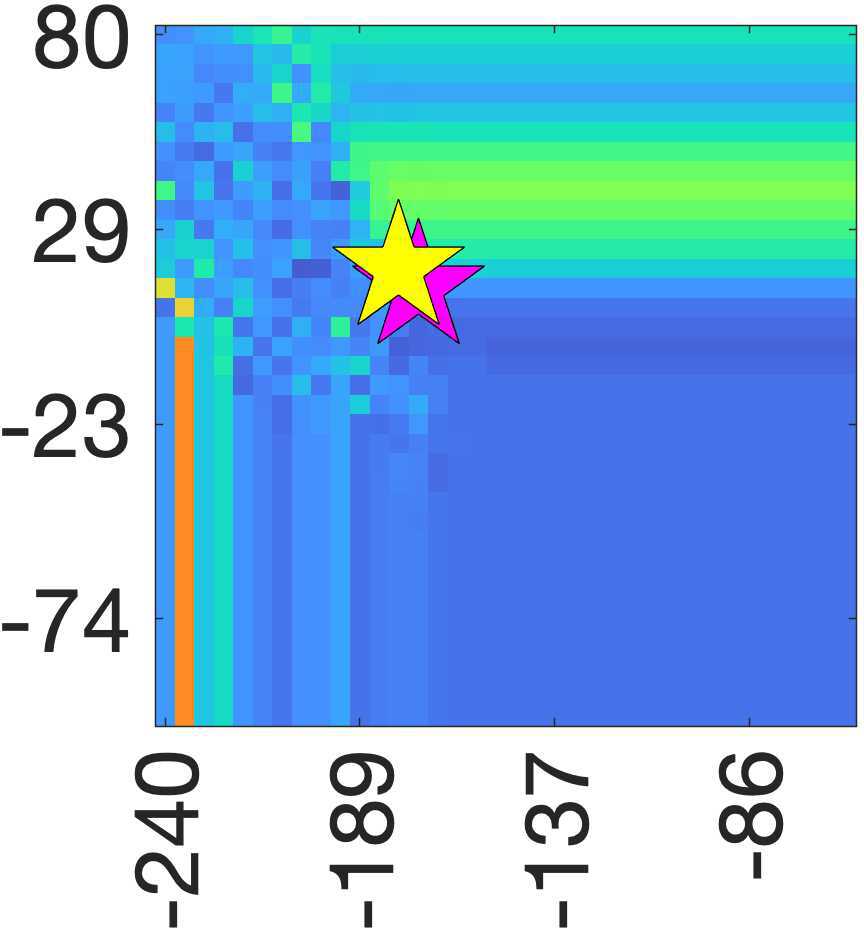}
            \end{minipage}
            \begin{minipage}{2.4cm}
                \centering
                \includegraphics[height=2.0cm]{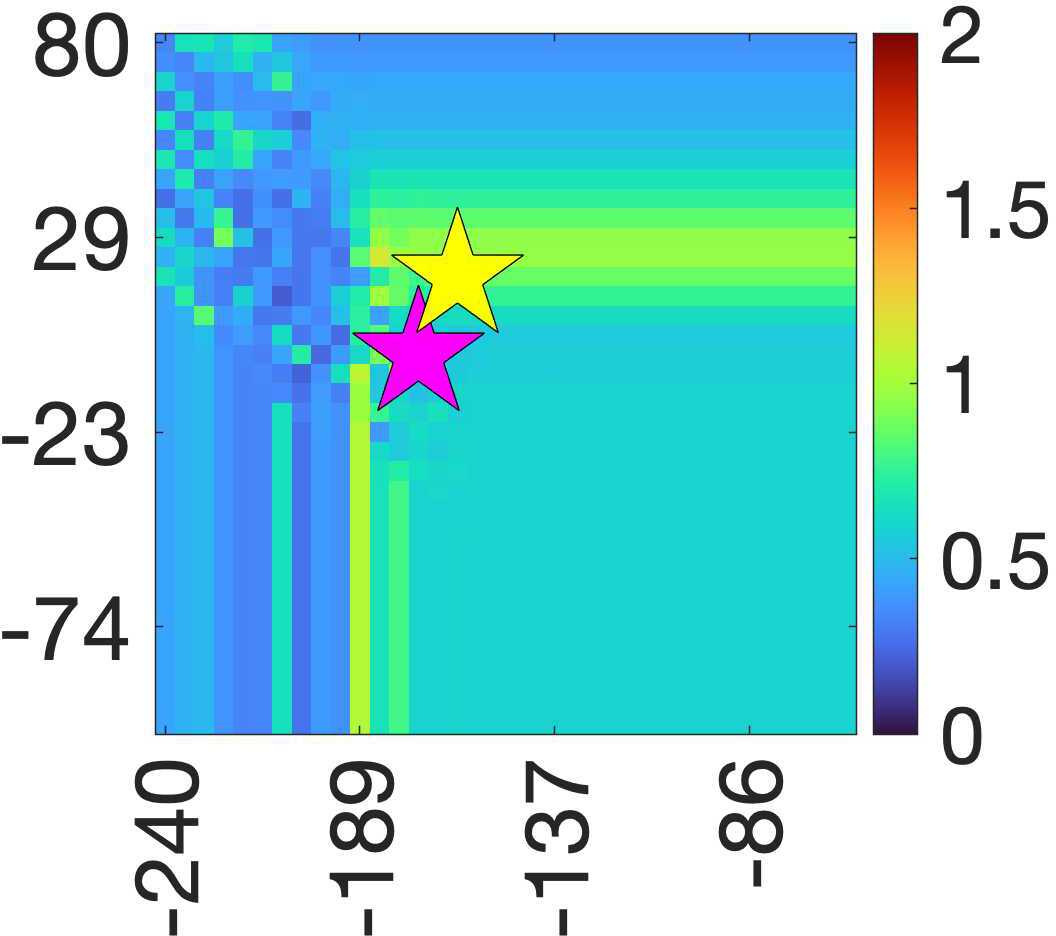}
            \end{minipage}
            \vskip0.1cm
            \hrule
        \end{minipage}
    \end{scriptsize}
\caption{Tikhonov regularized least-squares (TLS) charts showing $\Gamma$ (A/m\textsuperscript{2}), $\Theta$ (relative), AD (deg), and $\| {\bf y} \|_\infty$ (mA) for each point of the search lattice in the first run of the two-stage metaheuristic optimization process. Horizontal axis corresponds to the regularization parameter $\alpha$ and the vertical to the nuisance field weight $\delta$. The optimal solution for the case (A) is represented by a purple star, while for the case (B) as a yellow star. Axis are decibel (dB) scaled.
\label{fig:imagesc_TLS}}
\end{figure}

\begin{figure}[h!]
    \centering
    \begin{scriptsize} 
        \mbox{}
        \hskip-0.25cm
        \includegraphics[height=1.7cm]{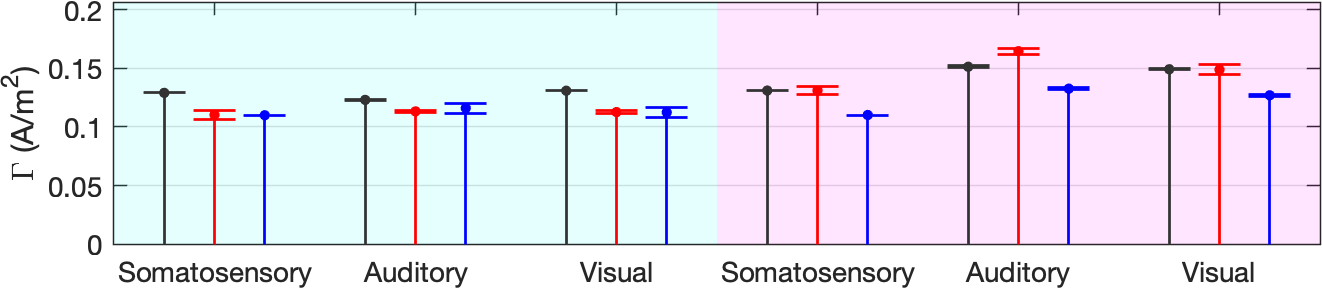} \\  \vskip0.1cm
        \includegraphics[height=1.7cm]{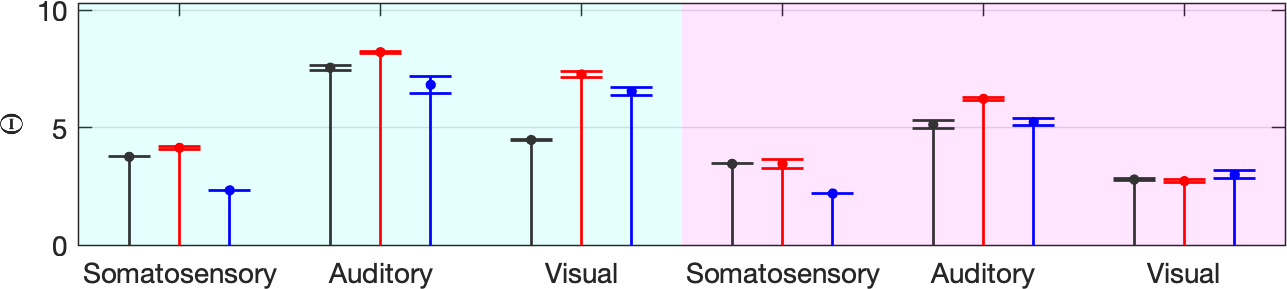} \\ \vskip0.1cm
        \includegraphics[height=1.7cm]{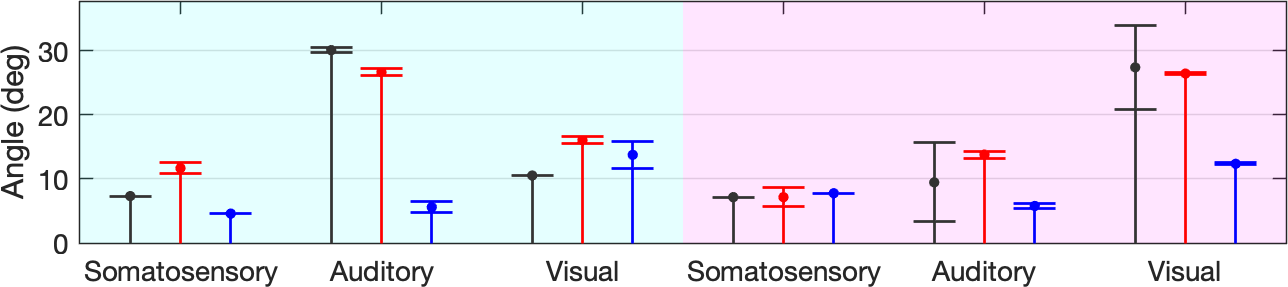} \\ \vskip0.1cm
        \includegraphics[height=1.7cm]{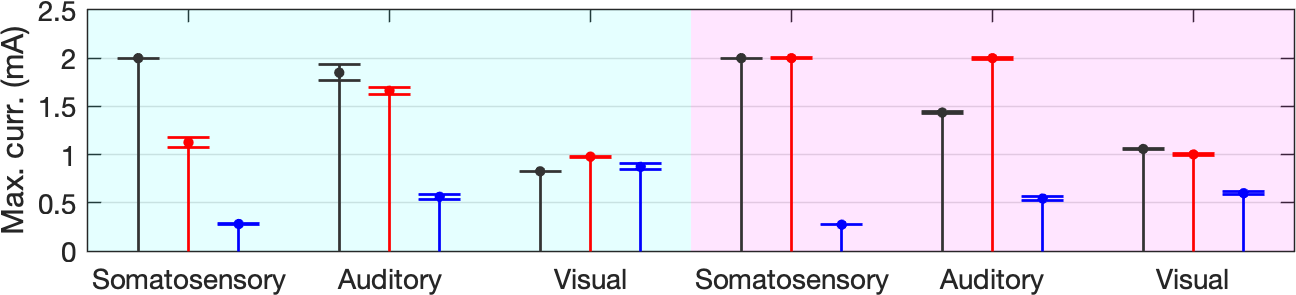} \\
        \includegraphics[width=4cm]{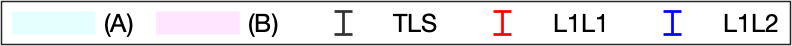} \\
        \textbf{20-electrode montage} \\ \vskip0.50cm  
        \mbox{} \hskip-0.25cm
        \includegraphics[height=1.7cm]{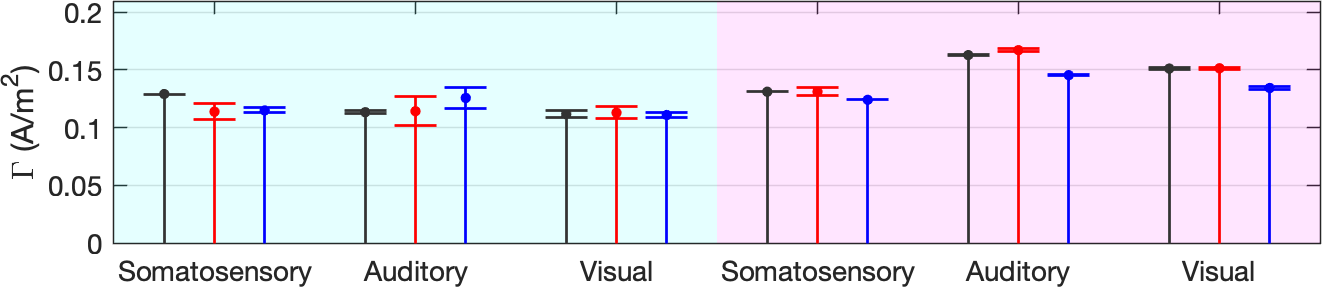} \\  \vskip0.1cm
        \includegraphics[height=1.7cm]{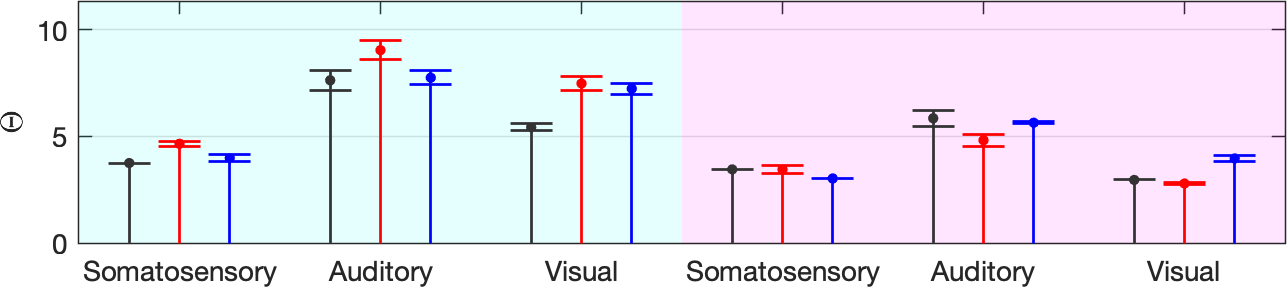} \\ \vskip0.1cm
        \includegraphics[height=1.7cm]{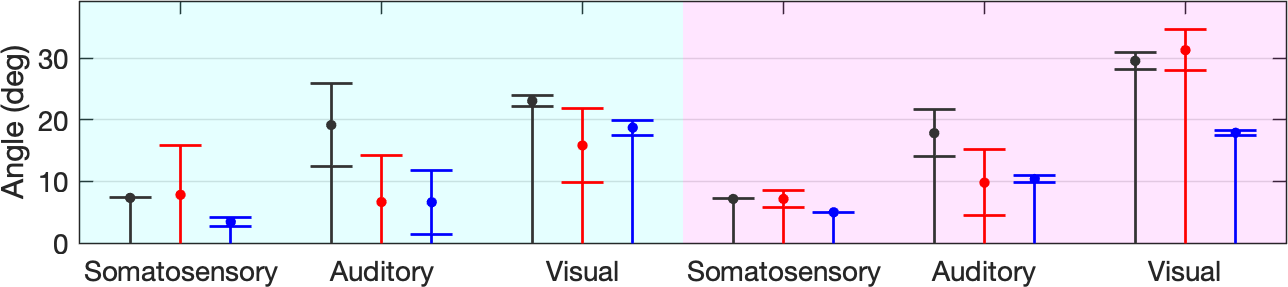} \\ \vskip0.1cm
        \includegraphics[height=1.7cm]{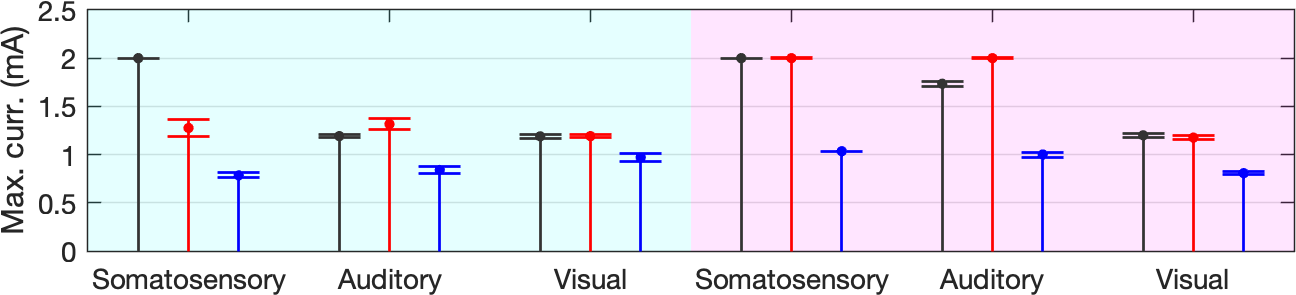} \\ \vskip0.1cm
        \includegraphics[width=4cm]{images_jpg/revision/result_legend.png} \\
        \textbf{8-electrode montage}
    \end{scriptsize}
\caption{A graphical illustration optimized values of $\Gamma$ (A/m\textsuperscript{2}), $\Theta$ (relative), AD (deg), and $\| {\bf y} \|_\infty$ (mA) and their maximal estimated lattice-induced deviations. On each case the dot on the stem shows the optimizer while the whiskers illustrate the estimated limits for the deviation. The TLS, L1L1 and L1L2 results correspond to dark grey, red and blue stems, and the optimization cases (A) and (B) to the cyan and magenta background color, respectively.}
\label{fig:stem_comparison}
\end{figure}

For each optimization method, the range of the metaheuristic lattice search was effective enough to cover those regions wherein the magnitude of the focused current density $\Gamma$ was close to its maximum. These regions were surrounded by comparably smooth transition zones from high to low values, showing the regularity of the optimization process with respect to parameter variation. The solutions obtained through the metaheuristic search case (A), i.e., the maximum of $\Theta$ for $\Gamma \geq 0.11$ A/m\textsuperscript{2}, were generally found from these transition zones. In the case (B), the maximum of $\Gamma$ was found from a regular region, where the variation of $\Gamma$ stayed on a comparably low level. The charts showing $\Gamma$, $\Theta$, AD, and $\| {\bf y} \|_\infty$ for the lattices of the first metaheuristic search run (with 128 electrode positions) together with the corresponding estimates found in the cases (A) and (B) have been included in Figures \ref{fig:imagesc_L1L1}, \ref{fig:imagesc_L1L2} and \ref{fig:imagesc_TLS}. The charts illustrated are shown with respect to the 20 most intense electric potential channels in the current pattern. The optimizer found in the case (A), where $\Theta$ is maximized in the second optimization stage for those lattice points which satisfy $\Gamma \geq 0.11$ A/m\textsuperscript{2} in the first stage, represented by a purple star, while in the case (B), in which the global maximizer of $\Gamma$ is found, as yellow star. Differences between the first and second runs, as well as electrode montage setups were minor, other charts are not shown.

The quantities $\Gamma$, $\Theta$, AD, and $\| {\bf y} \|_\infty$ and their estimated maximal deviations corresponding to the optimizers found are shown in Table \ref{table:method_comparison} and Figure \ref{fig:stem_comparison}. Compared to L1L2 and TLS, L1L1 yielded a greater or equal value of $\Theta$ and $\Gamma$ in cases (A) and (B), respectively. Agreeing with our initial hypothesis, the extra gain provided by L1L1 as compared to L1L2 was observed to be systematic and the most pronounced in case (A), where the optimized $\Theta$-value of L1L1 was 1.6 and 1.4 times that of L1L2. Considering TLS, L1L1 yielded a greater maximum current $\| {\bf y} \|_\infty$ in all the examined cases, while the optimizers tend to deviate overall somewhat more than with TLS. The greatest difference between the L1L1- and TLS-optimized focal current density $\Gamma$ was obtained for the auditory target current with both 20- and 8-electrode montage for which this difference was 0.032 and 0.022 A/m\textsuperscript{2}, respectively.

In most cases, the estimated maximum level of deviation does not exceed the observed differences between the optimized values, confirming that the mutual performance differences between the optimization methods. The results concerning AD are somewhat more obscure than for the other examined quantities due to comparably larger estimates for the maximal deviation. The tendency of TLS  to result in a  smaller AD compared to L1L1 and L1L2 seems obvious. Notably, in a mutual comparison between the results obtained with the 20- and 8-electrode montage, the latter was observed to result in overall smaller mutual differences between the methods.

The regularity of the candidate solutions was found to decrease towards the boundaries of the search lattice, which was reflected as a somewhat elevated deviation estimate, when an optimizer was found close to a boundary. In particular, the L1L1 method did not found an optimizer, when the nuisance field weight was greater than one, which can be observed from the charts in Figure \ref{fig:stem_comparison}.

\subsection{Current pattern and volume density}
\label{sec:Results_ROI}

\begin{figure*}[h!]
\centering
    \begin{scriptsize}
    \centering
    \begin{minipage}{18.0cm}
        \centering
        \begin{minipage}{0.2cm}
        \centering
           \mbox{}
        \end{minipage}
        \begin{minipage}{16.50cm}
        \centering
            \begin{minipage}{0.25cm}
                \mbox{}
            \end{minipage}
            \begin{minipage}{5.30cm}
            \centering
                \textbf{L1-norm regularized \\ L1-norm fitting (L1L1)}
            \end{minipage}
            \begin{minipage}{5.30cm}
                \centering
                \textbf{L1-norm regularized \\ L2-norm fitting (L1L2)}
            \end{minipage}
            \begin{minipage}{5.30cm}
                \centering
                \textbf{Tikhonov regularized \\ least-squares (TLS)}
            \end{minipage}
        \end{minipage}
    \end{minipage}
    \vskip0.1cm \hrule \vskip0.1cm
    \begin{minipage}{18.0cm}
    \centering
        \begin{minipage}{0.2cm}
        \centering
            \rotatebox{90}{\textbf{Somatosensory}}
        \end{minipage}
        \begin{minipage}{16.50cm}
        \centering
            \begin{minipage}{0.25cm}
            \centering
                \rotatebox{90}{\textbf{Case A}}
            \end{minipage}
            \begin{subfigure}{5.3cm}
            \centering
                \begin{minipage}{2.6cm}
                    \centering
                    \includegraphics[height = 1.75cm]{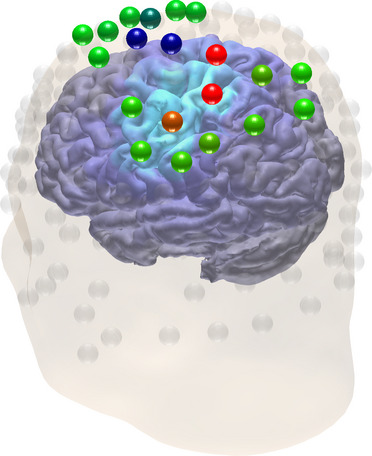} \\
                    \includegraphics[width = 2.30cm,height=1.2cm]{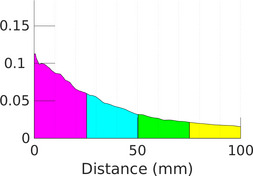} 
                \end{minipage} 
                \begin{minipage}{2.6cm}
                \centering 
                    \includegraphics[height = 2.55cm]{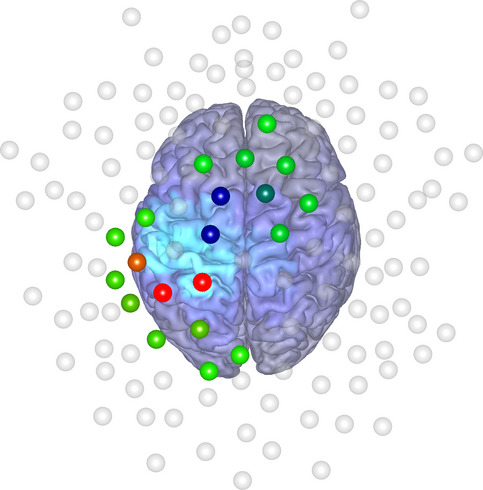}
                \end{minipage} \\
                \label{fig:postcentral_LP_1_20_results}
            \end{subfigure}
            \begin{subfigure}{5.3cm}
            \centering
                \begin{minipage}{2.60cm}
                    \centering
                    \includegraphics[height = 1.75cm]{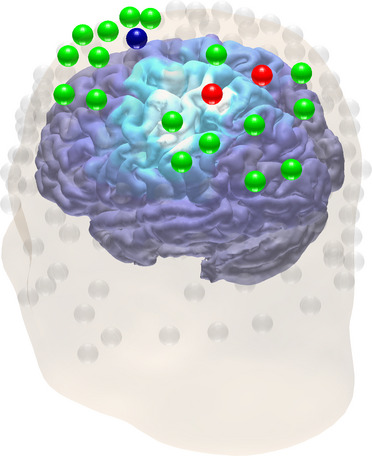} \\
                    \includegraphics[width = 2.30cm,height=1.2cm]{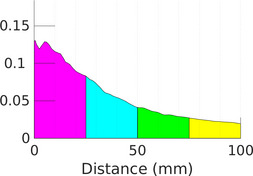} 
                \end{minipage} 
                \begin{minipage}{2.60cm}
                \centering 
                    \includegraphics[height = 2.55cm]{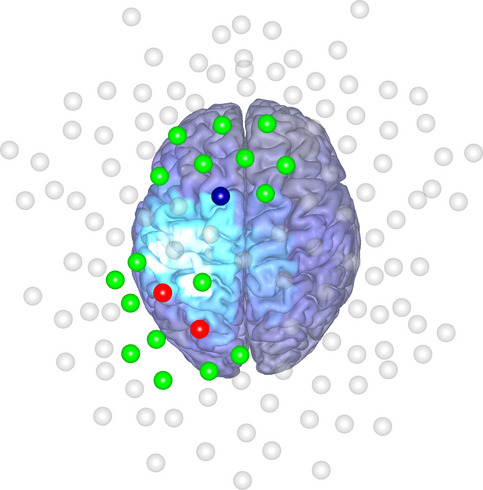}
                \end{minipage}
                \label{fig:postcentral_SDP_1_20_results}
            \end{subfigure}
            \begin{subfigure}{5.3cm}
            \centering
                \begin{minipage}{2.60cm}
                    \centering
                    \includegraphics[height = 1.75cm]{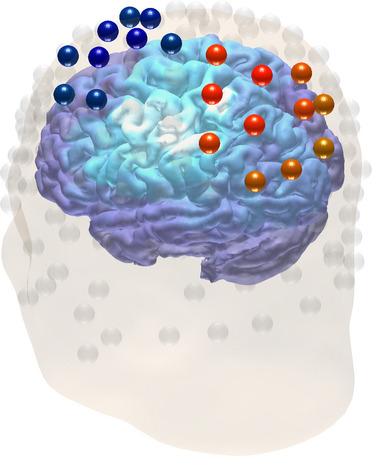} \\
                    \includegraphics[width = 2.30cm,height=1.2cm]{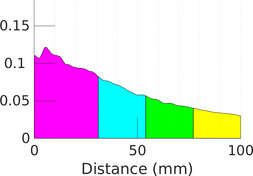}
                \end{minipage} 
                \begin{minipage}{2.60cm}
                \centering 
                    \includegraphics[height = 2.55cm]{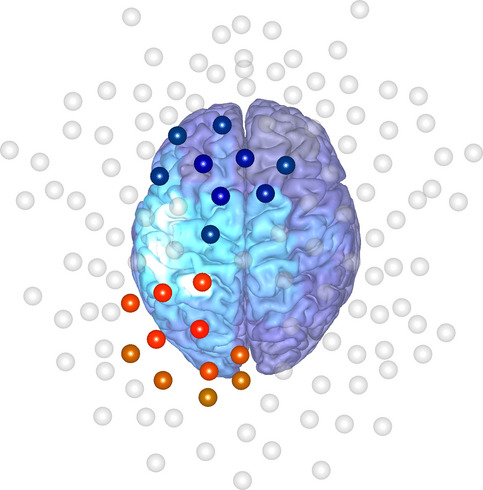}
                \end{minipage}
                \label{fig:postcentral_TLS_1_20_results}
            \end{subfigure}
            \vskip0.15cm
            \begin{minipage}{0.25cm}
            \centering
                \rotatebox{90}{\textbf{Case B}}
            \end{minipage}
            \begin{subfigure}{5.3cm}
            \centering
                \begin{minipage}{2.60cm}
                \centering
                    \includegraphics[height = 1.75cm]{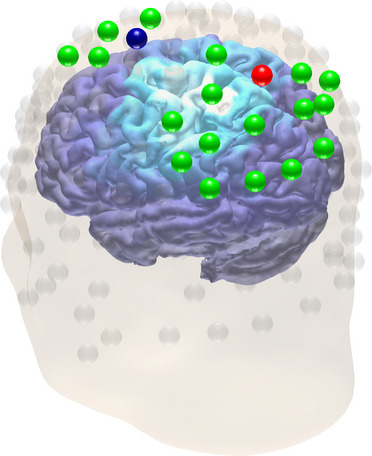} \\
                    \includegraphics[width = 2.30cm,height=1.2cm]{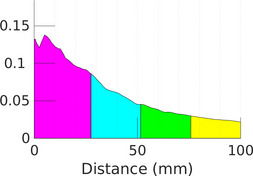}
                \end{minipage}
                \begin{minipage}{2.6cm}
                \centering 
                    \includegraphics[height = 2.55cm]{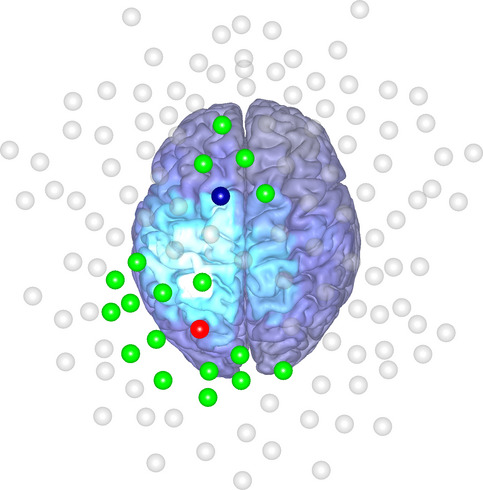}
                \end{minipage}
                \label{fig:postcentral_LP_2_20_results}
            \end{subfigure}
            \begin{subfigure}{5.3cm}
            \centering
                \begin{minipage}{2.60cm}
                    \centering
                    \includegraphics[height = 1.75cm]{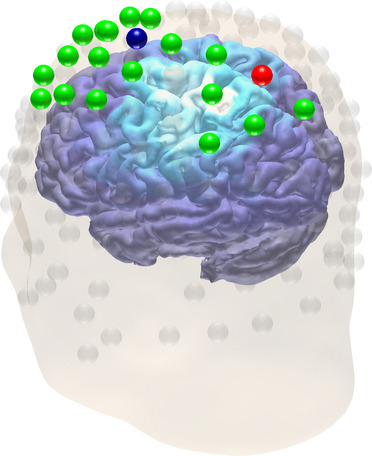} \\
                    \includegraphics[width = 2.30cm,height=1.2cm]{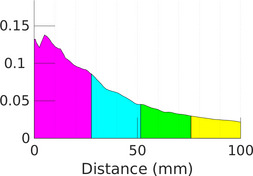}
                \end{minipage} 
                \begin{minipage}{2.60cm}
                \centering 
                    \includegraphics[height = 2.55cm]{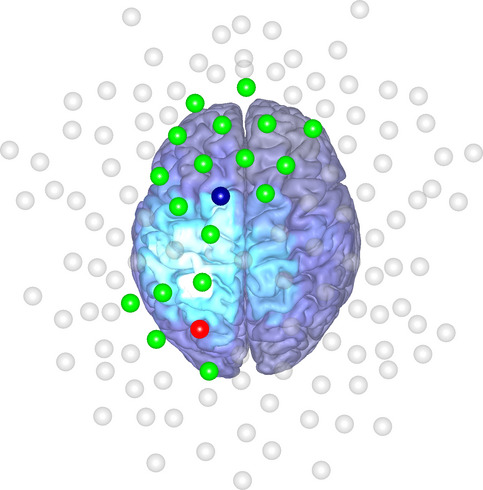}
                \end{minipage}
                \label{fig:postcentral_SDP_2_20_results}
            \end{subfigure}
            \begin{subfigure}{5.3cm}
            \centering
                \begin{minipage}{2.60cm}
                    \centering
                    \includegraphics[height = 1.75cm]{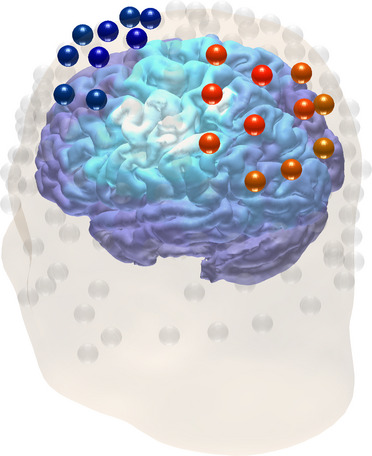} \\
                     \includegraphics[width = 2.30cm,height=1.2cm]{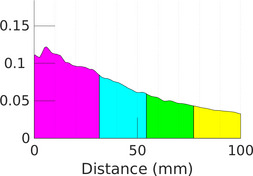} 
                \end{minipage} 
                \begin{minipage}{2.60cm}
                \centering 
                    \includegraphics[height = 2.55cm]{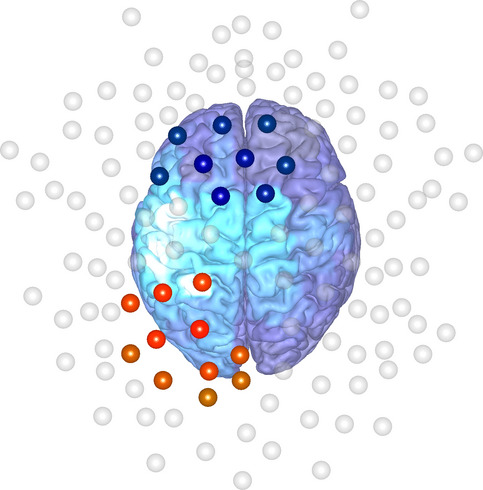}
                \end{minipage}
                \label{fig:postcentral_TLS_2_20_results}
            \end{subfigure}
        \end{minipage}
    \end{minipage}
    \vskip0.1cm \hrule \vskip0.1cm
    \begin{minipage}{18.0cm}
    \centering
        \begin{minipage}{0.2cm}
        \centering
            \rotatebox{90}{\textbf{Auditory}}
        \end{minipage}
        \begin{minipage}{16.50cm}
        \centering
            \begin{minipage}{0.25cm}
            \centering
                \rotatebox{90}{\textbf{Case A}}
            \end{minipage}
            \begin{subfigure}{5.3cm}
            \centering
                \begin{minipage}{2.6cm}
                    \centering
                    \includegraphics[height = 1.75cm]{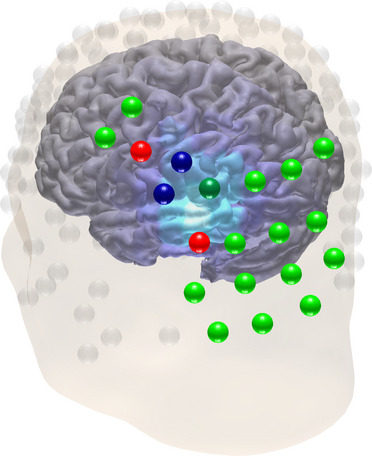} \\
                    \includegraphics[width = 2.30cm,height=1.2cm]{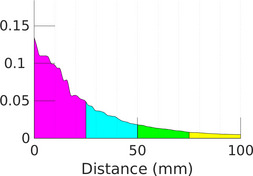} 
                \end{minipage} 
                \begin{minipage}{2.6cm}
                \centering 
                    \includegraphics[height = 2.55cm]{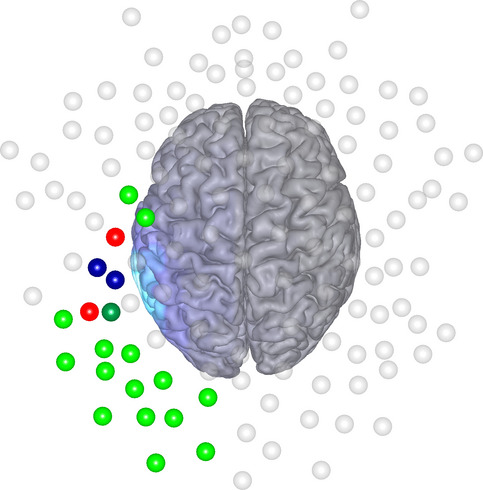}
                \end{minipage}
                \label{fig:auditory_LP_1_20_results}
            \end{subfigure}
            \begin{subfigure}{5.3cm}
            \centering
                \begin{minipage}{2.60cm}
                    \centering
                    \includegraphics[height = 1.75cm]{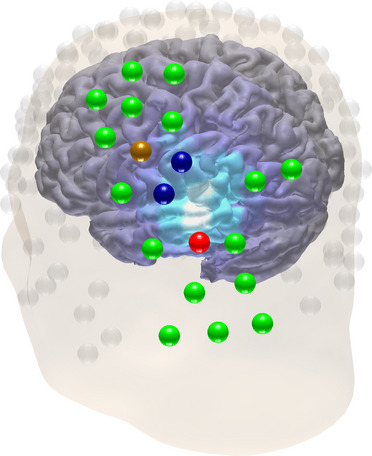} \\
                    \includegraphics[width = 2.30cm,height=1.2cm]{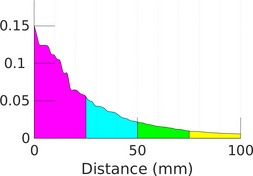} 
                \end{minipage}
                \begin{minipage}{2.60cm}
                \centering 
                    \includegraphics[height = 2.55cm]{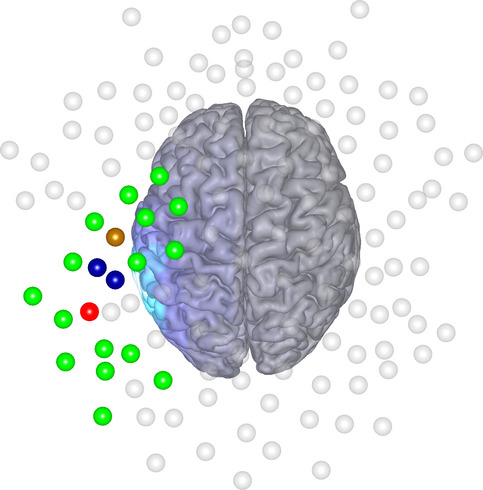}
                \end{minipage}
                \label{fig:auditory_SDP_1_20_results}
            \end{subfigure}
            \begin{subfigure}{5.3cm}
            \centering
                \begin{minipage}{2.60cm}
                    \centering
                    \includegraphics[height = 1.5cm]{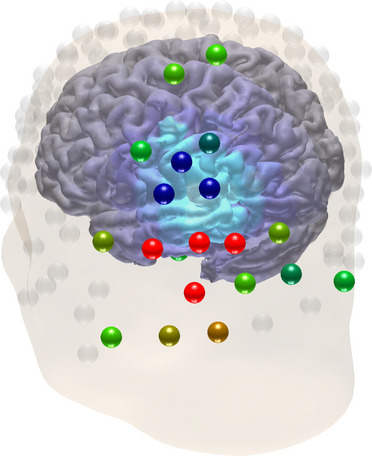} \\
                    \includegraphics[width = 2.6cm,height=1.2cm]{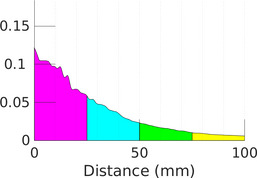}
                \end{minipage} 
                \begin{minipage}{2.60cm}
                \centering 
                    \includegraphics[height = 2.55cm]{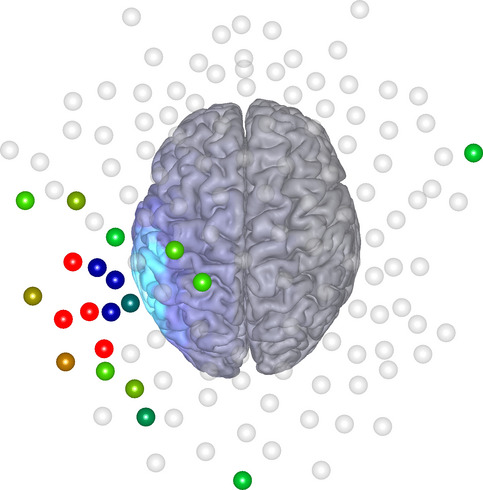}
                \end{minipage}
                \label{fig:auditory_TLS_1_20_results}
            \end{subfigure}
            \vskip0.15cm
            \begin{minipage}{0.25cm}
            \centering
                \rotatebox{90}{\textbf{Case B}}
            \end{minipage}
            \begin{subfigure}{5.3cm}
            \centering
                \begin{minipage}{2.60cm}
                    \centering
                    \includegraphics[height = 1.75cm]{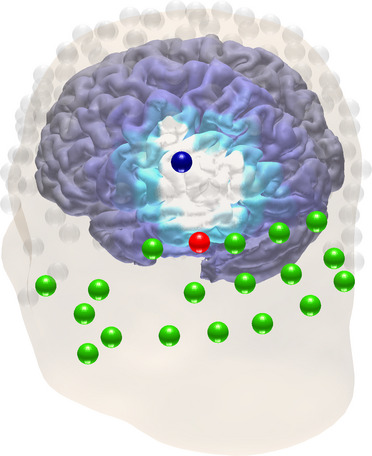} \\
                    \includegraphics[width = 2.30cm,height=1.2cm]{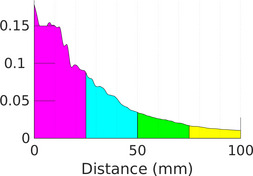}
                \end{minipage}
                \begin{minipage}{2.60cm}
                \centering 
                    \includegraphics[height = 2.55cm]{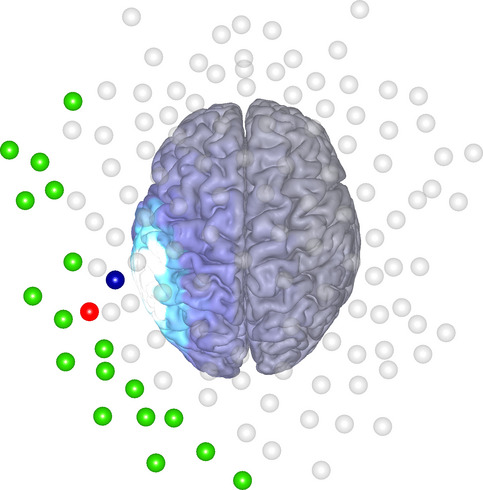}
                \end{minipage}
                \label{fig:auditory_LP_2_20_results}
            \end{subfigure}
            \begin{subfigure}{5.3cm}
            \centering
                \begin{minipage}{2.60cm}
                    \centering
                    \includegraphics[height = 1.75cm]{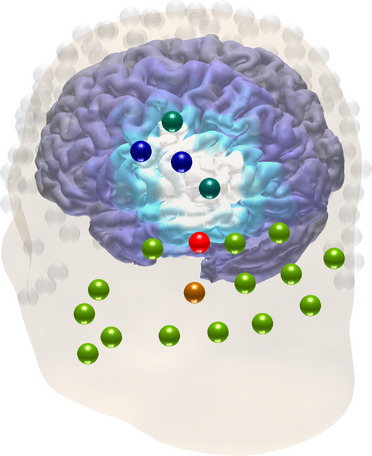} \\
                    \includegraphics[width = 2.30cm,height=1.2cm]{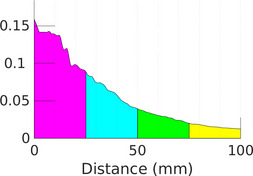}
                \end{minipage} 
                \begin{minipage}{2.60cm}
                \centering 
                    \includegraphics[height = 2.55cm]{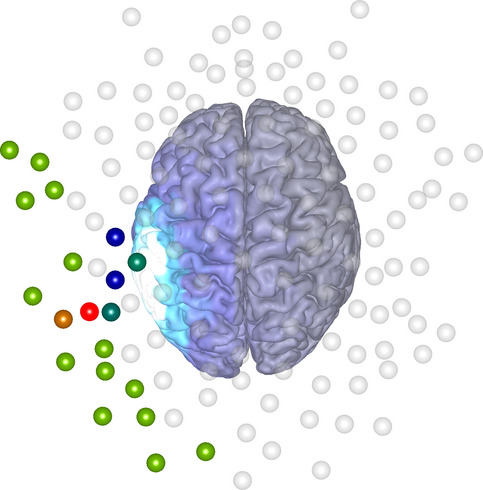}
                \end{minipage}
                \label{fig:auditory_SDP_2_20_results}
            \end{subfigure}
            \begin{subfigure}{5.3cm}
            \centering
                \begin{minipage}{2.60cm}
                    \centering
                    \includegraphics[height = 1.75cm]{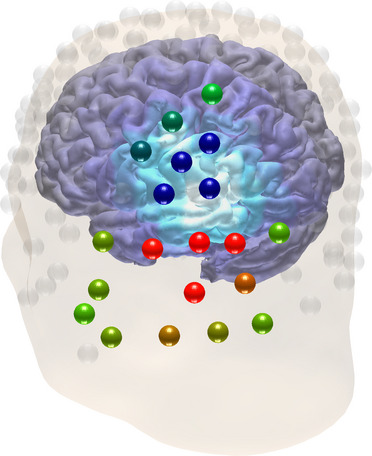} \\
                    \includegraphics[width = 2.30cm,height=1.2cm]{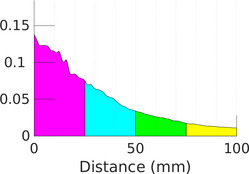} 
                \end{minipage} 
                \begin{minipage}{2.60cm}
                \centering 
                    \includegraphics[height = 2.55cm]{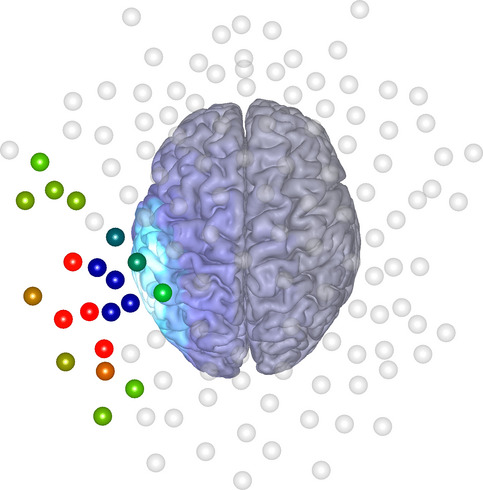}
                \end{minipage}
                \label{fig:auditory_TLS_2_20_results}
            \end{subfigure}
        \end{minipage}
    \end{minipage}
    \vskip0.1cm \hrule \vskip0.1cm
    \begin{minipage}{18.0cm}
    \centering
        \begin{minipage}{0.2cm}
        \centering
            \rotatebox{90}{\textbf{Visual}}
        \end{minipage}
        \begin{minipage}{16.50cm}
        \centering
            \begin{minipage}{0.25cm}
            \centering
                \rotatebox{90}{\textbf{Case A}}
            \end{minipage}
            \begin{subfigure}{5.3cm}
            \centering
                \begin{minipage}{2.60cm}
                \centering
                    \includegraphics[height = 1.75cm]{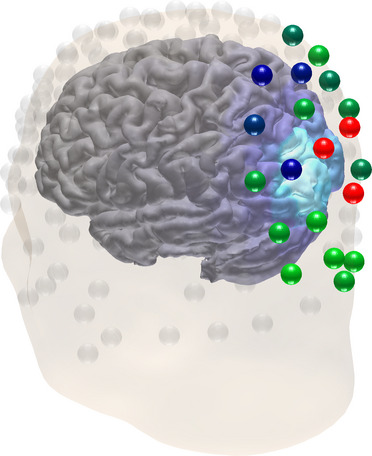} \\
                    \includegraphics[width = 2.30cm,height=1.2cm]{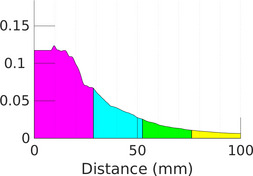} 
                \end{minipage} 
                \begin{minipage}{2.60cm}
                \centering 
                    \includegraphics[height = 2.55cm]{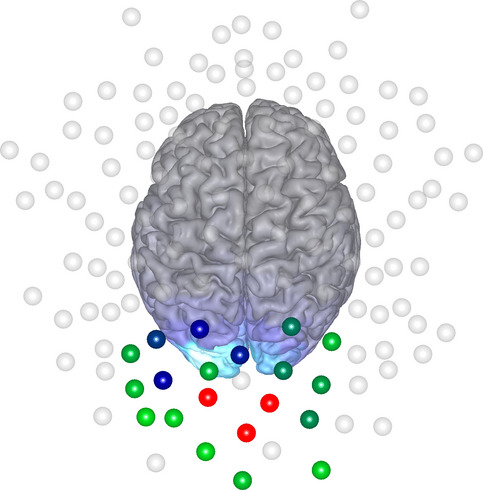}
                \end{minipage}
                \label{fig:occipital_LP_1_20_results}
            \end{subfigure}
            \begin{subfigure}{5.3cm}
            \centering
                \begin{minipage}{2.60cm}
                    \centering
                    \includegraphics[height = 1.75cm]{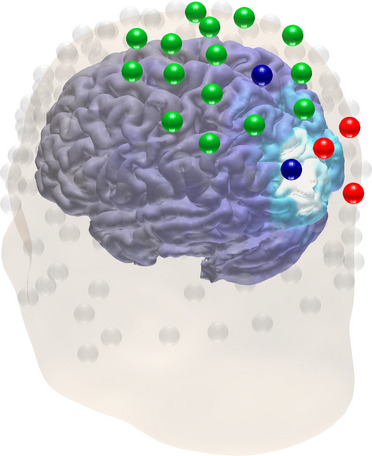} \\
                    \includegraphics[width = 2.30cm,height=1.2cm]{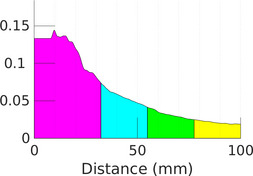} 
                \end{minipage} 
                \begin{minipage}{2.60cm}
                \centering 
                    \includegraphics[height = 2.55cm]{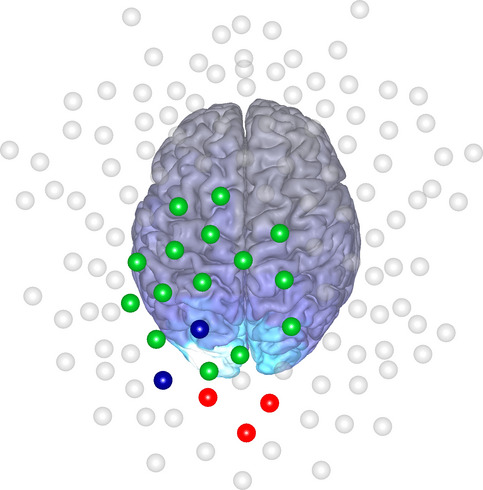}
                \end{minipage}
                \label{fig:occipital_SDP_1_20_results}
            \end{subfigure}
            \begin{subfigure}{5.3cm}
            \centering
                \begin{minipage}{2.60cm}
                    \centering
                    \includegraphics[height = 1.75cm]{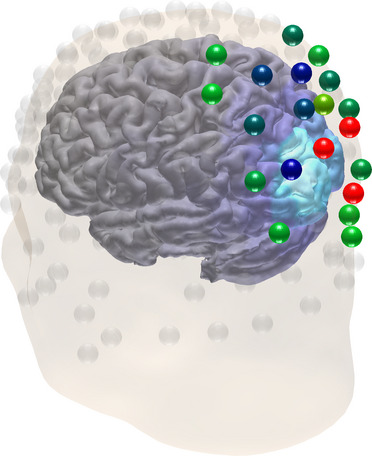} \\
                    \includegraphics[width = 2.30cm,height=1.2cm]{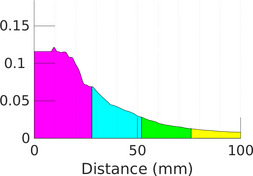}
                \end{minipage} 
                \begin{minipage}{2.60cm}
                \centering 
                    \includegraphics[height = 2.55cm]{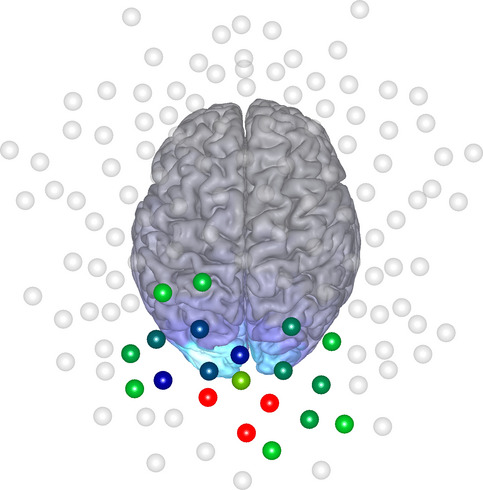}
                \end{minipage}
                \label{fig:occipital_TLS_1_20_results}
            \end{subfigure}
            \vskip0.15cm
            \begin{minipage}{0.25cm}
            \centering
                \rotatebox{90}{\textbf{Case B}}
            \end{minipage}
            \begin{subfigure}{5.3cm}
            \centering
                \begin{minipage}{2.60cm}
                    \centering
                    \includegraphics[height = 1.75cm]{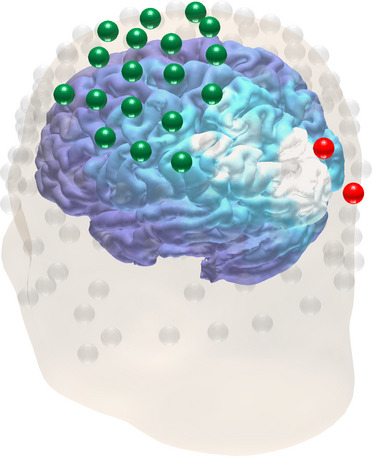} \\
                    \includegraphics[width = 2.30cm,height=1.2cm]{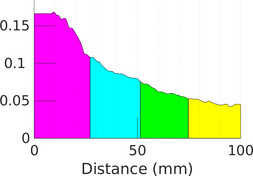}
                \end{minipage}
                \begin{minipage}{2.60cm}
                \centering 
                    \includegraphics[height = 2.55cm]{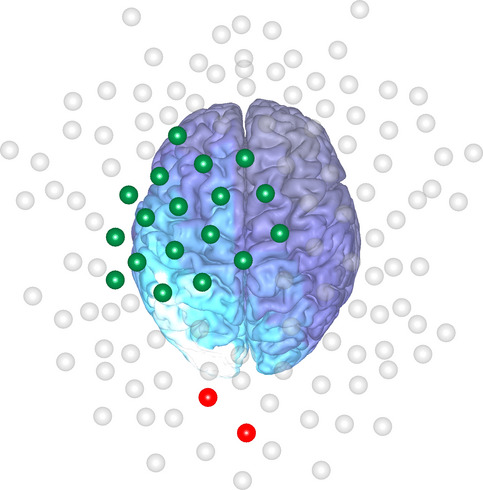}
                \end{minipage}
                \label{fig:occipital_LP_2_20_results}
            \end{subfigure}
            \begin{subfigure}{5.3cm}
            \centering
                \begin{minipage}{2.60cm}
                    \centering
                    \includegraphics[height = 1.75cm]{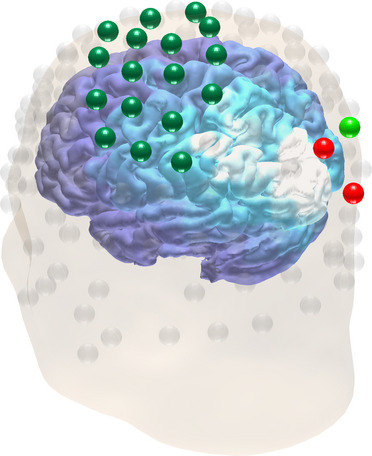} \\
                    \includegraphics[width = 2.30cm,height=1.2cm]{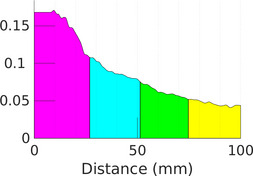}
                \end{minipage} 
                \begin{minipage}{2.60cm}
                \centering 
                    \includegraphics[height = 2.55cm]{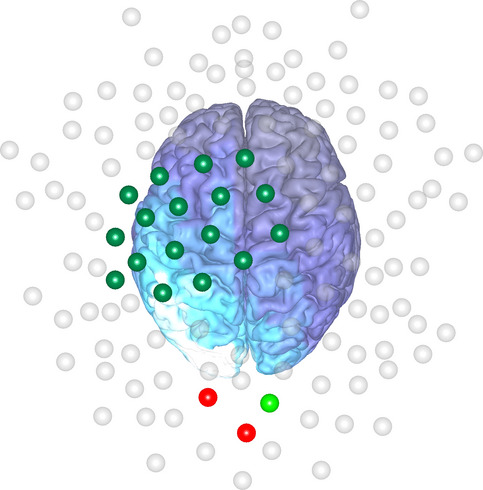}
                \end{minipage}
                \label{fig:occipital_SDP_2_20_results}
            \end{subfigure}
            \begin{subfigure}{5.3cm}
            \centering
                \begin{minipage}{2.60cm}
                \centering
                    \includegraphics[height = 1.75cm]{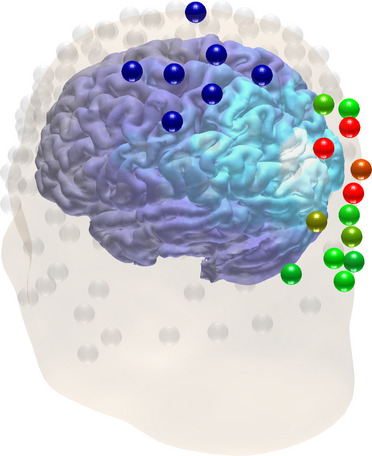} \\
                    \includegraphics[width = 2.30cm,height=1.2cm]{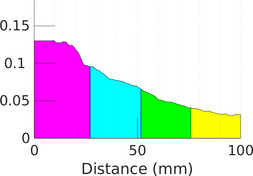} 
                \end{minipage} 
                \begin{minipage}{2.60cm}
                \centering 
                    \includegraphics[height = 2.55cm]{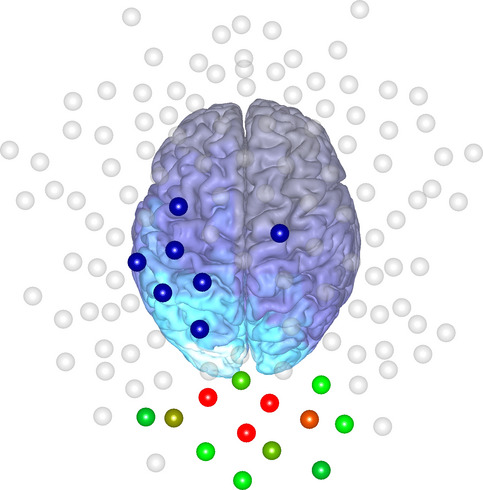}
                \end{minipage}
                \label{fig:occipital_TLS_2_20_results}
            \end{subfigure}
        \end{minipage}
    \end{minipage}
    \vskip0.1cm \hrule \vskip0.1cm
    \begin{minipage}{18.0cm}
    \centering
        \hskip0.75cm
        \includegraphics[height = 0.65cm, width = 3.5cm]{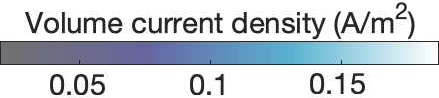}
        \hskip4.25cm
        \includegraphics[height = 0.65cm, width = 3.5cm]{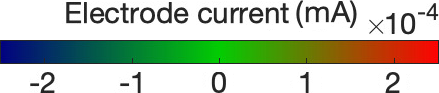}
    \end{minipage}
\end{scriptsize}
\caption{The 20-electrode montage and current pattern (mA) together with corresponding volume current density A/m\textsuperscript{2} obtained through two consecutive runs of the two-stage metaheuristic optimization process with respect to non-fixed and fixed montage of 20 active electrodes, respectively. Average magnitude in the direction of the target dipole is shown as a function of distance (mm) from the dipole position. The colorbar of the current pattern shows a color gradient for the interval from -0.25 to 0.25 mA to enhance the visibility of small variations in the pattern.}
\label{fig:results_20_channel}.
\end{figure*}

Based on the results, the dependence of the optimization accuracy and deviation on the spatial position and orientation of the target dipole becomes evident while the current patterns obtained via each applied method maintain their general characteristics regardless of the positioning of the target dipole. L1L1 and L1L2 tend to find a pattern where a large part of the stimulus current is driven through a comparably few electrodes in the current pattern as compared to TLS where the current amplitudes have smoother transitions between the electrodes. Consequently, L1L1- and L1L2-optimized current patterns are also likely to have a greater maximum current $\| {\bf y} \|_\infty$. Moreover, L1L1 and L1L2 solutions tend to include relatively many low-amplitude currents with close-to-equal amplitudes, which distributes the nuisance field current density over a large area  decreasing its amplitude. This can be interpreted as a consequence of the relatively large maximum current and is demonstrated by the results of the 20-electrode montage in Figure \ref{fig:results_20_channel}, while it is a somewhat less pronounced phenomenon with 8-electrode montages.

In L1L1 and L1L2, the anodal and cathodal electrodes tend to be further apart and the whole pattern is likely to be wider than in TLS, especially, in the pattern obtained via optimization strategy (B). Thus, such patterns can be interpreted as beneficial for enhancing the focused current density $\Gamma$ which is maximized in (B). Finally, L1L1 seems to find the most focal current density with a comparably large threshold for $\Gamma$, allowing for finding a focal stimulus with a relatively high focal current density. This is reflected by the results obtained via optimization strategy (A), where the current ratio $\Theta$ was greater for L1L1 than for L1L2 or TLS. The improved focality obtained with L1L1 as compared to TLS was observed with 20-electrode montage for somatosensory and auditory dipole, where the L1L2 and TLS solutions for cases (A) and  (B) are mutually similar  and less focal than the L1L1 solution for (A). For the 8-electrode montage the difference between L1L1 (A) and TLS (A) is less distinct, while L1L1 (A) yields a greater $\Theta$ with a larger marginal than the estimated maximal deviation regardless the electrode count.

\section{Discussion}
\label{sec:Discussion}
This study considered L1-norm data fitting via L1-norm regularized convex optimization (L1L1) as a potential alternative optimization approach for finding stimulus currents for a multi-channel tES exercise in comparison with both L2-norm data fitting counterpart (L1L2) \cite{wagner2016optimization}, and the weighted Tikhonov regularized least squares method (TLS) \cite{dmochowski2017optimal}. L1L1 has been earlier suggested as a tool to maximize the focused current for a given location and a given set of electrodes \cite{dmochowski2011optimized,fernandez2020unification}. We proposed the means to optimize global current density fitting to obtain the best possible localization of the stimulus current distribution and pattern.

The present topic is important due to the general tendency of L1-norm fitted solutions to be sparse compared to those obtained via L2-norm methods \cite{dmochowski2011optimized, dmochowski2017optimal, wagner2016optimization, khanindividually}. In tES, this means a greater focused current density driven in the targeted brain region. We applied L1-norm in both fitting and regularization, i.e., penalization of the non-zero entries in the current pattern, hypothesizing that the resulting volume current distribution and the current pattern are sparse. We considered this approach necessary to obtain the best possible fit for a given, user-defined number of active electrodes, from a montage with a minimal set of two channels, i.e., the standard two-patch tDCS \cite{Nitsche_Paulus2000, nitsche2003}, to a significantly higher version \cite{fernandez2020unification}. We apply our study into montages of 8 and 20 active electrodes based on commercial clinically-applied tES systems \cite{roy2019integration,tost2021choosing}. We defined an explicit parameter to steer the weight of the nuisance field distribution in L1L1 and L1L2, the reason being that with zero-weighting the maximum current in the targeted location is lower than with an appropriately chosen weight. In TLS, the weight is incorporated as the multiplier of the nuisance field term as shown in \ref{effect_of_weighting}.

The results suggest that the present L1L1 optimization technique performs appropriately for different target dipole locations and active electrode counts. L1L1 allowed finding a focal current pattern and well-localized stimulus current density in the brain with a focused maximum greater than 0.11 A/m\textsuperscript{2} which can be considered as  adequate in tES with the maximum current dose 2 mA. Compared to L1L2 and TLS, an enhanced focused vs.\ nuisance field current ratio $\Theta$ was obtained. This result is in accordance with our initial hypothesis on the potential benefit of L1L1 in localizing the stimulus current and the general knowledge that L1-norm optimization is advantageous to enhance contrasts between different solution components \cite{kaipio2006statistical,bertero2020introduction}. The comparison between L1L1 and L1L2 was particularly important to enlighten the enhanced potential of L1-norm fitting in maximizing the ratio $\Theta$, i.e., in the suppression of the nuisance currents.

The search case (B) suggest that the present L1L1 yields a greater focused current density {\em per se} compared to TLS, which is in agreement with the earlier observations \cite{dmochowski2011optimized}. The maximum obtained with L1L1 being was also systematically greater than that of L1L2, however, by a small margin. Notably, in this study, the greatest focused current amplitude is found with a non-zero nuisance field weight, which highlights the feasibility and importance of the present convex optimization approach, where both the current density field anywhere in the brain can contribute the solution of the optimization process: fitting the focused field alone will not yield the best possible optimization outcome.

Akin to the volume fields in the brain, the current patterns found using L1L1 were more concentrated and had, overall, a greater contrast than those obtained using L1L2 or TLS. In L1L1 and L1L2, the anodal and cathodal electrodes have greater current amplitudes and tend to be further apart from each other and less clustered than in TLS. In particular, L1L1 was shown to find a current pattern with a steep contrast between the anodal and cathodal electrodes while suppressing the nuisance currents in the brain, hence, providing a potential alternative to modulate the effects of the stimulation, e.g., the sensation experienced by the subject. These observations might be significant regarding physiological impediments, hardware constrains, or concentrations of skin irritability \cite{FITZ201461,PANERI2016740, Furubayashi2008} that the subject experiences during the stimulation session. Thus, an L1L1-based electrode montage and current pattern might  provide a potential alternative, if there is a need to modulate the stimulation configuration due to various effects it may cause.

The optimization parameter ranges covered in the metaheuristic optimization process included the neighborhoods of the global maximizers for $\Gamma$ and $\Theta$, while  the $36 \times 36$ search lattice allowed finding the full set of candidate optimizers in a relatively short computing time which was 7390, 11134, and 138 seconds with L1L1, L1L2, and TLS, respectively. Notably, our L1L1 implementation ran faster than L1L2, obviously, as L1L1 does not include additional quadratic constraints which are necessary in a semidefinite formulation of a linearly constrained problem \cite{tutuncu2003solving}. To improve the computational performance of the metaheuristic L1L1 or L1L2 process, one can consider computing multiple candidate solutions simultaneously, as the CVX optimizer is a single-thread process and allows for a straightforward parallelization in a multi-core processing environment. In addition to being a simpler method, TLS is automatically parallelized by Matlab's interpreter, as it includes only full matrices algebraic simple linear algebraic operations, which in part explains the computing time differences to L1L1 and L1L2 optimization.

Our estimates for the maximum lattice-based deviation of the optimized quantities suggest that the current computing accuracy is high enough to demonstrate the major differences between the L1L1, L1L2 and TLS method and to verify our initial hypotheses on the performance of L1L1. Evidently, the applied linear programming algorithm itself \cite{Boyd2004} might also affect the optimization outcome. The primal-dual interior-point algorithm \cite{tutuncu2003solving} of the SDPT3 solver was found to perform robustly in this study but, for example, MATLAB's built-in dual-simplex algorithm, which  was also tested in solving the L1L1 task, was found to result in a less robust outcome.

Potential future work would be to compare the present metaheuristic and CVX/SDPT3-based L1L1 and L1L2 implementations with an alternative optimization algorithm, most prominently, the Alternating Direction Method of Multipliers (ADMM) \cite{wagner2016optimization, khanindividually}. Our metaheuristic optimizers, which are available in ZI \cite{he2019zeffiro}, can also potentially be extended to include other non-invasive and invasive brain stimulation modalities to enhance the present electrical stimulation toolbox, for instance, Deep Brain Stimulation (DBS) applications \cite{jarvenpaa2020optimization}. Further investigations of the relationship between the explicit nuisance field weight $\varepsilon$ (tolerance) and modelling or other uncertainties is an important topic, as the uncertainty may be expected to limit the maximal obtainable accuracy for the nuisance field. Based on the present results, the current ratio $\Theta$ is optimized with a considerably lower value of $\varepsilon$ as compared to the focused current density $\Gamma$. Thus, the maximal obtainable $\Theta$ can be expected to decrease along with an increasing uncertainty, while the maximal $\Gamma$ can be assumed to be less affected by that. Finally, experimental work will obviously be needed to learn about other than the mathematical or computational aspects of L1L1-optimized current patterns in practice.

\section{Acknowledgments}
\label{sec:ack}
FGP, AR, MS, and SP were supported by the Academy of Finland Center of Excellence in Inverse Modelling and Imaging 2018--2025, DAAD project (334465) and by the ERA PerMed project PerEpi (344712). AR was supported by the Alfred Kordelini Foundation and MR was supported by PerEpi.

\subsection{Conflict of Interest}
\label{sec:conflict}
The authors certify that this study is a result of purely academic, open, and independent research. They have no affiliations with or involvement in any organization or entity with financial interest, or non-financial interest such as personal or professional relationships, affiliations, knowledge or beliefs, in the subject matter or materials discussed in this manuscript.

\appendix
\label{sec:Appendix}

\section{Forward Model}
\label{sec:BConditions}
The governing partial differential equation for the electric potential in the head model $\Omega$ is of the form
\begin{equation}
\label{eq:Max}
    \nabla \cdot (\sigma \nabla u) = 0.
\end{equation}

\label{sec:CEM}
The head model $\bf{\Omega}$ is stimulated through a montage of $(e_\ell)_{\ell=1}^L$  electrodes of size $\lvert e_\ell \rvert$. We denote the current applied on the ${\bf \ell}$-th electrode by $\textbf{I}_\ell$, electrode potential $\bf{U}_\ell$, and impedance $\bf{Z}_\ell$. The boundary conditions for the Complete Electrode Model (CEM) are the following:
\begin{align}
    \label{CEM1}
    0 &= \sigma \frac{\partial u}{\partial n}(\vec{r})  , &  \hbox{for }  \vec{r} \in  \partial \Omega \backslash \cup_{\ell=1}^L e_\ell, \\
    \label{CEM2}
    I_\ell &= \int_{e_\ell}{\sigma \frac{\partial u}{\partial n}(\vec{r}) dS}, & \hbox{for } \ell \!= \!1, \ldots,  L,  \\
    \label{CEM3}
    {U_\ell} &= u(\vec{r})  \! + \! \tilde{Z}_\ell \sigma \frac{\partial u}{\partial n}(\vec{r}),  \, &  \hbox{ for }  \vec{r} \in   e_\ell, \, \ell \!=\! 1,  \ldots,  L.
\end{align}
The boundary condition \ref{CEM1} describes that no current is flowing inside nor outside of head; \ref{CEM2} describes that the total current flux through the $\ell$-th electrode equals to the applied current $I_\ell$; \ref{CEM3} describes the relationship between the ungrounded electrode potential $U_\ell$ and the potential $u$ underneath the electrode;  By assuming that the effective contact impedance is $\tilde{Z}_\ell = Z_\ell \lvert e_\ell \rvert$, we can rewrite \ref{CEM3} as
\begin{equation}
    {U_\ell}  =  \frac{\int_{e_\ell}{u \, dS}}{\lvert e_\ell \rvert} + Z_\ell I_\ell.
\end{equation}

\subsection{Weak Form}
\label{app:forward_model}
A general weak form for electric potential field ${\bf u} \in {H}^1(\Omega)$ can be obtained integrating by parts. Here, ${ H}^1(\Omega)$ denotes a \textit{Sobolev space} of square integrable ($\int_\Omega |u|^2 \, dV < \infty$) functions with square integrable partial derivatives. By multiplying \ref{eq:Max} with a smooth enough test function $v \in S$, where ${S}$ is a subspace of ${H}^1(\Omega)$, it follows that
\begin{eqnarray}
\label{eq:integration}
     0 &=& -\int_\Omega \nabla \cdot (\sigma \nabla u) v \, dV, \nonumber \\
                                                &=& \int_\Omega \sigma \nabla u \cdot \nabla v \, dV - \int_{\partial \Omega} \sigma \frac{\partial u}{\partial n} v \, dS, \nonumber \\
                                                &=& \int_\Omega \sigma \nabla u \cdot \nabla v \, dV - \sum_{\ell = 1}^L \int_{e_\ell} \sigma \frac{\partial u}{\partial n} v \, dS.
\end{eqnarray}
In addition, we have the following equations:
{\setlength\arraycolsep{2 pt} \begin{eqnarray}
  - \sum_{\ell = 1}^L \int_{e_\ell} \sigma \frac{\partial u}{\partial n} v \, dS &= & - \sum_{\ell = 1}^L\,\frac{U_\ell}{Z_\ell |e_\ell|} \int_{e_\ell}  v \, dS \nonumber\\ & & +\sum_{\ell = 1}^L\,\frac{1}{Z_\ell |e_\ell|} \int_{e_\ell} u v \, dS .
\end{eqnarray}}
As a result, we may rewrite the formula $(A.6)$ as follows:
\begin{eqnarray}
\label{eq:weakform}
  0 & = & \int_\Omega \sigma \nabla u \cdot \nabla v \, dV - \sum_{\ell = 1}^L \frac{I_\ell}{| e_\ell |} \int_{e_\ell} v \, dS \nonumber \\ 
    & & - \sum_{\ell = 1}^L \frac{1}{{Z}_\ell | e_\ell |^2} { \int_{e_\ell} u \, dS \int_{e_\ell} v \, dS}  \nonumber \\ 
    & & + \sum_{\ell = 1}^L \frac{1}{{Z}_\ell | e_\ell |}  \int_{e_\ell} u v \, dS,
\end{eqnarray}
for all $v \in S$. The left-side of \ref{eq:weakform} defines a \textit{diffusion operator}. On the right-side, the first term corresponds to \textit{neural activity}, the second term to the \textit{targeted stimulus}, the third and fourth terms describe the \textit{shunting effects}.

\subsection{Resistivity Matrix}
\label{app:res_matrix}
Given the scalar valued functions $\psi_1, \psi_2, \ldots, \psi_{N} \in \mathcal{S}$, the potential ${\bf u}$ can be approximated as the finite sum $\mathbf{u} = \sum_{i = 1}^{N} z_i \psi_i$. Denoting by ${\bf z} = (z_1, z_2,\ldots, z_{N})$ the coordinate vector of the discretized potential, by ${\bf w} = (w_1, w_2,\ldots, w_{L})$ the (ungrounded) electrode voltages, and by ${\bf y} = (y_1, y_2, \ldots, y_L)$ as the injected current pattern, the weak form \ref{eq:weakform} is given by 
\begin{equation}
\label{eq:u_system}
     \left( \begin{array}{cc} {\bf A} & -{\bf B} \\
    -{\bf B}^T & {\bf C}
    \end{array} \right) \left( \begin{array}[c]{@{}l@{}} {\bf z}  \\
    {{\bf w}}
    \end{array} \right) = \left( \begin{array}[c]{@{}l@{}} {\bf 0}  \\
    {\bf y}
    \end{array} \right).
\end{equation}
Here, ${\bf A}$ is of the form
\begin{equation}
    a_{i, j} = \int_{\Omega} \sigma \nabla \psi_i \cdot \nabla \psi_j \, d V + \sum_{\ell = 1}^L \frac{1}{Z_\ell |e_\ell|} \int_{e_\ell} \psi_i \psi_j \, dS, 
\end{equation}
and the entries of ${\bf B}$ ($N$-by-$L$) and ${\bf C}$ ($L$-by-$L$) are given by 
\begin{eqnarray} 
\label{eq:fem_system} 
    b_{i, \ell}  & = & \frac{1}{Z_\ell |{e}_\ell|} \int_{e_{\ell}} \psi_i \, dS  , \label{hupu} \\
    c_{\ell, \ell}  & = & \frac{1}{Z_\ell}\,.
\end{eqnarray}
Consequently, the resistivity matrix ${\bf R}$ satisfying ${\bf z} = {\bf R} {\bf y}$ can be expressed as
\begin{equation}
    \label{eq:resistivity}
    {\bf R}  = {\bf A}^{-1} {\bf B} ( {\bf C}  - {\bf B}^{T} {\bf A }^{-1} {\bf B} )^{-1} .
\end{equation} 
The ungrounded electrode voltages ${\bf w}$ can be obtained by referring to the bottom row of \ref{eq:u_system}, i.e., ${\bf y} = {\bf-B^T z +Cw}$.

\subsection{Lead Field Matrix}
\label{app:Lead_Field}
By $\textbf{F}^{(k)}$ we denote a matrix which evaluates the $k$-th Cartesian component of the volume current density $- \sigma \nabla u$ when multiplied by the coordinate vector $\textbf{z}$ of the discretized electrical potential distribution $\textbf{u}$. The entries of this matrix are given by
\begin{equation}
 f_{i,\ell}^{(k)}=
 \begin{cases}
    -( \sigma_{ij} (\nabla \psi_\ell)_j)^k, & \hbox{if} \quad \hbox{supp}\{\psi_\ell\} \cap  T_{i} \neq \emptyset \\
    0 , & \text{otherwise},
    \end{cases}
\end{equation}
for $i, j, \ell=1, \ldots, N$ where subsets $T_i$, $i = 1, 2, \ldots, N$ form a partitioning of $\Omega$ for a user-defined dimension $N$. The $k$-th Cartesian component of the discretized volume current distribution given the stimulating current pattern can be obtained as follows ${\bf F}^{(k)}{\bf R}{\bf y}$, where
${\bf F}^{(k)}=(f_{i,\ell}^{(k)})^T, k=1, 2, 3$. Further, we define lead field matrix ${\bf L}$ as
\begin{equation}
\label{eq:L_Matrix}
    {\bf L} = \begin{pmatrix} {\bf F}^{(1)} \\ {\bf F}^{(2)} \\  {\bf F}^{(3)} \end{pmatrix} {\bf R} ={\bf F}{\bf R},
\end{equation}
where ${\bf F}=({\bf F}^{(1)},{\bf F}^{(2)}, {\bf F }^{(3)})^T$ and ${\bf L}=({\bf L}^{(1)},{\bf L}^{(2)}, {\bf L}^{(3)})^T$ with components ${\bf L}^{(k)}=({\bf L}_1^{(k)}, {\bf L}_2^{(k)})^T, k=1, 2, 3$. Formula (\ref{eq:L_Matrix}) can be considered as the \textit{forward mapping} in the process of optimizing the current pattern. 

\section{The effect of weighting in TLS}
\label{effect_of_weighting}
Denoting ${\bf W} = ( {\bf L}_1^T {\bf L}_1  \! + \!  \alpha^2 \sigma^2 {\bf I} )^{-1}$ and by $\tilde{\bf y} = {\bf W}  {\bf L}_1^T {\bf x}_1$ the special solution of (\ref{TLS_formula}) with $\delta = 0$, the general solution of (\ref{TLS_formula}) can be written as
{\setlength\arraycolsep{2 pt} \begin{eqnarray}
  {\bf y} &  = &   ( {\bf I}  + \delta^2 \alpha^2 {\bf W} {\bf L}_2^T {\bf L}_2 )^{-1} \tilde{{\bf y}} =  \tilde{\bf y} -  \delta^2 \alpha^2  {\bf W} {\bf L}_2^T {\bf L}_2  \tilde{\bf y} \nonumber\\ & & +  O(\delta^4), 
 \end{eqnarray}}
following  from the geometric series formula
{\setlength\arraycolsep{2 pt}  \begin{eqnarray}
({\bf I} \! + \!  \delta^2 \alpha^2 {\bf W} {\bf L}_2^T {\bf L}_2)^{-1} & = & {\bf I} \! - \! \delta^2 \alpha^2 {\bf W} {\bf L}_2^T {\bf L}_2 \nonumber\\ & & \! + \!  \delta^4 \alpha^4 {\bf W}^2 ( {\bf L}_2^T {\bf L}_2)^2 \! - \!  \cdots
\end{eqnarray}}
with assumption that $\delta < 1$, so that $\delta^2 \alpha^2 \| {\bf W} {\bf L}_2^T {\bf L}_2\|_2^2 \leq 1$. When $\delta = 0$, the special solution $\tilde{\bf y}$ gives the most intense current distribution to the direction of ${\bf x}$ at the targeted location, while the ratio between the focused field and nuisance field increases along with the value of $\delta \geq 0$. Namely, the inner product between the focused field and the targeted stimulus is of the form
 {\setlength\arraycolsep{2 pt} \begin{eqnarray}
{\bf x}_1^T {\bf L}_1 {\bf y}  & = &  {\bf x}_1^T {\bf L}_1 \tilde{\bf y} -   \delta^2 \alpha^2 {\bf x}_1^T  {\bf L}_1 {\bf W}{\bf L}_2^T    {\bf L}_2 \tilde{\bf y} +  O(\delta^4) \nonumber \\ & = &  {\bf x}_1^T {\bf L}_1 \tilde{\bf y} -   \delta^2 \alpha^2 \tilde{\bf y}^T {\bf L}_2^T    {\bf L}_2 \tilde{\bf y} +  O(\delta^4)
 \end{eqnarray}}
which can be further written as the following ratio between $\Gamma$ defined in (\ref{gamma_def}) and $\tilde{\Gamma} = \frac{{\bf x}_1^T {\bf L}_1 {\bf \tilde{y}}}{\| {\bf x}_1 \|_2}$, respectively:
  {\setlength\arraycolsep{2 pt} \begin{eqnarray} 
  \label{focused_field}
 \frac{\Gamma}{\tilde{\Gamma}} & = & \frac{{\bf x}_1^T {\bf L}_1 {\bf y}}{ {\bf x}_1^T {\bf L}_1 \tilde{\bf y}} = 1 - \delta^2 \alpha^2  \frac{  \| {\bf L}_2 \tilde{\bf y} \|_2^2}{{\bf x}_1^T {\bf L}_1 \tilde{\bf y}}  +  O(\delta^4).
 \end{eqnarray}}
The squared norm of the nuisance field can be written as ${\bf L}_2 {\bf y} = {\bf L}_2 \tilde{\bf y} \! - \!  \delta^2 \alpha^2 {\bf L}_2 {\bf W} {\bf L}_2^T    {\bf L}_2 \tilde{\bf y}+  O(\delta^4)
$, which yields the ratio
 {\setlength\arraycolsep{2 pt} \begin{eqnarray}
   \label{squared_ratio}
 \frac{\| {\bf L}_2 {\bf y} \|^2_2}{\| {\bf L}_2 \tilde{\bf y} \|^2_2}  & = &  1   - 2\,\delta^2 \alpha^2  \frac{ \| {\bf L}_2^T   {\bf L}_2 \tilde{\bf y} \|^2_{W}}{\| {\bf L}_2 \tilde{\bf y}\|^2_2 }+  O(\delta^4), 
 \end{eqnarray}}
where for convenience we have used the following norm definition $\| {\bf z} \|^2_W := {\bf z}^T {\bf W} {\bf z}$\,. The square root of (\ref{squared_ratio}) is of the form
 {\setlength\arraycolsep{2 pt} \begin{eqnarray}
   \label{offfield}
 \frac{\| {\bf L}_2 {\bf y} \|_2}{\| {\bf L}_2 \tilde{\bf y} \|_2}  & = &  1   -  \delta \alpha  \frac{ \| {\bf L}_2^T   {\bf L}_2 \tilde{\bf y} \|_{W}}{\| {\bf L}_2 \tilde{\bf y}\|_2 }+  O(\delta^2),
 \end{eqnarray}}
following from the Maclaurin series of the function $h(\tau) = (1+\tau)^{1/2}$. Formulas (\ref{focused_field}) and (\ref{offfield}) show that as $0 < \delta <1$ increases, the focused field intensity decreases linearly, i.e., with a slower rate than the quadratically decreasing nuisance field norm. Hence the ratio $\Theta$ defined in (\ref{theta_def}) increases. The validity of the total dose and maximum current constraint is taken care of by scaling the solution after the minimization process in the respective order.

\bibliographystyle{elsarticle-num}
\bibliography{bibliography_1.bib}

\end{document}